\pdfoutput=1
\RequirePackage{ifpdf}
\ifpdf % We are running pdfTeX in pdf mode
\documentclass[pdftex]{sigma}
\else
\documentclass{sigma}
\fi

\numberwithin{equation}{section}

\newtheorem{Theorem}{Theorem}[section]
\newtheorem{Corollary}[Theorem]{Corollary}
\newtheorem{Conjecture}[Theorem]{Conjecture}
\newtheorem{Lemma}[Theorem]{Lemma}
\newtheorem{Proposition}[Theorem]{Proposition}
 { \theoremstyle{definition}
\newtheorem{Remark}[Theorem]{Remark} }

\begin{document}
\allowdisplaybreaks

\newcommand{\arXivNumber}{1809.00122}

\renewcommand{\PaperNumber}{046}

\FirstPageHeading

\ShortArticleName{The Degenerate Third Painlev\'e Equation}

\ArticleName{Meromorphic Solution of the Degenerate Third\\ Painlev\'e Equation Vanishing at the Origin}

\Author{Alexander V.~KITAEV}

\AuthorNameForHeading{A.V.~Kitaev}

\Address{Steklov Mathematical Institute, Fontanka 27, St.~Petersburg, 191023, Russia}
\Email{\href{mailto:kitaev@pdmi.ras.ru}{kitaev@pdmi.ras.ru}}

\ArticleDates{Received November 13, 2018, in final form May 30, 2019; Published online June 18, 2019}

\Abstract{We prove that there exists the unique odd meromorphic solution of dP3, $u(\tau)$ such that $u(0)=0$,
and study some of its properties, mainly: the coefficients of its Taylor expansion at the origin and
asymptotic behaviour as $\tau\to+\infty$.}


\Classification{34M40; 33E17; 34M50; 34M55; 34M60}

\tableofcontents

\section{Introduction}\label{sec:Introduction}
The degenerate third Painlev\'{e} equation can be written in the following form,
\begin{gather}\label{eq:dp3}
u^{\prime \prime}(\tau) =\frac{(u^{\prime}(\tau))^{2}}{u(\tau)} - \frac{u^{\prime}(\tau)}{\tau} + \frac{1}{\tau}\big({-}8 \epsilon u^{2}(\tau) + 2ab\big) + \frac{b^{2}}{u(\tau)},
\end{gather}
where $\epsilon,a,b \in \mathbb{C}$. We recall that in the case $\epsilon b=0$ equation~\eqref{eq:dp3} can be integrated in elementary functions,
Otherwise, $\epsilon b\neq0$, both parameters, $\epsilon$ and~$b$, can be fixed arbitrarily in $\mathbb C{\setminus}0$
by rescaling variables $u$ and $\tau$, so that equation~\eqref{eq:dp3} possess only one essential parameter $a$.
In our previous works with A.H.~Vartanian (see~\cite{KV2004} and~\cite{KV2010}) we used both parameters
$\epsilon$ and $b$ because asymptotic results crucially depend on $\arg\tau$, so that it is convenient to identify
suitable case of equation~\eqref{eq:dp3} without making the rescaling, because in many cases it changes $\arg\tau$.
Here, in Sections~\ref{sec:Taylor}--\ref{sec:polynom-conjectures} I discuss explicit results and it is more
convenient, to simplify corresponding expressions, to put
\begin{gather}\label{eq:epsilon}
\epsilon=+1.
\end{gather}
Moreover, in these sections we keep ``normalization'' condition $b=a$. We turn back to equation~\eqref{eq:dp3}
(with $b\in\mathbb R{\setminus}0$) and complex $a$) in the Sections~\ref{sec:monodromy} and \ref{sec:asymptotics},
where we discuss isomonodromy deformations and large $\tau$ asymptotics. Even in these sections we keep
condition~\eqref{eq:epsilon}, because the solution we study has additional symmetries~\eqref{eq:u-symmetries},
therefore for this solution it is easy to change $\epsilon=+1\to\epsilon=-1$ in asymptotic results:
\begin{gather*}
u_{\epsilon=-1}(\tau)=iu_{\epsilon=+1}(i\tau).
\end{gather*}

Equation~\eqref{eq:dp3} has two singular points: the regular one at $\tau=0$ and irregular at $\tau=\infty$.
In this paper we consider a particular solution which is holomorphic and vanishing at $\tau=0$. In other words,
we study a meromorphic solution of equation~\eqref{eq:dp3} with the additional condition $u(0)=0$.
In Section~\ref{sec:Taylor} we prove that for all $a\notin i\mathbb Z/2$ this solution exists and is unique;
In the case $a-i/2\in i\mathbb Z$ this solution exists and is unique under the additional assumption that it is
an odd function of $\tau$.
It is this solution we denote throughout the paper $u(\tau)$, i.e., the last notation does not mean in this
paper the general or any other solution of equation~\eqref{eq:dp3}. The only exception is
Section~\ref{sec:monodromy}, where we recall a relation of equation~\eqref{eq:dp3} with the theory
of isomonodromy deformations and where the usage of $u(\tau)$ in a more general sense is specially stated.

My attention to this solution was attracted by B.I.~Suleimanov. In \cite{S} he studied some asymptotics
of a function which can be identified as a special solution of the second member of the hierarchy of the
third Painlev\'e equation (${}_2P_3$-function), which generally depends on two variables. When one of these
variables vanishes then ${}_2P_3$-equation reduces to the third Painlev\'e equations ($P_3$): one of them is
equivalent to the well-known similarity reduction of the $\sin$-Gordon equation, and the other one to the
solution $u(\tau)$ for the following values of the coefficients:
\begin{gather*}
\epsilon=\pm1,\qquad
a=\frac{i}2,\qquad
b=64/k^3,
\end{gather*}
where $k>0$ is a parameter from~\cite{S}. In fact, Suleimanov derived a different form of equation~\eqref{eq:dp3},
see equation~(25) in~\cite{S}, which differs by a change of the term $u^2(\tau)$ in equation~\eqref{eq:dp3}
by $u^2(\tau)/\tau$, so that at first glance one may think that these two equations are not equivalent.
However, R.~Garnier in the footnote on p.~52 of his paper~\cite{G} explained that there exists a simple point
transformation of variables mapping these two equations to each other. I also found convenient to use this (Garnier)
form of equation~\eqref{eq:dp3}, but in a slightly modified shape (see below equation~\eqref{eq:trans-Garnier}).
As I explain in Section~\ref{sec:Taylor} specifically for the Suleimanov case we have to add the condition that
$u(\tau)$ is an odd function, because the odd solutions of equation~\eqref{eq:dp3} corresponds to holomorphic
solutions of the Garnier equation.

Suleimanov's study has some physical motivation: self-focusing phenomenon for the one-dimensional
propagation of the light in unstable quasi gaseous media. He argues that under the assumption of small dispersion
in the case when the wave propagation can be described by the integrable nonlinear Schr\"odinger equation the
leading term of the intensity of a signal in the neighborhood of the focusing points under a number of scaling
transformations is described by a~special case of ${}_2P_3$-function. This is a new application of the higher
Painlev\'e functions, so that it would be interesting to see a more detailed explanation of these ideas.

The question important in view of the above studies concerns asymptotics as $\tau\to\infty$, speci\-fi\-cal\-ly,
the so-called connection problem for the Suleimanov solution of equation~\eqref{eq:dp3}.
Such problem with the help of the isomonodromy deformation method (IDM) for the general solution of
equation~\eqref{eq:dp3} was solved in \cite{KV2004} and \cite{KV2010}. Moreover, with the help of the
results obtained in these works an experienced reader by making certain tricks (analytic continuation,
rescaling and finding limits of the monodromy parameters, applying symmetries and transformations) can
find asymptotics of all solutions of equation~\eqref{eq:dp3} either for $\tau\to0$ and $\tau\to\infty$ with
$\arg\tau=\pi k/2$, $k\in\mathbb Z$.
The standard way to solve the asymptotic problem with the help of IDM
is the following: using the known behaviour of a solution at $\tau=0$, find the corresponding monodromy
parameters, then find the asymptotics as $\tau\to\infty$ in terms of
the monodromy parameters. The problem involved in
application of the results obtained in~\cite{KV2004} to the Suleimanov solution is that the
results are presented for leading terms of the asymptotic behavior of the general solution at $\tau=0$,
where equation~\eqref{eq:dp3} has the regular singularity, while the Suleimanov solution at this point is regular
and moreover vanishes, so it is a very special ``degenerate'' solution. One cannot obtain the complete set of
the monodromy data for this solution by just comparing its behaviour at the origin with the corresponding
behaviour of the general solutions. So, to get the complete set of the monodromy data we have to apply one of
the tricks mentioned above. Since such tricks might be too involved for the reader who just need to use the result
I dedicated Section~\ref{sec:monodromy} for the detailed derivation of the monodromy data.
In Section~\ref{sec:asymptotics} I use the monodromy data obtained in Section~\ref{sec:monodromy} and
the results given in papers~\cite{KV2004} and~\cite{KV2010} to get asymptotics of $u(\tau)$ for the large values
of $\tau$. I also present there a few plots of $u(\tau)$ and its large asymptotics for the pure imaginary and
real negative values of the parameter $a$. While writing Sections~\ref{sec:monodromy} and~\ref{sec:asymptotics}
I~noticed and corrected a few (nondramatic) faults in our works~\cite{KV2004} and~\cite{KV2010}, this information
will be useful for the readers interested in the results concerning any other solutions of equation~\eqref{eq:dp3}.

Looking at the plots of $u(\tau)$ presented in Section~\ref{sec:asymptotics} one can make a reasonable hypothesis
concerning the behaviour of $u(\tau)$ for positive finite values of $\tau$, these issues are discussed in
Section~\ref{sec:positive}.

The major part of this paper concerns the study of the coefficients of the Taylor expansion of the
solution holomorphic at $\tau=0$ for general value of the parameter $a$. These coefficients are rational
functions of $a^2$ with some interesting properties, which I was not able to prove directly with the help of
the recurrence relation. I~formulated corresponding conjectures in
Sections~\ref{sec:Taylor} and~\ref{sec:polynom-conjectures}.

Trying to prove these properties in Sections~\ref{sec:generating-A}--\ref{sec:residues} I develop a technique
of generating functions for various parameters defining the Taylor coefficients: I introduce
specially rescaled solutions of equation~\eqref{eq:trans-Garnier}, which I call super generating functions (SGF).
The scaling parameter is a proper function of $a^2$, say, in the simplest case it is just $a^2$. Then we develop
SGF into the Laurent series with respect to the scaling parameter, the coefficients of these expansions appear
to be functions of the rescaled variable $\tau$, which we denote as $z$, the latter functions of $z$ appear to be
generating functions for Laurent expansions in the scaling parameter of the original Taylor coefficients of the
solution holomorphic at $\tau=0$. On this way we find some non-trivial explicit general formulae and check for
them conjectures made in Section~\ref{sec:Taylor}. However, the full proof of the conjectures presented in
Section~\ref{sec:Taylor} requires additional ideas. This technique of generating functions has quite general nature
and can be applied not only for general solutions of equation~\eqref{eq:dp3} but also in analogous cases for the
other Painlev\'e equations.

In Section~\ref{sec:polynom-conjectures} we consider divisibility properties of the polynomials defining the
numerators of the Taylor coefficients. It is written in a fully conjectural manner. Due to the appearance of
the number $3$ in the recurrence relation for the Taylor coefficients, one may suspect that the set of coefficients
of the polynomial might suffer some special divisibility properties with respect to~$3$. After extensive numerical
studies I was able to formulate a conjecture concerning the greater common divisor (g.c.d.) of the coefficients.
In particular, g.c.d.\ contains the powers of~$3$. These powers vary depending on the coefficient in an amazing
way. If we plot the dependence of the power of~$3$ in g.c.d.\ of the $n$-th coefficient as a function of $n$
and connect neighboring points (with the abscissas $n$ and $n+1$) for all natural $n$ by line segments, then we get
a plot which I call a quasiperiodic fence. The fence is built on the positive horizontal semi-axis and,
architecturally, can be divided by segments. Each segment consists of two parts: the first one resembles all the
previously built fence and the second part is a newly constructed one. The rules of the construction have periodic
properties described in the section.

In view of the paper \cite{GIL2013} reporting combinatorial properties of small-argument expansions of
the $\boldsymbol{\tau}$-functions for the Pinlev\'e equations, it would be interesting to study whether the
combinatorics would be useful for proofs of the conjectures made in this paper. One of the obstacles on this way is that the
properties of the Taylor coefficients which are discussed in this paper do not preserve under transformations of
the series when turning from the $\boldsymbol{\tau}$-function to $u(\tau)$. At the same time it would be
interesting to check whether the properties of the coefficients that we discuss in this paper or some of their
analogues take place for the corresponding expansions of the $\boldsymbol{\tau}$-functions. On the other hand,
the technique of SGF that I develop in this work might be helpful for counting of the combinatorial objects
introduced in~\cite{GIL2013}.

During preparation of this work appeared a preprint~\cite{BLM} where exactly the same solution
of the nonlinear Schr\"odinger equation as the one studied by Suleimanov manifests itself in an absolutely
different physical context, namely, in the study of the so-called rogue waves of infinite order. The main subject of this paper is
the asymptotic study via the Deift--Zhou steepest descent method of the special case of ${}_2P_3$-function mentioned
above and a discussion of the results obtained to the description of the rogue waves. The ${}_2P_3$-function in
the notation of~\cite{BLM} depends on two variables, $X$ and $T$. In some sense it can be viewed as a deformation
from the regular at the origin solution of the third Painlev\'e equation (so-called~$D_8$ case of~$P_3$) at $T=0$
to the analogous solution of the degenerate third Painlev\'e equation (also known as~$D_7$ case of~$P_3$)
with $a=i/2$ at $X=0$. Equation~(223) of \cite{BLM} presents the leading term of the large-$T$ asymptotics of function $\Psi^+(0,T)$,
which in our notation up to a scalar constant factor coincides with ${\rm e}^{i\varphi(\tau)}$, $\tau^2=T$ for
$a=i/2$ and a proper value of parameter~$b$ (see Section~\ref{sec:monodromy}). We do not discuss asymptotics
of $\varphi$ here. More information about asymptotics of this function for general values of the parame\-ters~$a$ and~$b$ can be found in~\cite{KVdP3int}.

\section[Expansion as $\tau\to0$]{Expansion as $\boldsymbol{\tau\to0}$}\label{sec:Taylor}
In this section we assume that both parameters $a,b\in\mathbb{C}$. We begin with the following lemma which looks
special however, in fact, resembles the general case.
\begin{Lemma}\label{lem:odd}
If there exists a holomorphic solution of equation~\eqref{eq:dp3} for $a=b$, then for $4a^2+(2k+1)^2\neq0$,
$k\in\mathbb Z$, it is an odd function of $\tau$. If $4a^2+(2k+1)^2=0$, then this solution depends on one complex parameter. This parameter
can be chosen such that this solution is odd. Moreover, for $a^2+k^2\neq0$ the odd solution is unique.
\end{Lemma}
\begin{proof}The reader may start with the general form of the Taylor expansion vanishing at $\tau=0$, substitute it into
equation~\eqref{eq:dp3}, and, equating the coefficients for the leading terms in $\tau$, prove that it can
be presented in the following form
\begin{gather}\label{eq:tau-Taylor}
u(\tau)=-\frac{\tau}2\left(1+\sum\limits_{k=1}^{\infty}b_k\tau^{2k}+\sum\limits_{k=0}^{\infty}c_k\tau^{2k+1}\right),
\end{gather}
where the coefficients $b_k$ and $c_k$ depend only on~$a$. Substantially, here we have determined the leading coefficient of the expansion, follows from the cancelation of
the last two terms in equation~\eqref{eq:dp3} ($a=b$!), and separated odd and even parts of the solution.
Now, substituting expansion~\eqref{eq:tau-Taylor} into equation~\eqref{eq:dp3}, taking the common denominator,
and equating zero the coefficients of powers $\tau^{2k}$, and~$\tau^{2k+1}$
we find the recurrence equations for the coefficients:
\begin{gather}
 c_0\big(4a^2+1\big)=0,\qquad b_1\big(a^2+1\big)=1,\label{eq:c0-b0}\\
\big(4a^2+(2k+1)^2\big)c_k =C_k(c_0,\ldots,c_{k-1},b_1,\ldots,b_{k-1}),\label{eq:ck}\\
\big(4a^2+(2k)^2\big)b_k =B_k(c_0,\ldots,c_{k-1},b_1,\ldots,b_{k-1}),\label{eq:bk}
\end{gather}
where $C_k$ and $B_k$ are polynomials in variables indicated above. Polynomials, $C_k$ have one important property,
which do no have polynomials $B_k$, namely, each term of these polynomials consists of odd number of factors
$c_l$ and some factors $b_l$, $l=0,\ldots,k-1$, where we put $b_0=1$. This observation follows from the fact that
equation~\eqref{eq:ck} is obtained by equating zero coefficients of even powers of $\tau$ and equation~\eqref{eq:bk}
of the odd, respectively. The coefficients of the even powers of $\tau$ are obtained with the help of the odd ones
by multiplication by $\tau$ which comes from the denominator. As a result such coefficients always have odd
number of factors of parame\-ters~$c_l$ with junior subscripts $l$.

Now it is easy to finish the proof. Assume $4a^2+(2k+1)^2\neq0$ for all integer $k$, then the first
equation~\eqref{eq:c0-b0} and equations~\eqref{eq:ck},~\eqref{eq:bk} imply that all coefficients $c_k$ vanish.
Thus the solution is odd. If for some integer $k=k_0$ equation $4a^2+(2k_0+1)^2=0$ is valid, then, obviously,
$c_l=0$ for $l<k_0$, and $c_{k_0}$ is a complex parameter. The coefficients $c_l$ with $l>k_0$ will be uniquely
determined. In case we put $c_{k_0}=0$ we obtain the odd solution. Finally, if for the odd solution
$a^2+k^2\neq0$ for all integers $k$ the construction is unique, because of equations~\eqref{eq:bk} and~\eqref{eq:ck}.
\end{proof}
\begin{Remark}\label{rem:a-imaginary-integer}We assumed that the holomorphic solution at $\tau=0$ exists. The proof given above does not say what
happens when $a^2+k_0^2=0$, for some integer $k_0$. By the examining cases $k_0=0$, $k_0=1$, etc one can find that
the solution does not exist for these values of $k_0$, however the ``magic'' cancelation of the r.h.s.\ of
equation~\eqref{eq:bk} cannot be excluded. In fact, it follows from the more explicit version of
equation~\eqref{eq:bk} considered below for the general case (all $c_k=0$) equations~\eqref{eq:u-2} and~\eqref{eq:u-2n-recurrence}, there are some magic cancelations at $a^2+k^2=0$ for integer~$k$, see
Conjecture~\ref{con:strange-divisibility} and the following examples. At the same time in Section~\ref{sec:residues}
we consider generating function, $V_k(z)$ for the residues of the coefficients at possible poles at $a^2+k^2=0$.
This function satisfies equation~\eqref{eq:Vk} and it does not admit solution~$V_k(z)\equiv0$ for any~$k$.
\end{Remark}

Assuming in equation~\eqref{eq:dp3} $b=a$ we make the following change of variables
\begin{gather}\label{eq:U-x}
u(\tau)=-\frac{\tau}2\left(1+U(x)\right),\qquad x=\tau^2.
\end{gather}
Then equation~\eqref{eq:dp3} can be rewritten as follows
\begin{gather}\label{eq:Garnier}
\delta_x^2\ln(1+U)=x(1+U)-\frac{a^2U}{(1+U)^2},\qquad \delta_x=x\frac{{\rm d}}{{\rm d}x}.
\end{gather}
We can further transform equation~\eqref{eq:Garnier}
\begin{gather}\label{eq:trans-Garnier}
\big(\delta_x^2+a^2\big)U=x(1+U)^3+(\delta_xU)^2-U\delta_x^2U.
\end{gather}
Equation~\eqref{eq:Garnier} we call the Garnier form of equations~\eqref{eq:dp3} and~\eqref{eq:trans-Garnier}
the modified Garnier form. Originally, R.~Garnier considered the equation equivalent to the ones written above
but for the function $y(x)=1+U(x)$.
\begin{Lemma}\label{lem:existence&uniqueness}
For any $b=a$ and $|a|<1$ there exists a unique solution of equation~\eqref{eq:trans-Garnier} holomorphic in some
neighborhood of $\tau=0$ satisfying the condition $u(0)=0$.
\end{Lemma}
\begin{proof}
Consider the equivalent form of equation~\eqref{eq:trans-Garnier},
\begin{gather*}
\big(\delta_x^2+a^2\big)U=x(1+U)^3+2(\delta_xU)^2-\delta_x(U\delta_xU).
\end{gather*}
Let us define new variables
\begin{gather*}
\delta_xU(x)=rV(x),\qquad V(x)=xV_1(x),\qquad U(x)=xU_1(x),
\end{gather*}
where $r\in(0,1)$ is a number,
and rewrite equation~\eqref{eq:trans-Garnier} as the following system
\begin{gather*}%\label{sys:U1-V1-diff}
\frac{{\rm d}}{{\rm d}x}xU_1=rV_1,\\
\frac{{\rm d}}{{\rm d}x}xV_1=-\frac{a^2}{r}U_1+\frac{1}{r}(1+xU_1)^3+2rxV_1^2-\frac{{\rm d}}{{\rm d}x}(xU_1xV_1).
\end{gather*}
After integration one obtains
\begin{gather}\label{sys:U1-V1eq1}
U_1(x) =\frac{r}{x}\int_0^xV_1(\xi){\rm d}\xi,\\
V_1(x) =-xU_1(x)V_1(x)+\frac1{x}\int_0^x\left( -\frac{a^2}{r}U_1(\xi)+\frac1{r}(1+\xi U_1(\xi))^3+2r\xi V_1^2(\xi)\right){\rm d}\xi.\label{sys:U1-V1eq2}
\end{gather}
Define vector-function $\vec{f}=(U_1,V_1)^t$, and denote the right-hand side of system~\eqref{sys:U1-V1eq1},
\eqref{sys:U1-V1eq2} as $\vec{F}(\vec{f}\,)$, then this system can be presented in the concise form
\begin{gather*}
\vec{f}=\vec{F}(\vec{f}\,).
\end{gather*}
Our goal is to show that $\vec{F}$ is a contraction map in a proper Banach space.
To prove this we consider the space of holomorphic vector-functions, $\vec{f}=(f_1,f_2)$ in the disc centered
at the origin with the radius $\epsilon$, $D_{\epsilon}$. Supply this space with the $\sup$-norm
\begin{gather*}
\|\vec{f}\, \|= \max_{x\in D_{\epsilon}}\{|f_1(x)|,|f_2(x)|\} .
\end{gather*}
Consider the Banach space of vector functions $\vec f$ holomorphic in the disc $D_{\epsilon}$
with the $\sup$-norm. Denote $\mathcal{X}$ the closed ball in this space defined by the inequality,
$\|\vec{f} \,\|\leq C$. We assume that $\epsilon$ and $C$ are defined such that the following inequalities
are valid
\begin{gather}\label{eq:def-epsilon-C}
|a|^2<r<1,\qquad
C>\frac{1+\delta}{r-|a|^2},\qquad
\delta>0,\qquad
\epsilon<\frac{\delta}{C\big(r^2C^2/9+2r(r+1)C+3\big)},
\end{gather}
where $\delta>0$ is chosen arbitrary. Note that inequalities~\eqref{eq:def-epsilon-C} imply that $\epsilon<\big(r-|a|^2\big)/3<\big(1-|a|^2\big)/3$.

We claim that conditions~\eqref{eq:def-epsilon-C} imply that $\vec{F}(\vec{f}\,)$ maps $\mathcal{X}$
into itself. We can estimate from above the r.h.s.\ of equation~\eqref{sys:U1-V1eq2}
\begin{gather*}
\frac{|a|^2}{r}C+\frac1{r}(1+\epsilon C)^3+\epsilon(1+2r)C^2.
\end{gather*}
Now we exploit inequalities~\eqref{eq:def-epsilon-C} to check that this expression is less than~$C$. Consider the difference
\begin{gather*}\begin{split}&
C-\frac{|a|^2}{r}C-\frac1{r}(1+\epsilon C)^3-\epsilon(1+2r)C^2 >\frac{1+\delta}r-\frac1r -\frac{\epsilon}{r}\big(3C+\big(r+3\epsilon+2r^2\big)C^2+\epsilon^2C^3\big)\\
&\hphantom{C-\frac{|a|^2}{r}C-\frac1{r}(1+\epsilon C)^3-\epsilon(1+2r)C^2}{} >\frac1r(\delta-\epsilon C\big(3+2r(r+1)C+r^2C^2/9\big)>0.\end{split}
\end{gather*}
In the analogous way one finds that
\begin{gather*}
\|\vec{F}(\vec{f}\,)-\vec{F}(\vec{g})\|<r_1\|\vec{f}-\vec{g}\|,
\end{gather*}
where
\begin{gather*}
r_1=\max\left\{r,\frac{|a|^2}{r}+\frac{\epsilon}{4r}\left(6+8(\epsilon C)+3(\epsilon C)^2\right)
+\left(\frac53+r\right)\epsilon C \right\}.
\end{gather*}
Obviously, making $\epsilon$ smaller, if necessary, one can reach the condition $r_1<1$ and thus $\vec{F}(\vec{f},x)$
is the contraction.
\end{proof}

As follows from Lemma~\ref{lem:existence&uniqueness} the function $U(x)$ for $|a|<1$ can be developed as
Taylor series convergent in some neighbourhood of $x=0$:
\begin{gather}\label{eq:U-series}
U(x)=\sum\limits_{n=1}^{\infty}u_{2n}(a)x^n.
\end{gather}
Substituting it into equation~\eqref{eq:trans-Garnier} one finds the following recurrence relation for the coefficients, $u_{2n}=u_{2n}(a)$:
\begin{gather}
\big(a^2+1\big)u_{2} =1,\label{eq:u-2}\\
\big(a^2+n^2\big)u_{2n} =3u_{2(n-1)}+3\sum\limits_{j_1+j_2=n-1}u_{2j_1}u_{2j_2} +\sum\limits_{j_1+j_2+j_3=n-1}u_{2j_1}u_{2j_2}u_{2j_3}\nonumber\\
\hphantom{\big(a^2+n^2\big)u_{2n} =}{} -\sum\limits_{j_1=1}^{2j_1<n}(n-2j_1)^2u_{2j_1}u_{2(n-j_1)},\qquad n=2,3,\ldots.\label{eq:u-2n-recurrence}
\end{gather}
\begin{Remark}\label{rem:coeff3}
Appearance of coefficients $3$ in the first two terms of equation~\eqref{eq:u-2n-recurrence} has an interesting
consequence on the properties of $u_{2n}(a)$ discussed in Section~\ref{sec:polynom-conjectures}.
\end{Remark}

The first few terms of the sequence $u_{2n}$:
\begin{gather*}
u_{2}=\frac{1}{a^2+1},\qquad u_{4}=\frac{3}{\big(a^2+1\big)\big(a^2+4\big)},\\
 u_{8}=\frac{5 \big(11 a^2+36\big)}{\big(a^2+1\big)^2\big(a^2+4\big)\big(a^2+9\big)\big(a^2+16\big)},\qquad
u_{6}=\frac{6 \big(2 a^2+3\big)}{\big(a^2+1\big)^2\big(a^2+4\big)\big(a^2+9\big)},\\
u_{10}=\frac{3 \big(91 a^6+1115 a^4+4219 a^2+3600\big)}{(a^2+1)^3(a^2+4)^2\big(a^2+9\big)\big(a^2+16\big)\big(a^2+25\big)},\qquad \ldots.
\end{gather*}
\begin{Proposition}\label{prop:taylor-general}
The solution defined in Lemma~{\rm \ref{lem:existence&uniqueness}} can be uniquely continued, as a holomorphic function
of $a$, on any simply connected domain in $\mathbb C{\setminus}\{i\mathbb Z\}$.
\end{Proposition}

\begin{proof}
Clearly, equations~\eqref{eq:u-2} and \eqref{eq:u-2n-recurrence} allow one to construct uniquely coefficients $u_{2n}$ $n=1,2,\ldots$
for any complex values of $a$ provided $a\notin i\mathbb Z$. We are going to prove that the series~\eqref{eq:U-series} is convergent
for these values of $a$. Take arbitrary $L>0$, and $\epsilon_1>0$ and consider a~compact subset (a~cheese-like domain) of
the complex plane~$a$:
\begin{gather*}
\mathbb{D}_{L,\epsilon_1}:=
\big\{a\in\mathbb{C},\, \big|a^2\big|\leq L \, {\rm and} \, \big|n_1^2+a^2\big|\geq\epsilon_1, {\rm if} \, n_1\in\mathbb{Z}\, {\rm and} \, n_1^2\leq L\big\}.
\end{gather*}
We are going to prove the uniform convergence of series~\eqref{eq:U-series} in $\mathbb{D}_{L,\epsilon_1}$.
Now take $\epsilon$, $0<\epsilon<1$, and $N>0$ such that for all $|n|>N$ the following estimate is valid
\begin{gather*}
\left|\frac{n^2+a^2}{n^2}\right|>1-\epsilon.
\end{gather*}
It is easy to prove the following asymptotics
\begin{gather}
\sum\limits_{j_1+j_2=n-1}\frac1{j_1^2 j_2^2} \underset{n\to\infty}{=}\frac{\pi^2}{3 n^2}+\frac{4\ln n}{n^3}+O\left(\frac{1}{n^3}\right),\nonumber\\
\sum\limits_{j_1+j_2+j_3=n-1}\frac1{j_1^2 j_2^2 j_3^2} \underset{n\to\infty}{=}\frac{\pi^4}{12 n^2}+O\left(\frac{\ln n}{n^3}\right),\nonumber\\
\sum\limits_{j_1}^{2j_1<n}\frac{(n-2j_1)^2}{j_1^2(n-j_1)^2} \underset{n\to\infty}{=} \frac{\pi^2}{6}-\frac{2\ln n}{n}+O\left(\frac{1}{n}\right).\label{eqs:sums-asymptotics}
\end{gather}
If necessary we increase $N$ to ensure for $n>N$ the inequality
\begin{gather}\label{ineq:2choice-N}
\frac{3}{(1-\epsilon)(n-1)^2}<\frac14.
\end{gather}
Now we choose $\epsilon_2\in(0,1)$ such that
\begin{gather}\label{ineq:epsilon2}
\frac{\pi^2}{6}\frac{\epsilon_2}{(1-\epsilon)}<\frac14.
\end{gather}
Again, if necessary we increase $N$ such that for all $n>N$ both inequalities are valid
\begin{gather}\label{ineq:3choice-N}
\frac{\pi^2\epsilon_2}{(1-\epsilon)n^2}<\frac14,\qquad \frac{\pi^4}{12}\frac{\epsilon_2^2}{(1-\epsilon)n^2}<\frac14.
\end{gather}
Because of asymptotics~\eqref{eqs:sums-asymptotics} and inequalities~\eqref{ineq:epsilon2} and \eqref{ineq:3choice-N}
we can, if necessary, increase for the last time~$N$ to get for all $n>N$ the following inequality
\begin{gather}\label{ineq:sum}
3\epsilon_2\sum\limits_{j_1+j_2=n-1}\frac1{j_1^2 j_2^2}+\epsilon_2^2\sum\limits_{j_1+j_2+j_3=n-1}\frac1{j_1^2 j_2^2 j_3^2}+
\epsilon_2\sum\limits_{j_1}^{2j_1<n}\frac{(n-2j_1)^2}{j_1^2(n-j_1)^2}<\frac34.
\end{gather}
After we fixed $N$ we choose the number $C_1>1$ such that for all $n=1,\ldots,N$ $u_{2n}(a)< \epsilon_2C_1^n/n^2$.
Now, recurrence relation~\eqref{eq:u-2n-recurrence} with the help~\eqref{ineq:sum} and equation~\eqref{ineq:2choice-N}
via the mathematical induction implies that $u_{2n}(a)< \epsilon_2C_1^n/n^2$ for $n>N$. This completes the proof of
uniform convergence of series~\eqref{eq:U-series} in some neighbourhood of $x=0$ for $a$ in any compact subset of~$\mathbb{C}{\setminus} i\mathbb{Z}$.
\end{proof}

\begin{Remark}One actually can prove a sharper estimate for coefficients $u_{2n}$, say, the same scheme implies an estimate $u_{2n}(a)< \epsilon_2C_1^n/n^p$ for any $p>1$.
\end{Remark}

We can summarize the previous studies as the following theorem.
\begin{Theorem}\label{th:existence & uniqueness}
For all $a,b\in\mathbb C$, $a\neq in$, and $a\neq i(n+1/2)$ with $n\in\mathbb Z$ there exists a~unique
meromorphic solution of equation~\eqref{eq:dp3} vanishing at the origin. For $a=i(n+1)/2$ there exists a~unique odd meromorphic solution of this equation vanishing at the origin.
The Taylor expansion of this solution at $\tau=0$ reads
\begin{gather}\label{eq:zero-expansion}
u(\tau)=-\frac{b}{2a}\tau\left(1+\sum\limits_{n=1}^{\infty}u_{2n}(a)\left(\frac{b}{a}\tau^2\right)^n\right).
\end{gather}
\end{Theorem}
\begin{proof}Convergence and uniqueness of expansion~\eqref{eq:zero-expansion} follows from Proposition~\ref{prop:taylor-general}
by resca\-ling of equations~\eqref{eq:U-x} and \eqref{eq:Garnier}. The meromorphicity is the consequence of
the Painlev\'e property for equation~\eqref{eq:dp3}.
\end{proof}

\begin{Remark}\label{rem:a-imaginary-integer-final}
For $a=in$ with $n\in\mathbb Z$ meromorphic solution of equation~\eqref{eq:dp3} vanishing at the origin does not
exists. For small values of $n$ it is obvious from the explicit formulae presented after Remark~\ref{rem:coeff3}.
For general $n$ it follows from the monodromy theory (see Proposition~\ref{prop:monodromy-final}) and super
generating function considered in Section~\ref{sec:residues} (see Remark~\ref{rem:a-imaginary-integer}).
\end{Remark}

We formulate further properties of the coefficients $u_{2n}(a)$ as the following conjecture.
\begin{Conjecture}\label{con:expansion-u-main}
\begin{gather}\label{eq:u2ka}
u_{2n}(a)=\frac{P_{m(n)}\big(a^2\big)}{\prod\limits_{k=1}^{n}\left(a^2+k^2\right)^{n_k}},\\
n_k=\left[\frac{n+1}{k+1}\right],\qquad
m(n)=\sum\limits_{k=1}^{n}n_k-n,\label{eq:n-k-m(n)}
\end{gather}
notation $ [\cdot ]$ means the integer part of the enclosed number
and $P_m(x)$ is the polynomial\footnote{Irreducible over $Q[x]$ for all $m\in\mathbb{Z_+}$.} of $x$ with $m\equiv m(n)=\deg P_m(x)$ and positive integer coefficients
\begin{gather}\label{eq:polynom-P}
P_{m(n)}\big(a^2\big)=\sum\limits_{k=0}^{m(n)} p_k(n)a^{2k},\qquad p_k(n)\in\mathbb{Z_+}.
\end{gather}
\end{Conjecture}
\begin{Remark}\label{rem:m(n)}
The first members of the sequence $m(n)$ are
\begin{gather*}
m(1)=m(2)=0,\qquad m(3)=m(4)=1,\qquad m(5)=m(6)=3,\qquad m(7)=5, \\
m(8)=6, \qquad m(9)=m(10)=8, \qquad m(11)=m(12)=12, \qquad m(13)=14, \\
m(14)=16,\qquad\ldots.
\end{gather*}
\end{Remark}

\begin{Remark}\label{rem:coeff-numer-conjecture}
The fact that coefficients $p_k(n)$ are integers is obvious because of recurrence relation~\eqref{eq:u-2} and~\eqref{eq:u-2n-recurrence}. At the same time it is not that easy to prove that these integers are positive because of the minus in the recurrence relation~\eqref{eq:u-2n-recurrence}.
\end{Remark}

\begin{Remark}\label{rem:coeff-denom-conjecture}
The denominator of equation~\eqref{eq:u2ka} is easy to confirm, modulo explicit expression~\eqref{eq:n-k-m(n)}
for the numbers~$n_k$. The direct proof, based on recurrence relation~\eqref{eq:u-2}
and~\eqref{eq:u-2n-recurrence}, requires confirmation of the following conjecture.
\end{Remark}

\begin{Conjecture}\label{con:strange-divisibility}
For all natural numbers $k$ and $l$, such that: $k\geq4$, $l\geq2$, $k+2\geq3l$, and $k+2$ is divisible by~$l$,
the numerator of the rational function
\begin{gather*}
\sum\limits_{\frac{k+2}l=m_1+m_2+m_3}u_{2(lm_1-1)}u_{2(lm_2-1)}u_{2(lm_3-1)}\nonumber\\
\qquad{} -\sum\limits_{\frac{k+2}l=m_4+m_5, m_4<\frac{k+2}{2l}}\left[l(m_5-m_4)\right]^2u_{2(lm_4-1)}u_{2(lm_5-1)}%\label{eq:strange-divisibility}
\end{gather*}
is divisible by $\big(a^2+(l-1)^2\big)$. Note, that to simplify our notation here and sometimes below, we omit the dependence of the coefficients $u_{2k}(a)$ of $a$, so that $u_{2k}\equiv u_{2k}(a)$.
\end{Conjecture}

We finish this section by providing examples of Conjecture~\ref{con:strange-divisibility}:
\begin{enumerate}\itemsep=0pt
\item[1)] $k=4$, $l=2$,
\begin{gather*}
\mathrm{numerator}\big(u_2\big(u_2^2-4u_6\big)\big)=\big(a^2+1\big)\big(a^2-36\big);
\end{gather*}
\item[2)] $k=6$, $l=2$,
\begin{gather*}
\mathrm{numerator}(u_2(3u_2u_6-16u_{10}))=6\big(a^2+1\big)\big(6a^6-455a^4-4676a^2-14400\big);
\end{gather*}
\item[3)] $k=8$, $l=2$,
\begin{gather*}
\mathrm{numerator}
\big(3u_2\big(u_2u_{10}+u_6^2\big)-4u_6u_{10}-36u_2u_{14}\big)=9\big(a^2+1\big)\big(139a^{12}-16186a^{10}\\
\qquad{}-833966a^8-15895545a^6-128899248a^4-449762544a^2-533433600\big);
\end{gather*}
\item[4)] $k=7$, $l=3$,
\begin{gather*}
\mathrm{numerator}\big(u_4\big(u_4^2-9u_{10}\big)\big)=27\big(a^2+4\big)\big(a^6-226a^4-1622a^2-1800\big);
\end{gather*}
\item[5)] $k=10$, $l=3$,
\begin{gather*}
\mathrm{numerator}\big(3u_4\big(u_4u_{10}-12u_{16}\big)\big)=81\big(a^2+4\big)\big(91a^{12}-42555a^{10}-2464380a^8\\
\qquad{} -55847687a^6-573508161a^4-1948922208a^2-1828915200\big);
\end{gather*}
\item[6)]$k=13$, $l=3$,
\begin{gather*}
\mathrm{numerator}\big(3\big(u_4^2u_{16}+u_4u_{10}^2-27u_4u_{22}-3u_{10}u_{16}\big)\big)=81\big(a^2+4\big)\big(22702a^{24}\\
\qquad{} -16646090a^{22}-4061032152a^{20}-413489537329a^{18}-23690569569496a^{16}\\
\qquad{}- 816781188263163a^{14}-17400650459323535a^{12}-229588162659563852a^{10}\\
\qquad{} -1844596326528992619a^8-8649917000534607066a^6\\
\qquad{} -21696167625164762400a^4 -25866244844475840000a^2\\
\qquad{}-11292874661376000000\big);
\end{gather*}
\item[7)] $k=10$, $l=4$,
\begin{gather*}
\mathrm{numerator}\big(u_6\big(u_6^2-16u_{14}\big)\big)=72\big(a^2+9\big)\big(2a^2+3\big) \\
\qquad{}\times \big(12a^{10}-8896a^8-272369a^6-2858377a^4-8718516a^2-6350400\big);
\end{gather*}
\item[8)] $k=14$, $l=4$,
\begin{gather*}
\mathrm{numerator}\big(u_6(3u_6u_{14}-64u_{22})\big)=432\big(a^2+9\big)\big(2a^2+3\big) \big(3876a^{24}-5756420a^{22}\\
\qquad{} -1235083643a^{20}-114944445444a^{18}-6103132228087a^{16}\\
\qquad{}- 197871121155883a^{14}-3930475972840326a^{12}-47045222366439497a^{10}\\
\qquad{}- 336583346858652920a^8-1396061649915602256a^6-3170843975740838400a^4\\
\qquad{}- 3496575097981440000a^2-1434015830016000000\big).
\end{gather*}
\end{enumerate}
\begin{Remark}
Note that $\mathrm{numerator}(u_2)=1$, $\mathrm{numerator}(u_4)=3$, $\mathrm{numerator}(u_6)=6\big(2a^2+3\big)$, so that corresponding factors does not effect on the divisibility by $\big(a^2+(l-1)^2\big)$.
\end{Remark}

\section[Generating function $A(a,z)$]{Generating function $\boldsymbol{A(a,z)}$}\label{sec:generating-A}

As follows from Section~\ref{sec:Taylor} the solution defined in Lemma~\ref{lem:existence&uniqueness} has
singular points at $a\in i\mathbb{Z}$ and, surely, at $a=\infty$. In this and subsequent sections we obtain
asymptotic expansions of the solution at these points. These asymptotic expansions we call
{\it super generating functions}, because their members are generating functions for infinite sequences of numbers
related with the coefficients~$u_{2n}$. Clearly, by definition, all so-called super generating functions, as the
functions in classical understanding of the notion ``function'', just coincide, modulo a rescaling, with
the solution~\eqref{eq:U-x} defined in Lemma~\ref{lem:existence&uniqueness}. However, we consider them as the
functions of the scaling parameter, which is the local parameter in the neighbourhood of singularities with respect
to variable~$a$, rather than~$\tau$ (or~$x$), moreover, the perspective we use these functions makes it convenient
to give them a different name and assign special notation.

In this section we study the super generating function at $a=\infty$ which we denote $A(a,z)$.
Define $A(a,z)$ as the following formal expansion
\begin{gather}\label{eq:A-series}
A(a,z)=\sum\limits_{k=0}^{\infty}\frac{A_k(z)}{a^{2k}},
\qquad
x=a^2z,
\end{gather}
where the coefficients $A_k(z)$ are to be defined by substitution of $A=A(a,z)$ into ODE
\begin{gather}\label{eq:A-ode}
a^2\big(z(1+A)^3-A\big)=(1+A)\delta^2_zA-(\delta_zA)^2,
\end{gather}
where equation~\eqref{eq:A-ode} is obtained by substitution $U=A$, $x=a^2z$ into equation~\eqref{eq:Garnier}.
\begin{Proposition}\label{prop:Ak-formal}
For $k=0,1,\ldots$ there exists a unique sequence $A_k(z)$ of functions holomorphic in the disc $|z|<\frac{2^2}{3^3}$
such that series~\eqref{eq:A-series} formally solves equation~\eqref{eq:A-ode}. Moreover, all func\-tions~$A_k(z)$ are rational functions of $A_0(z)$, and $A_k(0)=0$.
\end{Proposition}
\begin{proof}Substituting series~\eqref{eq:A-series} into equation~\eqref{eq:A-ode} one finds
\begin{gather}\label{eq:A0}
z(1+A_0(z))^3=A_0(z),
\end{gather}
and the following recurrence relation for $k=0,1,\ldots$
\begin{gather}
A_{k+1}(z) =\frac{1+A_0(z)}{2A_0(z)-1}\left(\delta^2_zA_k(z)+\sum\limits_{i+j=k}\left(
A_i(z)\delta^2_zA_j(z)-\delta_zA_i(z)\delta_zA_j(z)\right)\right.\nonumber\\
 \left.\hphantom{A_{k+1}(z) =}{} -3z\sum\limits_{i+j=k+1}^{i,j\leq k}A_i(z)A_j(z)-z\sum\limits_{j_1+j_2+j_3=k+1}^{j_1,j_2,j_3\leq k}
A_{j_1}(z)A_{j_2}(z)A_{j_3}(z)\right).\label{eq:recAk}
\end{gather}
Consider equation~\eqref{eq:A0}. Differentiating it one finds
\begin{gather}\label{eq:diffA0}
\delta_zA_0(z)=A_0(z) \frac{1+A_0(z)}{1-2A_0(z)},\qquad
\delta_z^2\ln(1+A_0(z))=A_0(z) \frac{1+A_0(z)}{(1-2A_0(z))^3}.
\end{gather}
Note that variable $z$ can be excluded from equation~\eqref{eq:recAk} with the help of equation~\eqref{eq:A0}. Now, assuming
that $A_0(z)$ is given and putting successively into equation~\eqref{eq:recAk} $k=0$, $1$, and $2$, one, with the help of
equations~\eqref{eq:diffA0}, finds
\begin{gather}
A_1(z) =-A_0(z)\left(\frac{1+A_0(z)}{1-2A_0(z)}\right)^4,\label{eq:A1}\\
A_2(z) =A_0(z)\frac{(1+A_0(z))^6}{(1-2A_0(z))^9}\big(1+36A_0(z)+135A_0^2(z)+19A_0^3(z)\big),\label{eq:A2}\\
A_3(z) =-A_0(z)\frac{(1+A_0(z))^8}{(1-2A_0(z))^{14}}\big(1+216A_0(z)+5952A_0^2(z)+40875A_0^3(z) \nonumber\\
\hphantom{A_3(z) =}{} +77922A_0^4(z)+25821A_0^5(z)+1262A_0^6(z)\big).\label{eq:A3}
\end{gather}
Obviously, if $A_0(z)\neq1/2$, then recurrence relation~\eqref{eq:recAk} allows one to uniquely construct the sequence $A_k(z)$ for a given function $A_0(z)$. Moreover, by induction it is easy to prove that
\begin{gather}\label{eq:Ak-gen}
A_k(z)=-A_0(z)\frac{(1+A_0(z))^{2(k+1)}}{(-1+2A_0(z))^{5k-1}}R_{3(k-1)}(A_0(z)),
\end{gather}
where $R_{3(k-1)}(A_0)$ is a polynomial in $A_0$ with $\deg R_{3(k-1)}=3(k-1)$ and integer coefficients. One can further
\begin{Conjecture}For $k\geq2$: the coefficients of $R_{3(k-1)}(A_0)$ are positive integers,
$R_{3(k-1)}(A_0)=1+6^kA_0+\cdots$, and $R_{3(k-1)}(x)$ is irreducible over $Q[A_0]$.\footnote{This conjecture does not used in the following proof.}
\end{Conjecture}

Substitute $A_0(z)=1/2$ into equation~\eqref{eq:A0}, then we find $z=2^2/3^3$. Therefore $A_0(z)\neq1/2$ for the
function $A_0(z)$ defined by equation~\eqref{eq:A0} and $|z|< 2^2/3^3$. We are interested in a regular in
a neighbourhood of $z=0$ solution of equation~\eqref{eq:A0}. There are three solutions of this equation two of them
are singular at $z=0$: the latter solutions can be constructed as convergent series in $\sqrt{z}$ with the leading
terms $\pm1/\sqrt{z}$. The regular solution which vanishes at $z=0$ can be constructed as convergent Taylor series.
The solutions of the cubic equation can be also in the standard way expressed in terms of the trigonometric or
hyperbolic functions. Say, in our case, it is convenient to use hyperbolic functions
\begin{gather*}
A_0(z)+1\equiv p\sinh(\varphi) \quad \Rightarrow \quad
zp^3\sinh^3(\varphi)=p\sinh(\varphi)-1,\qquad zp^3\equiv 4q, \qquad
p\equiv-3q.
\end{gather*}
Thus,
\begin{gather*}
p=\frac2{\sqrt{-3z}},\qquad
-\frac1q=\frac{3\sqrt{-3z}}{2},\qquad
q\sinh(3\varphi)=-1 \quad \Rightarrow \\
3\varphi_k=\mathrm{arcsinh}\left(-\frac1q\right)+2\pi ik,\qquad k=0,\pm1.
\end{gather*}
To get the regular solution we have to take $k=0$ in the formula above,
\begin{gather}\label{eq:A0-sinh}
A_0(z)=-1+\frac2{\sqrt{-3z}}\sinh\left(\frac13 \mathrm{arcsinh}\left(\frac32\sqrt{-3z}\right)\right).
\end{gather}
Expanding equation~\eqref{eq:A0-sinh} into the Taylor series we get
\begin{gather}
A_0(z)=z+3z^2+12z^3+55z^4+273z^5+1428z^6+7752z^7\nonumber\\
\hphantom{A_0(z)=}{} +43263z^8+246675z^9+\cdots.\label{eq:A0-Taylor}
\end{gather}
Since the radius of convergence of Taylor series for the function $\mathrm{arcsine}(x)$ is~$1$,
we get that series~\eqref{eq:A0-Taylor} converges for $\left|\frac32\sqrt{-3z}\right|<1$, i.e., $|z|<2^2/3^3$.
Obviously, Taylor expansions for all functions $A_k(z)$, because of equation~\eqref{eq:Ak-gen}, converges in the
same disk as the Taylor expansion for $A_0(z)$.
\end{proof}

Denote as $A_k[n]$ coefficients of the Taylor expansions of the functions $A_k(z)$:
\begin{gather}\label{eq:Akn-def}
A_k(z)=(-1)^k\sum\limits_{n=1}^{\infty} A_k[n]z^n.
\end{gather}
Above we used the fact that $A_k(0)=0$, it follows from $A_0[0]=0$, see equation~\eqref{eq:A0-Taylor} and
equation~\eqref{eq:Ak-gen}; for $k\geq1$ this fact is also easy to establish directly from recurrence relation~\eqref{eq:recAk}.
\begin{Proposition}\label{prop:Akn-positive}
For all $k\geq0$ and $n\geq1$ the numbers $A_k[n]$ are positive integers. Moreover, $A_k[1]=1$.
\end{Proposition}
\begin{proof}
For $k=0$, $1$, $2$, $3$ the statement, obviously, follows from explicit formulae~\eqref{eq:A1}--\eqref{eq:A3}.
For the other values of $k$ and $n$ it can be established by mathematical induction by substitution of definition~\eqref{eq:Akn-def}
into recurrence relation~\eqref{eq:recAk}. In this proof one has to combine together two terms
\begin{gather*}
A_i(z)\delta^2_zA_j(z)-\delta_zA_i(z)\delta_zA_j(z)\qquad
\mathrm{and}\qquad
A_j(z)\delta^2_zA_i(z)-\delta_zA_j(z)\delta_zA_i(z),
\end{gather*}
to get sign definite contribution.

To prove that $A_k[1]=1$ is also easy by mathematical induction with the help of recurrence relation~\eqref{eq:recAk}.
Using this relation we see that $A_{k+1}[1]=A_k[1]$, because the terms proportional to $z$ come only from the first term,
$\delta^2_zA_k(z)$ of this relation.
\end{proof}

\begin{Remark}\label{rem:sequences-Ak}
Using explicit formulae~\eqref{eq:A1}--\eqref{eq:A3} with the help of Maple code we find for $n=1,2,\ldots,11,\ldots,$
the following first terms of the corresponding integer sequences:
\begin{align*}
A_0[n]\colon\quad
&1,3,12,55,273,1428,7752,43263,246675,1430715,8414640,\ldots,\\
A_1[n]\colon\quad
&1,15,162,1525,13308,110691,890724,6996474,53953605,410084004,\\
&3080715624,\ldots,\\
A_2[n]\colon\quad
&1,63,1674,30610,452619,5832225,68232648,743146326,7659571500,\\
&75562845204,719340288408,\ldots,\\
A_3[n]\colon\quad
&1,255,15924,546950,13372449,262072839,4394608056,65619977445,\\
&895717557900,11382479204349,136443463958412,\ldots.
\end{align*}
Since it is integer sequences it is natural to check
``The on-line encyclopedia of integer sequences''~\cite{OEIS}.
Actually, the sequence $A_0[n]$ is sequence $A001764$ of the encyclopedia, which
is the famous sequence of Fuss--Catalan numbers
\begin{gather}\label{eq:Fuss-Catalan}
A_0[n]=\frac{\dbinom{3n}{n}}{2n+1}.
\end{gather}
The other sequences at the time being are not included there.
\end{Remark}

\begin{Proposition}\label{prop:Akn-monotonicity}
For all $k=0,1,2,\ldots$, $A_k[n]$ are monotonically increasing sequences of $n$.
\end{Proposition}
\begin{proof}
The proof is achieved via mathematical induction. For $k=0$ the proof follows from
the explicit formula (see equation~\eqref{eq:Fuss-Catalan}). For sequences $A_1[n]$,
$A_2[n]$, and $A_3[n]$ the proof is based on explicit equations~\eqref{eq:A1}--\eqref{eq:A3}:

To do it one have to notice that, say, for $A_1[n]$ any two successive coefficients of the terms~$z^n$ and~$z^{n+1}$, i.e., $A_1[n]$ and $A_1[n+1]$, can be presented as the sum of positive terms.
Each term is a~product of the coefficients of series, with positive monotonically increasing terms,
corresponding to the powers of $z$ lower than~$n$ or $n+1$, respectively. I~mean the series~$A_0(x)$,
$(1+A_0(x))^4$, and $(1-2A_0(x))^{-4}$. Comparing these sums we find that the number of the terms in
the sum for the senior coefficient is larger, and for each term in the sum for the junior coefficient there
is a term in the senior sum which is the product of the same coefficients except one. That last coefficient
corresponds to the same series as the coefficient from the junior sum but has a~subscript greater by~$1$
comparing to it. The last fact means that the product from the senior sum is greater than the corresponding product
of the junior one. Analogous idea works for the sequences $A_2[n]$ and $A_3[n]$ and, after the inductive assumption,
for $A_{k+1}[n]$ with the help of recurrence relation~\eqref{eq:recAk}.
\end{proof}

\begin{Proposition}For $k=0,1,2,\ldots$
\begin{gather}\label{eq:Akn-n-asymptotics}
A_k[n]\underset{n\to+\infty}{=}
\frac{3^{(5-k)/2}}{2^{2k+3}}\frac{R_{3(k-1)}\left(\frac12\right)}{\Gamma\left(\frac{5k-1}{2}\right)}
n^{(5k-3)/2}\left(\frac{3^3}{2^2}\right)^n\left(1+o(1)\right),
\end{gather}
where, for $k=0$ we put formally $R_{-3}\left(\frac12\right)\equiv-\left(\frac23\right)^2$, polynomials $R_{3(k-1)}(A_0)$
are defined in equation~\eqref{eq:Ak-gen}, and $\Gamma(\cdot)$ is the Euler Gamma-function.
\end{Proposition}
\begin{proof}For the case $k=0$ the sequence $A_0[n]$ represents the Fuss--Catalan numbers (see Remark~\ref{rem:sequences-Ak}).
Therefore asymptotics for this sequence is easy to establish by applying the Stirling formula to
the explicit expression~\eqref{eq:Fuss-Catalan}.

The case $k\geq1$ is more deliberate. Using equation~\eqref{eq:A0} one proves the following estimate,
\begin{gather*}
A_0(z)\underset{z\to1}=\frac12 +\left(\frac32\right)^2\sqrt{1-z}+O\left(1-z\right).
\end{gather*}
Substituting this estimate into equation~\eqref{eq:Ak-gen}, having in mind definition~\eqref{eq:Akn-def}, and
Propositions~\ref{prop:Akn-positive} and \ref{prop:Akn-monotonicity} we find ourself in position to apply the
well-known Tauberian theorem by Hardy--Littlewood~\cite{HL} which implies the result stated in
equation~\eqref{eq:Akn-n-asymptotics}.
\end{proof}
\begin{Conjecture}
For $k=0,1,2,\dots$ the numbers $A_k[n]$ as $n\to+\infty$ approach to their asymptotic
value monotonically growing,
\begin{gather*}
\frac{A_k[n]}{n^{(5k-3)/2}}\left(\frac{3^3}{2^2}\right)^{-n}\nearrow
\frac{3^{(5-k)/2}}{2^{2k+3}}\frac{R_{3(k-1)}\left(\frac12\right)}{\Gamma\left(\frac{5k-1}{2}\right)},
\qquad\mathrm{as}\quad n\nearrow+\infty,
\end{gather*}
in particular, the error estimate in equation~\eqref{eq:Akn-n-asymptotics} is a negative number.
\end{Conjecture}
\begin{Remark}
According to equations~\eqref{eq:A1}--\eqref{eq:A3}
\begin{gather*}
R_0\left(\frac12\right)=1,\qquad
R_3\left(\frac12\right)=\frac{441}{8}=\frac{3^27^2}{2^3},\qquad
R_6\left(\frac12\right)=\frac{99225}{8}=\frac{3^45^27^2}{2^3}.
\end{gather*}
\end{Remark}

\begin{Remark}
Explicit formula for the Fuss--Catalan numbers~\eqref{eq:Fuss-Catalan}, $A_0[n]$, allows one
to find successively for $k=1,2,\ldots$ explicit expressions for sequences $A_k[n]$.
It would be interesting to find a general formula for these sequences for all $k$. In the following
Proposition we show how one can get formula for $A_1[n]$. The proof makes it clear how to extend this procedure
and successively obtain explicit formulae for $A_2[n]$, $A_3[n]$,\ldots.
\end{Remark}

\begin{Proposition}
\begin{gather}\label{eq:A1n-explicit}
A_1[n]=\frac{n+1}{18}\dbinom{3n+4}{n+1}\left(F(1,-n-1;2n+4;-2)-\frac{4n+6}{3n+4}\right),\qquad
n=1,2,\ldots,
\end{gather}
where $\tbinom{3n+4}{n+1}$ is the binomial coefficient and $F(a,b;c;z)$ is the Gauss hypergeometric function.
\end{Proposition}
\begin{proof}
We begin with the following observation, the first equation~\eqref{eq:diffA0} can be rewritten as
\begin{gather}\label{eq:1-2A0z-inv}
1+\frac23\left(2\delta_zA_0(z)+A_0(z)\right)=\frac1{1-2A_0(z)}.
\end{gather}
Thus one can obtain explicit formula for the coefficients of the Taylor expansion at $z=0$ of the function in the
r.h.s.\ of equation~\eqref{eq:1-2A0z-inv} in terms of the Fuss--Catalan numbers, defining corresponding Taylor
expansion in its l.h.s. Namely, after a simple calculation one finds
\begin{gather}\label{eq:genfunbinom2n-1:n}
\frac1{1-2A_0(z)}=\sum\limits_{n=0}^{\infty}\dbinom{3n-1}{n}z^n.
\end{gather}
This is not absolutely trivial result as it is equivalent to the following identity for the Fuss--Catalan
numbers~\eqref{eq:Fuss-Catalan}, $A_0[k_i]$,
\begin{gather*}
\dbinom{3n-1}{n}=\sum\limits_{l=1}^{l=n}2^l\sum\limits_{\substack{k_1+\ldots k_l=n\\k_1\geq1,\ldots,k_l\geq1}}
A_0[k_1]\cdot A_0[k_1]\cdots A_0[k_l].
\end{gather*}
L.h.s.\ of equation~\eqref{eq:genfunbinom2n-1:n} is the generating function for sequence~A165817 of~\cite{OEIS}.

It will be more convenient to consider the function
\begin{gather}\label{eq:series1standard}
\frac{1+A_0(z)}{1-2A_0(z)}=\frac12\left(\frac3{1-2A_0(z)}-1\right)=\sum\limits_{n=0}^{\infty}\dbinom{3n}{n}z^n.
\end{gather}
This is the generating function for sequence A005809 of \cite{OEIS}. The next step is the function
\begin{gather}\label{eq:series2standard}
\left(\frac{1+A_0(z)}{1-2A_0(z)}\right)^2=\sum\limits_{n=0}^{\infty}\sum\limits_{k=0}^{n}\dbinom{3k}{k}\dbinom{3(n-k)}{n-k}z^n.
\end{gather}
The coefficients of this expansion constitutes sequence A006256 of \cite{OEIS}. The sums with factorials often can be
presented in terms of special values of the hypergeometric functions
\begin{align*}
\sum\limits_{k=0}^{n}\dbinom{3k}{k}\dbinom{3(n-k)}{n-k}&=\dbinom{3n}{n}
{}_4F_3\begin{bmatrix}1/3, 2/3, 1/2-n, -n;&1\\ 1/2, 1/3-n, 2/3-n&\end{bmatrix}\\
&=\dbinom{3n+1}{n} {}_2F_1\begin{bmatrix}1,-n;&-2\\2n+2&\end{bmatrix}.
\end{align*}
The last two equalities can be found in A006256 of \cite{OEIS}.\footnote{Due to the contributions of
Jean-Fran\c cois Alcover and Peter Luschny to \cite{OEIS}.}
Now we can introduce the auxiliary function we need for our proof
\begin{gather*}%\label{eq:A-caligraphic}
\mathcal{A}(z)\equiv A_0(z)\left(\frac{1+A_0(z)}{1-2A_0(z)}\right)^2.
\end{gather*}
To find Taylor expansion of $\mathcal{A}(z)$ at $z=0$ one have to consider decomposition of this function
in partial fractions and apply the results obtained above,
\begin{align*}
\mathcal{A}(z)&=\frac34+\frac{A_0(z)}4-\frac{15}8\frac1{1-2A_0(z)}+\frac98\frac1{(1-2A_0(z))^2}\\
&=\frac58+\frac{A_0(z)}4-\frac98\frac1{1-2A_0(z)}+\frac12\left(\frac{1+A_0(z)}{1-2A_0(z)}\right)^2.
\end{align*}
Since the Taylor series of all functions in the r.h.s.\ of this equation are obtained above
we arrive at the following result
\begin{gather}\label{eq:mathcal-A-gauss}
\mathcal{A}(z)=\sum\limits_{n=0}^{\infty}\mathcal{A}_nz^{n+1},\qquad
\mathcal{A}_n=\frac12\binom{3n+4}{n+1}\left(F(1,-n-1;2n+4;-2)-1\right).
\end{gather}
The function $\mathcal{A}(z)$ is the generating function of sequence A075045 of \cite{OEIS}.

Differentiating $\mathcal{A}(z)$ we find
\begin{align}
\delta_z\mathcal{A}(z)&=\frac{9A_0^2(z)(1+A_0(z))}{(1-2A_0(z))^4}+\frac{3A_0(z)(1+A_0(z))}{4(1-2A_0(z))^3}
+\frac{A_0(z)(1+A_0(z))}{4(1-2A_0(z))}\nonumber\\
&=9zA_0(z)\left(\frac{1+A_0(z)}{1-2A_0(z)}\right)^4+\frac12\delta_z^2A_0(z)+\frac12\delta_zA_0(z)\nonumber\\
&=-9zA_1(z)+\frac12(\delta_z^2+\delta_z)A_0(z).\label{eq:A1z-linearization}
\end{align}
Comparing coefficients of the terms $z^{n+1}$ in equation~\eqref{eq:A1z-linearization} one proves
\begin{gather*}
A_1[n]=\frac{n+1}{18}\left(2\mathcal{A}_n-(n+2)A_0[n+1]\right).
\end{gather*}
After a simple calculation the last equation together with the second equation~\eqref{eq:mathcal-A-gauss} implies equation~\eqref{eq:A1n-explicit}.
\end{proof}
\begin{Proposition}
\begin{gather}
A_2[n] =\frac{n+1}{128}\dbinom{3n+1}{n}\left(\frac{168n^3+846n^2+1211n+510}5\right.\nonumber\\
\hphantom{A_2[n] =}{} -3(n+1)(25n+34)F(1,-n;2n+2;-2)\bigg),\qquad n=1,2,\ldots,\label{eq:A2n-explicit}
\end{gather}
where $\tbinom{3n+1}{n}$ is the binomial coefficient and $F(a,b;c;z)$ is the Gauss hypergeometric function.
\end{Proposition}
\begin{proof}Here, we present the proof in a more algorithmic way comparing with the one for the previous proposition.
This proof is easy to generalize for sequences $A_k[n]$ with $k=3,4,\ldots$.

We have two ``standard'' series~\eqref{eq:series1standard} and \eqref{eq:series2standard}, which we denote
$F_1(z)$ and $F_2(z)$, respectively. The idea is to present equation~\eqref{eq:A2} for $A_2[n]$ in terms of these series and their derivatives.
This can be done (easily with the help of Maple code) by decomposition of $A_2[n]$ on partial fractions
\begin{gather*}
A_2(z)=-\frac{19}{512}A_0(z)-\frac{19}{512}-\left(\frac{51}{64}+\frac{279}{128}\delta_z+
\frac{63}{32}\delta_z^2+\frac{75}{128}\delta_z^3\right)F_2(z)\\
\hphantom{A_2(z)=}{}+\left(\frac{427}{512}+\frac{4367}{1280}\delta_z+\frac{2853}{640}\delta_z^2+\frac{1479}{640}\delta_z^3
+\frac{63}{160}\delta_z^4\right)F_1(z).
\end{gather*}
Substituting \looseness=1 into the above equation known series for $A_0(z)$, $F_1(z)$, and $F_2(z)$ and equa\-ting~corres\-ponding coefficients one arrives, after the straightforward calculations, to equa\-tion \eqref{eq:A2n-explicit}.
\end{proof}
\begin{Proposition}\label{prop:U-A} Let $U\equiv U(a,x)$ be given by equations~\eqref{eq:U-series}--\eqref{eq:u-2n-recurrence}. For $a\in\mathbb{R}$
\begin{gather*}%\label{eq:asympU-A}
U\big(a,a^2z\big)\underset{a\to+\infty}\sim A(a,z).
\end{gather*}
\end{Proposition}
\begin{proof}Recall that function $U(a,x)$ is related, via equation~\eqref{eq:U-x} for $b=a$, with the solu\-tion~$u(\tau)$
defined in Theorem~\ref{th:existence & uniqueness}. Thus, it is a meromorphic function of $x$ and holomorphic of~$a$ in any simply connected domain of $\mathbb{C}{\setminus} i\mathbb{Z}$.

Mathematical induction with the help of recurrence relation~\eqref{eq:u-2n-recurrence} and initial coefficients
given in Remark~\ref{rem:coeff3} allows one to prove that $u_{2n}(a)=O\left(1/a^{2n}\right)$ as $a\to+\infty$.
We can surely prove even more, that this estimate is valid for $|\arg a|<\pi/2-\varepsilon$ for some
sufficiently small $\varepsilon>0$ and the analogous result for the left complex semiplane. After that it is
natural to consider the change of independent variable, $x=a^2z$, and develop the rational function
$u_{2n}(a)a^{2n}$ into the asymptotic series as $a\to\infty$ in the corresponding semiplane. Thus we get
\begin{gather*}
u_{2n}(a)a^{2n}=\sum\limits^\infty_{k=0}\frac{u^k_{2n}}{a^{2k}}.
\end{gather*}
Consider the difference
\begin{gather*}
U\big(a,a^2z\big)-\sum\limits^\infty_{n=0}u^0_{2n}z^n=\sum\limits^\infty_{n=0}\big(u_{2n}(a)a^{2n}-u^0_{2n}\big)z^n.
\end{gather*}
Both series in r.h.s.\ of this equation are convergent, it follows from the proof of
Proposition~\ref{prop:taylor-general}, see the choice of the constant $C_1$ underneath inequality~\eqref{ineq:sum}:
obviously this constant is proportional to $1/|a|^2$ for the large values of~$|a|$. The last fact leads to the
finite radius of convergence of our series. Since $a^2\big(u_{2n}(a)a^{2n}-u^0_{2n}\big)=O(1)$ as the function of~$a$,
by construction has the same radius of convergence as the one without $a^2$. Therefore,
\begin{gather*}
U\big(a,a^2z\big)-\sum\limits^\infty_{n=0}u^0_{2n}z^n
=\frac1{a^2}\sum\limits^\infty_{n=0}a^2\big(u_{2n}(a)a^{2n}-u^0_{2n}\big)z^n=O\left(\frac1{a^2}\right).
\end{gather*}
We can inductively continue this construction and arrive at the following asymptotic expansion
\begin{gather*}
U\big(a,a^2z\big)=\sum\limits^\infty_{n=0}\frac{a_k(z)}{a^{2k}},
\end{gather*}
where the functions $a_k(z)=\sum\limits^\infty_{n=0}u^k_{2n}z^k$ are given by the convergent series.
It follows from Proposition~\ref{prop:Ak-formal} that this expansion coincides with \eqref{eq:A-series}.
\end{proof}
\begin{Corollary}\label{cor:u-2n-Akn}
\begin{gather}\label{eq:u-2n-Akn}
u_{2n}(a)a^{2n}=\sum\limits_{k=0}^{\infty}(-1)^k\frac{A_k[n]}{a^{2k}},
\end{gather}
where the series is absolutely convergent for $|a|>n$.
\end{Corollary}
\begin{proof}
According to Proposition~\ref{prop:U-A} equation~\eqref{eq:u-2n-Akn} holds
for arbitrary $|z|<2^2/3^3$ and posi\-tive~$a$. In fact, it is valid for complex $a$ as follows from the monodromy lemma for
analytic functions.
The convergence of series~\eqref{eq:u-2n-Akn} is a consequence of the fact that $u_{2n}(a)$ is the rational function with
the largest poles at $a^2=-n^2$.
\end{proof}
\begin{Remark}Because of the convergence of series~\eqref{eq:u-2n-Akn} it is clear that $A_k[n] = O \big(n^{2k+o(k)} \big)$
as $k\to+\infty$.
Using definition~\eqref{eq:Akn-def} and recurrence relation~\eqref{eq:recAk}, it is easy to establish that $A_k[2]=2^{2k+2}-1$.
The last sequence is also presented in OEIS~\cite{OEIS} as the sequence $A024036$. The explicit formulae for $A_k[n]$ as $n>2$
is not that easy to establish. However, one can conjecture the following asymptotic estimate:
\begin{Conjecture}\label{con:A-kn-k-asymptotics}
\begin{gather}\label{eq:A-kn-k-asymptotics}
A_k[n]\underset{k\to+\infty}=\frac{n^{2k+3(n-1)}}{2^{n-1}\left((n-1)!\right)^3}
\left(1+O\left(\frac1{n^{2k/c_n}}\right)\right), \qquad n\geq2,
\end{gather}
and the numbers $c_n\nearrow\infty$ as $n\to+\infty$: $c_2=1$, $c_3\approx2.7$, and $c_n>n$ for $n\geq4$.
\end{Conjecture}
Thus, as $n$ growth, asymptotics~\eqref{eq:A-kn-k-asymptotics} provides a good numerical approximation of sequence
$A_k[n]$ for the larger values of~$k$.
\end{Remark}
\begin{Corollary}\label{cor:mn-sum-nk-n}
The numbers $m(n)$ defined in Conjecture~{\rm \ref{con:expansion-u-main}} satisfy the following relation
\begin{gather}\label{eq:mn-nk-n}
m(n)+n=\sum\limits_{k=1}^n n_k.
\end{gather}
\end{Corollary}
\begin{proof}
Equation~\eqref{eq:mn-nk-n} follows from Corollary~\ref{cor:u-2n-Akn} and equation~\eqref{eq:u2ka}. Note that
in fact it is proved in the beginning of the proof of Proposition~\ref{prop:U-A} via mathematical induction.

The nontrivial part of Conjecture~\ref{con:expansion-u-main} includes the explicit
expressions for $n_k$ (see the first equation~\eqref{eq:n-k-m(n)}), the second equation~\eqref{eq:n-k-m(n)}
is confirmed.
\end{proof}
\begin{Corollary}Assume Conjecture~{\rm \ref{con:A-kn-k-asymptotics}} is true. Then numbers $A_k[n]$ as $k\to+\infty$ are approaching
to their asymptotic value monotonically growing,
\begin{gather*}
\frac{A_k[n]}{n^{2k}}\nearrow\frac{n^{3(n-1)}}{2^{n-1}\left((n-1)!\right)^3},\qquad\mathrm{as}\quad k\nearrow+\infty,
\end{gather*}
in particular, the error estimate in equation~\eqref{eq:A-kn-k-asymptotics} is a negative number.
\end{Corollary}
\begin{proof}
The statement follows from Conjecture~\ref{con:A-kn-k-asymptotics} and the estimate
\begin{gather*}
A_{k+1}[n]>n^2A_k[n],
\end{gather*}
which is easy to deduce from recurrence relation~\eqref{eq:recAk}.
\end{proof}

Now, we consider application of equation~\eqref{eq:u-2n-Akn} for calculation of coefficients $p_k(n)$
defined in equation~\eqref{eq:polynom-P}. We begin with a practical comment, $\sum\limits_{k=1}^n n_k$
in r.h.s.\ of equation~\eqref{eq:mn-nk-n} coincides with the sum of the elements in the $n$-th row of the semi-infinite matrix $M$ constructed in the following way:
$M\equiv(M_1M_2\ldots M_k\ldots)$, where $M_k$ are the semi-infinite columns:
\begin{align*}
M_1&=(112233\ldots mm\ldots)^{\rm T},\qquad
M_2=(0111222333\ldots mmm\ldots)^{\rm T},\\
M_k&=(\underset{k-1}{\underbrace{0\ldots0}} \underset{k+1}{\underbrace{1\ldots1}} \underset{k+1}{\underbrace{2\ldots2}} \ldots
\underset{k+1}{\underbrace{m\ldots m}} \ldots)^{\rm T},\qquad\ldots.
\end{align*}
Now, comparing equations~\eqref{eq:u-2n-Akn}, \eqref{eq:u2ka}, and \eqref{eq:polynom-P} and denoting $\check{a}\equiv1/a^2$
we find
\begin{gather}\label{eq:p-k-A-k-generating-equation}
\sum\limits_{k=0}^{m(n)}p_{m(n)-k}\check{a}^k=\prod\limits_{k=0}^n\big(1+k^2\check{a}\big)^{n_k}
\sum\limits_{k=0}^{\infty}(-1)^kA_k[n]\check{a}^{k}.
\end{gather}
To use this relation it is convenient to introduce numbers $q_k(n)$, $k=0,1,\ldots, m(n)+n$, as the coefficients of the
polynomial,
\begin{gather}\label{def-qkn}
Q_n(\check{a})\equiv
\prod\limits_{k=0}^n\big(1+k^2\check{a}\big)^{n_k}=\sum\limits_{k=0}^{m(n)+n}q_k(n)\check{a}^k.
\end{gather}
Obviously,
\begin{gather*}
q_0(n)=1,\qquad
q_1(n)=\sum\limits_{k=1}^{n}n_kk^2.
\end{gather*}
More generally, denote
\begin{gather*}
S_l=\frac{(-1)^{l+1}}{l}\sum\limits_{k=1}^{n}n_kk^{2l},
\end{gather*}
so that
\begin{gather*}
q_k(n)=\sum\limits_{\substack{i_1+2i_2+\dots+ki_k=k\\i_1\geq0,i_2\geq0,\ldots,i_k\geq0}}
\frac{S_1^{i_1}S_2^{i_2}\cdots S_k^{i_k}}{i_1!i_2!\cdots i_k!}.
\end{gather*}
The first few polynomials $Q_n(\check{a})$ are as follows
\begin{gather*}
Q_1(\check{a})=1+\check{a},\qquad
Q_2(\check{a})=1+5\check{a}+4\check{a}^2,\qquad
Q_3(\check{a})=1+15\check{a}+63\check{a}^2+85\check{a}^3+36\check{a}^4, \\
Q_4(\check{a})=1+31\check{a}+303\check{a}^2+1093\check{a}^3+1396\check{a}^4+576\check{a}^5,\\
Q_5(\check{a})=1+61\check{a}+1362\check{a}^2 +14282\check{a}^3+76373\check{a}^4+213753\check{a}^5
+306664\check{a}^6\\
\hphantom{Q_5(\check{a})=}{} +213904\check{a}^7+57600\check{a}^8,\qquad\ldots.
\end{gather*}
Identity~\eqref{eq:p-k-A-k-generating-equation} implies the following equations:
\begin{gather}
p_{m(n)-k}(n) =\sum\limits_{i=0}^k(-1)^{k-i}q_i(n)A_{k-i}[n],\qquad k=0,1,\ldots,m(n),\label{eq:pmn-in-Ak}\\
0 =\sum\limits_{i=0}^k(-1)^{k-i}q_i(n)A_{k-i}[n],\qquad m(n)+1\leq k\leq m(n)+n,\nonumber\\
0 =\sum\limits_{i=0}^{m(n)+n}(-1)^{k-i}q_i(n)A_{k-i}[n],\qquad k\geq m(n)+n.\nonumber
\end{gather}
\begin{Corollary}
The numbers $p_{m(n)}(n)$ coincide with the Fuss--Catalan numbers
\begin{gather}
 p_{m(n)}(n)=A_0[n]=\frac{\dbinom{3n}{n}}{2n+1},\label{eq:p-mn-Fuss-Catalan}\\
 p_{m(n)-1}(n)=q_1(n)A_0[n]-A_1[n]\label{eq:pmn-1-via-Ak}\\
 =A_0[n] \left(\sum\limits_{k=0}^n n_kk^2
+\frac{(3n+1)(3n+2)}{6} \left(1-\frac{3n+4}{6(2n+3)}F(1,3n+5;2n+4;2/3)\right) \right)\label{eq:p-mn-1-hypergeo}\\
 =A_0[n]\left(\frac{(3n+1)(3n+2)}{6}+\sum\limits_{k=1}^{n}n_kk^2\right)
-\frac{n+1}{18}\sum\limits_{k=0}^{n+1}\binom{3n+4}{k}2^{n+1-k},\label{eq:p-mn-1sums}
\\
\label{eq:pmn-2-via-Ak}
p_{m(n)-2}(n) =q_2(n)A_0[n]-q_1(n)A_1[n]+A_2[n]\\
\hphantom{p_{m(n)-2}(n)}{} =\frac{A_0[n]}{2}\left(\left(\sum\limits_{k=0}^n n_kk^2\right)^2-\sum\limits_{k=0}^n n_kk^4\right.\nonumber\\
 \hphantom{p_{m(n)-2}(n)=}{}+\frac{(3n+1)(3n+2)}{3}
\left(1-\frac{3n+4}{2(2n+3)}F(1,-n-1;2n+4;-2)\right)\sum\limits_{k=0}^n n_kk^2\nonumber\\
 \hphantom{p_{m(n)-2}(n)=}{} +\frac{(3n+1)(n+1)}{64}\left(\frac{168n^3+846n^2+1211n+510}5\right.\nonumber\\
 \hphantom{p_{m(n)-2}(n)=}{} -3(n+1)(25n+34)F(1,-n;2n+2;-2)\bigg).\label{eq:p-mn-2-hypergeo}
\end{gather}
\end{Corollary}
\begin{proof}
The expressions of coefficients $p_{m(n)}$, $p_{m(n)-1}$, and $p_{m(n)-2}$ in terms of $A_k[n]$ are just special cases of
equation~\eqref{eq:pmn-in-Ak}. The explicit formula for the Fuss--Catalan numbers (see equation~\eqref{eq:p-mn-Fuss-Catalan})
is given in \cite[sequence~$A001764$]{OEIS}).

Equations~\eqref{eq:p-mn-1-hypergeo} and \eqref{eq:p-mn-1sums} follows from equation~\eqref{eq:A1n-explicit} with the help
of the well-known relations for the Gauss hypergeometric function (see~\cite{BE}):
\begin{gather*}
F(1,-n-1;2n+4;-2)=\frac13F(1,3n+5;2n+4;2/3)=\frac{2^{n+1}}{\binom{3n+4}{n+1}}\sum\limits_{k=0}^{n+1}\frac{\binom{3n+4}{k}}{2^k}.\tag*{\qed}
\end{gather*}\renewcommand{\qed}{}
\end{proof}

\begin{Remark}
The definition of numbers $q_k(n)$ (see equation~\eqref{def-qkn}) and Proposition~\ref{prop:Akn-positive} imply that
$q_k(n)$ and $A_k[n]$ are positive integers, therefore equations~\eqref{eq:pmn-1-via-Ak} and \eqref{eq:pmn-2-via-Ak}
show that $p_{m(n)-1}(n)$ and $p_{m(n)-2}(n)$ are integers. However, the fact that they are positive is not that obvious.
Moreover, it is not immediate to see that the explicit expressions for coefficients $p_{m(n)-1}(n)$ and $p_{m(n)-2}(n)$,
given by equations~\eqref{eq:p-mn-1sums} and \eqref{eq:p-mn-2-hypergeo}, are positive integers. Let us confirm
Conjecture~\ref{con:expansion-u-main} (see equation~\eqref{eq:polynom-P}) for $p_{m(n)-1}(n)$. The case $p_{m(n)-2}(n)$
can be studied analogously.

We recall that the numbers $p_k(n)$ are the coefficients of the polynomial $P_{m(n)}\big(a^2\big)$ (see equation~\eqref{eq:polynom-P})
so that they are not defined, or, formally, can be put equal to zero for $k<0$, or $k>m(n)$.
Since $m(1)=m(2)=0$ (see Remark~\ref{rem:m(n)}) we have $p_{m(1)-1}(1)=p_{m(2)-1}(2)=0$. Note that
expression~\eqref{eq:p-mn-1sums} vanishes for $n=1$ and $2$.
\end{Remark}
\begin{Proposition}\label{prop:pmn-1>0}
\begin{gather*}
p_{m(n)-1}(n)>0,\qquad\textrm{for}\quad n\geq3.
\end{gather*}
\end{Proposition}
\begin{proof}Consider equation~\eqref{eq:p-mn-1sums}. If $n$ is odd, then
\begin{gather*}
\sum\limits_{k=1}^nn_kk^2\geq\sum\limits_{k=1}^nk^2+\sum\limits_{k=1}^{(n-1)/2}k^2=
\frac{n(n+1)(3n+1)}8\geq\frac{n\big(3n^2+3n+2\big)}8.
\end{gather*}
If $n$ is even, then
\begin{gather*}
\sum\limits_{k=1}^nn_kk^2\geq\sum\limits_{k=1}^{n}k^2+\sum\limits_{k=1}^{n/2-1}k^2>\frac{n\big(3n^2+3n+2\big)}8.
\end{gather*}
To prove that $p_{m(n)-1}(n)>0$ it is enough to prove that
\begin{gather}\label{ineq:proof-pmn-1-positive}
A_0[n]\sum\limits_{k=1}^{n}n_kk^2>\frac{n+1}{18}2^{n+1}\left(1+\frac12\right)^{3n+4}>
\frac{n+1}{18}\sum\limits_{k=0}^{n+1}\binom{3n+4}{k}2^{n+1-k}.
\end{gather}
Since the above inequality looks cumbersome, it is natural to consider a simpler inequality
\begin{gather}\label{ineq:simpleproof-pmn-1-positive}
\binom{3n}{n-1}\big(3n^2+3n+2\big)>(n+1)\cdot\frac{3^{3n+2}}{2^{2n+1}},
\end{gather}
The last inequality with the help of equation~\eqref{eq:p-mn-Fuss-Catalan} implies
inequality~\eqref{ineq:proof-pmn-1-positive}. To study inequa\-li\-ty~\eqref{ineq:simpleproof-pmn-1-positive}, we
introduce variable
\begin{gather*}
X_n\equiv\binom{3n}{n-1} \frac{2^{2n+1}}{3^{3n+2}} \frac{3n^2+3n+2}{n+1},\qquad n=1,2,\ldots,
\end{gather*}
and prove that it is a monotonically increasing sequence. Actually, the straightforward calculation shows
that
\begin{gather*}
X_{n+1}=\frac{(n+1/3)(n+2/3)(n+1)\big(n^2+3n+8/3\big)}{n(n+3/2)(n+2)\big(n^2+n+2/3\big)}X_n.
\end{gather*}
Since
\begin{gather*}
(n+1/3)(n+2/3)(n+1)\big(n^2+3n+8/3\big)-n(n+3/2)(n+2)\big(n^2+n+2/3\big)\\
\qquad{}=\frac12n^4+\frac{49}{18}n^3+\frac{35}{9}n^2+\frac{52}{27}n+\frac{16}{27}>0,
\end{gather*}
we have established the monotonicity of $X_n$. In case $X_1\geq1$ we would finish the proof, however
\begin{gather*}
X_1=\frac{2^5}{3^5}=0,1316\ldots,\qquad
X_2=\frac{1280}{6561}=0,1950\ldots,\\
X_{38}=0,9888\ldots,\qquad X_{39}=1,002\ldots,
\end{gather*}
where the last two calculations have been done with the help of Maple code. Thus we see, that the above proof
works for $p_{m(n)-1}(n)$ with $n\geq39$. Positiveness of $p_{m(n)-1}(n)$ with $n<39$ should be established
directly. For $n\leq5$ it follows from explicit expressions presented right after Remark~\ref{rem:coeff3}.
Positiveness of $p_{m(n)-1}(n)$ for $n=6,\ldots,38$ should be checked directly with the help of Maple code and
equation~\eqref{eq:p-mn-1-hypergeo}.
\end{proof}

Of course, explicit calculation of so many coefficients rises a desire to improve the above proof. This refinement is presented below.
\begin{proof}
One writes a more accurate estimate of the r.h.s.\ of equation~\eqref{eq:p-mn-1sums}
\begin{gather}
p_{m(n)-1}(n) \geq \binom{3n}{n-1} \left(\frac{(3n+1)(3n+2)}{6n}+\frac{3n^2+3n+2}{8}\right)\nonumber\\
\hphantom{p_{m(n)-1}(n) \geq}{} -\frac{n+1}{18}\left(\frac{3^{3n+4}}{2^{2n+3}}
- \sum\limits_{k=n+2}^{n+6} \binom{3n+4}{k}2^{n+1-k} \right) .\label{ineq:pmn-1}
\end{gather}
Note that $n+6\leq3n+4$ for $n\geq1$, so that this estimate works for all natural~$n$. Moreover, equality in~\eqref{ineq:pmn-1} takes place only for $n=1$.

Our goal is to prove that the expression in r.h.s.\ of inequality~\eqref{ineq:pmn-1} is positive for all $n\geq3$.
For $n=1$, $2$, and $3$ this expression equals
\begin{gather*}
0,\qquad-\frac{7}{256}=-0,027\ldots,\qquad\frac{18109}{2304}=7,859\ldots,
\end{gather*}
respectively. As in the previous proof we are going to use mathematical induction. To this end we rewrite our
statement positiveness of r.h.s.\ of \eqref{ineq:pmn-1} in the following way
\begin{gather}\label{ineq:YZn+1}
Y_n+Z_n>n+1,
\end{gather}
where
\begin{gather*}
Z_n=\frac{8}{81}\left(\frac{4}{27}\right)^n(n+1)\sum\limits_{l=1}^5\binom{3n+4}{l+n+1}\frac1{2^l},
\end{gather*}
and
\begin{gather*}
Y_n=\frac{8}{9}\left(\frac{4}{27}\right)^n\binom{3n}{n-1}\left(\frac{(3n+1)(3n+2)}{3n}+\frac{3n^2+3n+2}{4}\right).
\end{gather*}
Our nearest steps towards the proof of inequality~\eqref{ineq:YZn+1} is to establish the monotonic growth of the
sequences $Z_n$ and $Y_n$ and study how their members changing with $n$.

Consider $Z_n$. After multiplication on the common factor in front of the sum we can consi\-der~$Z_n$ as the sum of 5
entries (labeled by~$l$). For each entry we consider the ratio of its successive values with the change of~$n$:
\begin{gather*}
\frac{(n+2)(n+7/2)(n+2)(n+5/2)}{(n+1)(n+2+l)(n+2-l/2)(n+5/2-l/2)}.
\end{gather*}
The difference of the numerator and denominator of the above ratio is positive for all natural $n$ and $l=1,2,3,4,5$:
\begin{gather*}
\frac{n^3}{2}+\left(\frac{61}{18}-\frac{l}{4}+\frac{3l^2}{4}\right)n^2+
\left(\frac{68}{9}-\frac{3l}{4}+\frac{5l^2}{2}-\frac{l^3}{4}\right)n+
\left(\frac{50}{9}-\frac{l}{2}+\frac{7l^2}{4}-\frac{l^3}{4}\right),
\end{gather*}
because each bracket (actually each difference) above is positive for $1\leq l\leq5$. Thus $Z_n$ is the sum of monotonically
growing sequences and therefore is monotonically growing itself. It is easy to find asymptotics
\begin{gather*}
Z_n\underset{n\to\infty}=5\sqrt{\frac{3n}{4\pi}}+O\left(\frac{1}{\sqrt{n}}\right).
\end{gather*}
Therefore, $Z_{n+1}-Z_n$ vanishes as $n\to\infty$.

Because of the last property of the sequence $Z_n$ we have to prove a stronger monotonicity property for the
sequence $Y_n$, namely,
\begin{gather}\label{eq:DeltaYn>1}
Y_{n+1}-Y_n\geq1,
\end{gather}
otherwise the mathematical induction process cannot be launched
\begin{gather*}
\Delta Y_n\equiv Y_{n+1}-Y_n=\frac{2\big(243n^4+1170n^3+1773n^2+1014n+200\big)}{243(n+1)(2n+1)(2n+3)}\binom{3n}{n}\left(\frac{4}{27}\right)^n>0.
\end{gather*}
Therefore, the sequence $Y_n$ is monotonically growing. Asymptotics
\begin{gather*}
\Delta Y_n\underset{n\to\infty}=\sqrt{\frac{3n}{16\pi}}+O\left(\frac{1}{\sqrt{n}}\right),
\end{gather*}
shows that for validity of condition~\eqref{eq:DeltaYn>1} (at least beginning from some rather large $n$)
it is enough to prove monotonicity of $\Delta Y_n$:
\begin{gather*}
\frac{\Delta Y_{n+1}}{\Delta Y_n}-1=
\frac{2187n^5+12312n^4+25497n^3+23616n^2+7908n-400}{9(n+2)(2n+5)(243n^4+1170n^3+1773n^2+1014n+200)}>0,
\end{gather*}
because $n>1$. Calculation with the help of Maple code shows that
\begin{gather*}
\Delta Y_{13}=0,992\ldots,\qquad\Delta Y_{14}=1,021\ldots.
\end{gather*}
Thus, to prove the base of induction we have to check the validity of inequality~\eqref{ineq:YZn+1} for
$n=3,4,\ldots,14$. These calculations surely can be done by hands, however, it is much faster to make them
with Maple code:
\begin{gather*}
\big(Y_n+Z_n-(n+1)\big)\Big\vert_{n=3}^{14}=0,045\ldots, 0,133\ldots, 0,258\ldots, 0,414\ldots, 0,599\ldots,0,810\ldots,\\
1,045\ldots, 1,303\ldots, 1,584\ldots, 1,886\ldots, 2,208\ldots, 2,552\ldots.\tag*{\qed}
\end{gather*}\renewcommand{\qed}{}
\end{proof}
\begin{Remark}Although both proofs of Proposition~\ref{prop:pmn-1>0} follow the same scheme, an interesting feature of the second proof
is that it avoids explicit calculation of the coefficients $p_{m(n)-1}$, while this calculation is needed for the first one.
\end{Remark}

\section[Generating function $B(a,x)$]{Generating function $\boldsymbol{B(a,x)}$}\label{sec:generating-B}
Consider Taylor expansion of the coefficients $u_{2n}(a)$:
\begin{gather}\label{eq:u{2n}(a)-taylor}
u_{2n}(a)=u_{2n}^0+u_{2n}^1a^2+u_{2n}^2a^4+O\big(a^6\big).
\end{gather}
\begin{Proposition}\label{prop:u-at0}
\begin{gather}\label{eq:u2k0}
u_{2n}^0=u_{2n}(0)=\frac{n+1}{2^n}.
\end{gather}
\end{Proposition}
\begin{proof}
Firstly, put in equation~\eqref{eq:dp3} and expansion~\eqref{eq:zero-expansion} $a=b$, secondly, $a=0$.
Then equation~\eqref{eq:dp3} reduces to its integrable version (recall our convention equation~\eqref{eq:epsilon})
\begin{gather}\label{eq:dP3a=b=0}
u^{\prime \prime}(\tau) =\frac{(u^{\prime}(\tau))^{2}}{u(\tau)} - \frac{u^{\prime}(\tau)}{\tau} -\frac{8u^2(\tau)}{\tau},
\end{gather}
and expansion~\eqref{eq:zero-expansion} reads
\begin{gather}\label{eq:tracesolution}
u(\tau)=-\frac\tau2\left(1+\sum\limits_{n=1}^{\infty}u_{2n}(0)\tau^{2n}\right).
\end{gather}
Equation~\eqref{eq:dP3a=b=0} have the following general
\begin{gather}\label{eq:generalsolution}
u(\tau)=-\frac{C_1C_2}{4}\frac{\tau^{\sqrt{C_1}-1}}{\big(1-C_2\tau^{\sqrt{C_1}}\big)^2}
\end{gather}
and special
\begin{gather*}
u(\tau)=\frac{1}{\tau(C_3+2i\ln\tau)^2}
\end{gather*}
solutions, where $C_1$, $C_2$, and $C_3$ are complex parameters.

Comparing expansion~\eqref{eq:zero-expansion} with equation~\eqref{eq:generalsolution} one finds
\begin{gather*}
\sqrt{C_1}-1=1 \quad\Rightarrow\quad C_1=4 \quad\Rightarrow\quad C_2=1/2.
\end{gather*}
For these values of the parameters, equation~\eqref{eq:generalsolution} takes the following form
\begin{gather}\label{eq:general-special}
u(\tau)=-\frac{\tau}{2}\frac1{(1-\tau^2/2)^2}.
\end{gather}
Expanding now equation~\eqref{eq:general-special} into the Taylor series at $\tau=0$ and comparing it with
the expansion~\eqref{eq:zero-expansion} one arrives at equation~\eqref{eq:u2k0}.
\end{proof}
\begin{Remark}It is interesting to notice that solution~\eqref{eq:tracesolution} of equation~\eqref{eq:dP3a=b=0}
is the memory of this limiting equation about the last two terms of equation~\eqref{eq:dp3} which
disappeared in the limit, $b=a\to0$. The original solution~\eqref{eq:zero-expansion} is defined
by the condition of cancelation of the singularity at $\tau=0$ related with the presence of these
two terms. It is clear that after the limit the terms that disappear cannot affect on the solutions
of the equation obtained in the limit, especially taking into account that the other members of
the equation remained unchanged. Nevertheless, among solutions of the limiting equation there is the
one which remember about the disappeared terms!
\end{Remark}
\begin{Corollary}
\begin{gather}\label{eq:p0n}
p_0(n)=P_{m(n)}(0)=\frac{n+1}{2^n}\prod\limits_{k=1}^nk^{2n_k}.
\end{gather}
\end{Corollary}
\begin{proof}
Follows from comparison of equation~\eqref{eq:u2k0} with expansion~\eqref{eq:u2ka}.
\end{proof}
\begin{Remark}
Sure equation~\eqref{eq:p0n} is proved modulo the explicit expressions for the numbers $n_k$
given in the first equation~\eqref{eq:n-k-m(n)} of Conjecture~\ref{con:expansion-u-main}. Analogous comment
concerns Corollaries~\ref{cor:p1n} and \ref{cor:p2n}.
\end{Remark}
\begin{Remark}
The first terms of the integer sequence $P_{m(n)}(0)$ are
\begin{gather*}
1, \ 3, \ 18, \ 180, \ 10800, \ 226800, \ \ldots.
\end{gather*}
\end{Remark}

We can generalize the idea employed in Proposition~\ref{prop:u-at0} and calculate generating functions for further
terms of the Taylor expansion of the coefficients $u_{2n}(a)$ at $a=0$. These functions allow one to calculate
integer sequences of coefficients of the polynomial $P_{m(n)}\big(a^2\big)$ and in that sense represent the generating
functions for these sequences. Below we consider this construction.

We put in equations~\eqref{eq:dp3} and~\eqref{eq:zero-expansion} as above $a=b$, and rearrange summation in the last equation
such that we can present it in the following form
\begin{gather}\label{eq:u-B}
u(\tau)=-\frac\tau2(1+B(a,x)),
\end{gather}
where
\begin{gather*}%\label{eq:B-expansion}
B\equiv B(a,x)=\sum\limits_{k=0}^{\infty}a^{2k}B_k(x),\qquad x=\tau^2.
\end{gather*}
We call $B$ the (super)generating function for the Taylor expansions of coefficients $u_{2n}(a)$.
In this notation Proposition~\ref{prop:u-at0} can be reformulated as
\begin{gather}\label{eq:B0-solution}
1+B_0(x)=\frac1{\big(1-\frac{x}2\big)^2}
\end{gather}
in accordance with equation~\eqref{eq:general-special}.

Our goal now is to calculate the further terms of this expansion. Substituting $u(\tau)$ given by equation~\eqref{eq:u-B} into equation~\ref{eq:dp3} we find ODE for $B$:
\begin{gather}\label{eq:diff-eq-B}
\delta_x^2\ln(1+B)=xB-\frac{a^2B}{(1+B)^2},\qquad\delta_x=x\frac{{\rm d}}{{\rm d}x}.
\end{gather}
Substituting expansion~\eqref{eq:u-B} into equation~\eqref{eq:diff-eq-B} one finds for $B_0$
\begin{gather*}
\delta_x^2\ln(1+B_0)=xB_0.
\end{gather*}
The appropriate solution for this equation is given by equation~\eqref{eq:B0-solution}. Then, for $B_1$
\begin{gather}\label{eq:B1}
\delta_x^2\ln\left(\frac{B_1}{1+B_0}\right)-xB_1=-\frac{B_0}{(1+B_0)^2}
\end{gather}

For the other coefficients we get the following recurrence system ($k=2,3,\ldots$) of inhomogeneous ODEs
\begin{gather}
 \delta_x^2\left(\frac{B_k}{1+B_0}\right)-xB_k\label{eq:B-recurrence-homo}\\
= \sum_{\substack{n_1i_1+\dots+n_pi_p=k\\k-1\geq i_1>i_2>\dots>i_p\geq1}}
(-1)^{n_1+\dots+n_p} \frac{(n_1+\dots+n_p-1)!}{n_1!\cdots n_p!}
\delta_x^2\left(\frac{B_{i_1}^{n_1}\cdots B_{i_p}^{n_p}}{(1+B_0)^{n_1+\dots+n_p}}\right)
\label{eq:B-recurrence-rhs1}\\
+\!\!\!\sum_{\substack{n_1i_1+\dots+n_pi_p=k-1\\k-1\geq i_1>i_2>\dots>i_p\geq1}}\!\!\!\!\!
(-1)^{n_1+\dots+n_p} \frac{(n_1+\dots+n_p)!}{n_1!\cdots n_p!}
\frac{(n_1+\dots+n_p-B_0)}{(1+B_0)^2}\frac{B_{i_1}^{n_1}\cdots B_{i_p}^{n_p}}{(1+B_0)^{n_1+\dots+n_p}}.\!\!\!\!\label{eq:B-recurrence-rhs2}
\end{gather}
The right-hand side of this equation consists of two terms~\eqref{eq:B-recurrence-rhs1} and \eqref{eq:B-recurrence-rhs2}.
The sum in the first term~\eqref{eq:B-recurrence-rhs1} is taken over Young diagrams representing the partitions of $k$ such that
the length of the rows does not exceed $k-1$. The sum in the second term~\eqref{eq:B-recurrence-rhs2} is taken over all
Young diagrams representing the partitions of $k-1$.

If we equate to $0$ the differential part~\eqref{eq:B-recurrence-homo} of the system~\eqref{eq:B-recurrence-homo}--\eqref{eq:B-recurrence-rhs2}, then we arrive (for all~$k$) to the so-called degenerate
case of the hypergeometric equation. The solution of this equation can be written as
\begin{gather*}
C_k \frac{x+2}{(x-2)^3}+D_k \frac{8+(x+2)\ln x}{(x-2)^3},
\end{gather*}
where $C_k$ and $D_k$ are constants of integration. Clearly, in our case always $D_k=0$ while the constant $C_k$
depends on $k$ and should be chosen with the help of the condition $B_k(0)=0$.
Therefore, the main problem in construction of $B_k$ is to find a particular rational solution of
ODE~\eqref{eq:B-recurrence-homo}--\eqref{eq:B-recurrence-rhs2}, which exists by construction for all $k$.
Below we present the results for $k=1$ and $k=2$.

Substituting $B_0$ from equation~\eqref{eq:B0-solution} into equation~\eqref{eq:B1} and reducing both parts by factor $x(1-x/2)$ one finds
\begin{gather*}
x(1-x/2)B_1''(x)+ (1-5x/2 )B_1'(x)-2B_1(x)=-(1-x/2)(1-x/4),
\end{gather*}
where the primes denote derivatives with respect to $x$. Unique rational solution of this equation corresponding the initial condition $B_1(0)=0$ reads
\begin{gather}\label{eq:B1x}
B_1(x)=\frac{1}{64}x^2-\frac{11}{72}x+\frac{61}{144}+\frac{61}{36} \frac{x+2}{(x-2)^3}.
\end{gather}
The first terms of the Taylor expansion at $x=0$ are
\begin{gather*}
B_1(x)=-x-\frac{15}{16}x^2-\frac{61}{72}x^3-\frac{1525}{2304}x^4-\frac{61}{128}x^5-\frac{2989}{9216}x^6-\frac{61}{288}x^7+O\big(x^8\big).
\end{gather*}
Thus, we arrive at the following
\begin{Proposition}\label{prop:u2n1}
\begin{gather}\label{eq:u2n1}
u_2^1=-1,\qquad
u_4^1=-\frac{15}{16},\qquad
u_{2n}^1=-\frac{61}{12^2}\frac{(n+1)^2}{2^n}, \qquad n\geq3,
\end{gather}
where the sequence $u_{2n}^1$ is defined in equation~\eqref{eq:u{2n}(a)-taylor}.
\end{Proposition}
\begin{Corollary}\label{cor:p1n}
\begin{gather}\label{eq:p1n}
p_1(n)=\frac{(n+1)^2}{2^n}\left(\frac{1}{n+1}\sum\limits_{k=1}^n\frac{n_k}{k^2}-\frac{61}{144}\right)
\prod\limits_{k=1}^{n}k^{2n_k},\qquad n\geq3.
\end{gather}
\end{Corollary}
\begin{proof}
The proof is straightforward: expand equation~\eqref{eq:u2ka}, with the help of equation~\eqref{eq:polynom-P} and
compare with equation~\eqref{eq:u2n1} and take into account equation~\eqref{eq:p0n}.
\end{proof}
\begin{Remark}\label{rem:p1n}
Since $m(0)=m(1)=0$ (see Remark~\ref{rem:m(n)}), the numbers $p_1(1)$ and $p_1(2)$ are not defined.
According to Conjecture~\ref{con:expansion-u-main} $p_1(n)$ is a sequence of positive integers. It can be established directly with the help of equation~\eqref{eq:p0n}.
Consider the sum
\begin{gather*}
\sum\limits_{k=1}^{\infty}\frac{1}{(k+1)} \frac{1}{k^2}=\frac{\pi^2}{6}-1.
\end{gather*}
An elementary estimate shows that
\begin{gather*}
\sum\limits_{k=1}^{n}\frac{1}{(k+1)} \frac{1}{k^2}>\frac{\pi^2}{6}-1-\frac{1}{2n^2},\\
\sum\limits_{k=1}^{n}\frac{1}{(k+1)} \frac{1}{k^2}-\frac{n}{(n+1)^2}<\frac{1}{n+1}\sum\limits_{k=1}^n\frac{n_k}{k^2}
<\sum\limits_{k=1}^{n}\frac{1}{(k+1)} \frac{1}{k^2}.
\end{gather*}
Thus the sum in the parentheses in equation~\eqref{eq:p1n} is larger than
\begin{gather*}
\frac{\pi^2}{6}-1-\frac{61}{144}-\frac{1}{2n^2}-\frac{n}{(n+1)^2}>
\frac{(3+0.1)^2}{6}-1-\frac{61}{144}-\frac{n}{(n+1)^2}-\frac{1}{2n^2}\underset{n\geq11}{>}0.
\end{gather*}
Therefore, it is enough to check that $p_1(n)>0$ for $3\leq n\leq10$:
\begin{gather*}
12, \ 55, \ 12657, \ 176022, \ 84817044, \ 10913409936, \ 11716666225920, \ 509615533152000.
\end{gather*}
Although the numbers above are large the limiting value of the expression in the parentheses in equation~\eqref{eq:p1n}
is $0.2213\ldots$ and it is substantially smaller for $n\in[3,10]$.
\end{Remark}

Consider now the simplest application of equations~\eqref{eq:B-recurrence-homo}--\eqref{eq:B-recurrence-rhs2}, $k=2$,
\begin{gather}\label{eq:B2-gen}
\delta_x^2\left(\frac{B_2}{1+B_0}\right)-xB_2=\frac12\delta_x^2\left(\frac{B_1}{1+B_0}\right)^2-\frac{(1-B_0)B_1}{(1+B_0)^3}.
\end{gather}
Substituting into equation~\eqref{eq:B2-gen} $B_0$ and $B_1$ defined by equations~\eqref{eq:B0-solution} and \eqref{eq:B1x}, respectively,
and dividing both parts by $x(1-x/2)$ we arrive at the following ODE
\begin{gather}
x(1-x/2)B_2''(x)+(1-5x/2)B_2'(x)-2B_2(x)\nonumber\\
\qquad{}= -\frac{49}{2^{13}3^2}x^5+\frac{551}{2^{12}3^2}x^4-\frac{5455}{2^93^4}x^3+\frac{1471}{2^83^2}x^2-\frac{6313}{2^93^2}x+\frac{17015}{2^83^4}\nonumber\\
\qquad\quad {} -\frac{(61)^2}{2^53^4}\left(\frac{1}{(x-2)^3}+\frac{12}{(x-2)^4}+\frac{24}{(x-2)^5}\right).\label{eq:B2x}
\end{gather}
The unique rational solution of equation~\eqref{eq:B2x} satisfying initial condition $B_2(0)=0$ reads
\begin{gather*}
B_2(x) =\frac{1}{36864}x^5-\frac{263}{331776}x^4+\frac{1643}{172800}x^3-\frac{15923}{230400}x^2+\frac{41993}{172800}x-\frac{74849}{259200}\\
\hphantom{B_2(x) =}{} -\frac{2099}{4800(x-2)^2}+\frac{2272}{2025(x-2)^3}+\frac{3721}{432(x-2)^4}.
\end{gather*}
The first terms of the Taylor expansion of $B_2(x)$ are
\begin{gather*}
B_2(x) =x+\frac{63}{64}x^2+\frac{2917}{2592}x^3+\frac{335485}{331776}x^4+\frac{382273}{460800}x^5\\
\hphantom{B_2(x) =}{} +\frac{21009877}{33177600}x^6+\frac{105619}{230400}x^7+\frac{260899}{819200}x^8+\frac{1136621}{5308416}x^9+O\big(x^{10}\big).
\end{gather*}
Thus we arrive at the following
\begin{Proposition}
\begin{gather*}
 u_2^2=1,\qquad
u_4^2=\frac{63}{64},\qquad
u_6^2=\frac{2917}{2592},\qquad
u_8^2=\frac{335485}{331776},\qquad
u_{10}^2=\frac{382273}{460800},\\
 u_{2n}^2=\frac{61^2}{12^4}\frac{(n+1)^2}{2^{n+1}}\left(n+\frac{2^229^289}{5^261^2}\right),\qquad n=6,7, \qquad \ldots.
\end{gather*}
\end{Proposition}
\begin{Corollary}\label{cor:p2n}
\begin{gather}
p_2(5) =3345,\qquad{\rm for} \quad n\geq6 \label{eq:p2n}\\
p_2(n) =\frac{(n+1)^3}{2^{n+1}}\prod\limits_{k=1}^n k^{2n_k} \left( \left(
\frac{1}{n+1}\sum\limits_{k=1}^n\frac{n_k}{k^2}-\frac{61}{12^2}\right)^2 +
\frac{11\cdot73\cdot257}{5^212^4(n+1)}-\frac{1}{(n+1)^2}\sum\limits_{k=1}^n\frac{n_k}{k^4}\right).\nonumber
\end{gather}
\end{Corollary}
\begin{proof}Straightforward calculation, which is very analogous to the proof of Corollary~\ref{cor:p1n}.
\end{proof}

\begin{Remark}As follows from Remark~\ref{rem:m(n)} numbers $p_2(1)$, $p_2(2)$, $p_2(3)$, and $p_2(4)$ are not defined.
To prove that $p_2(n)$ for $n\geq6$ is an integer sequence it is enough to notice that the terms with the
denominator $k^4$ originated from the first sum in equation~\eqref{eq:p2n} cancel with that in the last sum in
this equation provided $n_k=1$. For $n_k\geq2$ these denominators cancel with the corresponding factors in
the product in front of the brackets.

The proof that all numbers $p_2(n)$ are positive also goes analogously to the proof of positiveness of $p_1(n)$
given in Remark~\ref{rem:p1n}. The only new formula needed in the estimate is
\begin{gather*}
\sum\limits_{k=1}^{\infty}\frac{1}{(k+1)} \frac{1}{k^4}=\frac{\pi^4}{90}+\frac{\pi^2}{6}-\zeta(3)-1=0.5252003\ldots.
\end{gather*}
On this way one can estimate that the number in the external parentheses in equation~\eqref{eq:p2n} is positive for all $n\geq7$.
The limiting value of the number in the external parentheses (after the product) in equation~\eqref{eq:p2n} is only $0.04898\ldots$,
while it is much smaller for small values of $n$, so the fact that $p_2(n)>0$ is not immediately obvious from the explicit
equation~\eqref{eq:p2n}. At the same time even the first members of the sequence $p_2(n)$, for $n\geq5$ are quite large:
\begin{gather*}
3345, \ 27825, \ 35168472, \ 4617359640, \ 7902853050240, \ 260852007650256, \ \ldots.
\end{gather*}
\end{Remark}

\section{Generating functions for the residues of coefficients}\label{sec:residues}
Consider now the decomposition of the coefficients $u_{2n}(a)$ in partial fractions
\begin{gather}\label{eq:u-n-parfrac}
u_{2n}=\sum\limits_{k=1}^n\sum\limits_{i=1}^{n_k}\frac{\gamma_{k,i}(n)}{\big(a^2+k^2\big)^i},
\end{gather}
where $\gamma_{k,i}(n)\in\mathbb R$ are some rational numbers. The number $\gamma_{k,i}(n)$ can be treated as residue of
the function $\big(a^2+k^2\big)^{i-1}u_{2n}(n)$, so, for brevity, we call all these numbers as residues.
The total number of the partial fractions in sum~\eqref{eq:u-n-parfrac} equals
\begin{gather*}
\sum\limits_{k=1}^{n}n_k=\sum\limits_{k=0}^{n}n_k-(n+1)=\sum\limits_{k=1}^{n+1} d(k)-(n+1),
\end{gather*}
where we use notation from the book by Hardy and Wright~\cite{HW}, $d(k)$, which denotes the number of divisors of integer $k$
including $1$ and $k$. According to Dirichlet (see~\cite{HW}),
\begin{gather}\label{eq:Dirichlet}
\sum\limits_{k=1}^{n}d(k)=n\ln(n)+(2\gamma-1)n+O\big(n^{\theta}\big),\qquad\theta=\frac12,
\end{gather}
where $\gamma=0.577215664\ldots$ is the Euler constant.
\begin{Remark}
We recall that finding the minimal value of $\theta$ in equation~\eqref{eq:Dirichlet} constitutes the so-called
Dirichlet divisor problem (see~\cite{W}). According to~\cite{W} Hardy and Landau in 1916 proved that
$\theta\geq\frac14$, while Huxley in 2003 proved that $\theta\leq\frac{131}{416}\approx0.31490$. My numerics
shows that the error for all $n\leq10^5$ does not exceed $2.3\cdot(n\ln(n))^{1/4}$.
\end{Remark}

Anyway, it follows from equation~\eqref{eq:Dirichlet} that the total number of coefficients $\gamma_{ki}$
in equation~\eqref{eq:u-n-parfrac} approximately equals to
\begin{gather*}
(n+1)\big(\ln(n+1)+2(\gamma-1)\big)+O\big(n^{\theta}\big).
\end{gather*}
We recall that the total number of coefficients $p_k(n)$, defining the numerator of $u_{2n}$ is $m(n)+1$
which is less than the number of residues by~$n$ (see equations~\eqref{eq:u2ka}--\eqref{eq:polynom-P}).
Obviously, the residues can be expressed in terms of $p_k(n)$ via linear relations
(see, e.g., Corollary~\ref{cor:gamma1n-pkn} below). In case we would calculate, say, all residues corresponding to
the poles of order higher than $1$ and one more residue for a pole of the first order we would be able to find
a~general formulae for~$p_k(n)$.

\begin{Remark}We can surely express linear combinations of the residues in terms of the gene\-ra\-ting functions $A_k(z)$ and $B_k(x)$.
Comparing equations~\eqref{eq:u-2n-Akn} and~\eqref{eq:u-n-parfrac} we can get $n-1$ linear relations for the
residues $\gamma_{k,i}(n)$ which are free of the numbers~$A_k[n]$, the simplest one is
\begin{gather*}
\sum\limits_{k=1}^n\gamma_{k,1}(n)=0,\qquad \text{for}\quad n>1.
\end{gather*}
In principle, we can get enough linear equations for the residues $\gamma_{k,i}(n)$ to express them as linear combinations
of numbers $A_l[n]$, however, we can explicitly calculate the latter numbers only for several first values of $l$. So that
to get explicit formulae for the residues with arbitrary large $n$ is problematic with this approach.

Analogously, we can use the generating functions $B_k(x)$ to get some explicit formulae for linear combinations of the
residues, e.g., expanding equation~\eqref{eq:u-n-parfrac} in the Taylor series at $a^2=0$ comparing this expansion
with \eqref{eq:u{2n}(a)-taylor} and using Propositions~\ref{prop:u-at0} and \ref{prop:u2n1} we find
\begin{gather*}
\sum\limits_{k=1}^n\sum\limits_{i=1}^{n_k}\frac{\gamma_{k,i}(n)}{k^{2i}}=\frac{n+1}{2^n}, \qquad n\geq1,\\
\sum\limits_{k=1}^n\sum\limits_{i=1}^{n_k}i\frac{\gamma_{k,i}(n)}{k^{2i+2}}=\frac{61}{12^2}\frac{(n+1)^2}{2^n},\qquad n\geq3.
\end{gather*}
\end{Remark}

Our goal is to calculate residues $\gamma_{k,i}(n)$ with the help of generating functions. We begin with the construction
of super generating function for sequences $\gamma_{1,i}(n)$, i.e., the case $k=1$. Define parameter
\begin{gather}\label{eq:xi1}
\xi_1=\sqrt{a^2+1}
\end{gather}
and the super generating function $V_1(\xi_1,z)$ as the rescaling of the solution $U(a,x)$: $V_1(\xi_1,z)=U(a,\xi_1 z)$.
Thus the function $V_1$ solves the following ODE
\begin{gather}\label{eq:V1}
\big(\delta_z^2+\xi_1^2-1\big)V_1=\xi_1 z(1+V_1)^3+(\delta_zV_1)^2-V_1\delta_z^2V_1,
\end{gather}
which is obtained by the rescaling $x=\xi_1 z$ of equation~\eqref{eq:trans-Garnier}. The function $V_1$ has the following
asymptotics
\begin{gather}\label{eq:V1xi}
V_1\underset{\xi_1\to 0}=\sum\limits_{k=-1}^\infty v_{1,k}(z)\xi_1^k.
\end{gather}
The fact that after the rescaling function $V_1(\xi,z)$ should coincide with $U(a,x)$ can be reformulated without any reference
to $U(a,x)$: generating functions $v_{1,k}(z)$ are singlevalued (in fact rational!) functions of $z$ and
\begin{gather}
v_{1,-1}(z)\underset{z\to0}=z+O\big(z^3\big),\qquad
v_{1,2n}(z)\underset{z\to0}=O\big(z^2\big),\nonumber\\
v_{1,2n+1}(z)\underset{z\to0}=O\big(z^3\big),\qquad n=0,1,\ldots.\label{eq:v1-initialdata}
\end{gather}
Now, substituting expansion~\eqref{eq:V1xi} into equation~\eqref{eq:V1} and successively equating coefficients of powers
$\xi_1^k$, for $k=-2$,$-1$,$0$,$1,\ldots$ to zero we, putting $v_k\equiv v_{1,k}$, find
\begin{gather}
v_{-1}\delta_{z}^2v_{-1}-(\delta_{z}v_{-1})^2-zv_{-1}^3=0,\label{eq:v-1}\\
v_{-1}\delta_{z}^2v_n-2\delta_{z}v_{-1}\delta_{z}v_n+\big(\delta_{z}^2v_{-1}-3zv_{-1}^2\big)v_n+F_n(z;v_{-1},\ldots,v_{n-1}), \qquad n=0,1,\ldots\label{eq:v-n}
\end{gather}
where equation~\eqref{eq:v-n} represents the recurrence relation with a function $F_n$ which depends on $z$ and $n+1$
variables $v_{-1}$, $\ldots$, $v_{n-1}$ determined on the previous steps
\begin{gather*}
F_0 =\delta_z^2v_{-1}-v_{-1}-3zv_{-1}^2,\\ %\label{eq:F0}\\
F_1 =(1+v_0)\delta_z^2v_0-v_0-(\delta_zv_0)^2-3zv_{-1}(1+v_0)^2,\\ %\label{eq:F1}\\
F_2 =(1+v_0)\delta_z^2v_1-2\delta_zv_0\delta_zv_1+\big(\delta_z^2v_0-6zv_{-1}(v_0+1)-1\big)v_1+v_{-1}-z(1+v_0)^3,%\label{eq:F2}
\end{gather*}
 and for $n\geq3$
\begin{gather*} F_n =(1+v_0)\delta_z^2v_{n-1}-2\delta_zv_0\delta_zv_{n-1}+\big(\delta_z^2v_0-6zv_{-1}(v_0+1)-1\big)v_{n-1}
+v_{n-3}-3zv_{n-2}\nonumber\\
\hphantom{F_n =}{} -z\sum_{\substack{i+j+k=n-2\\i,j,k\geq0}} v_iv_jv_k-3z\sum_{\substack{i+j=n-2\\i,j,k\geq0}}v_iv_j+
\sum_{\substack{i+j=n-1\\i,j\geq1}}v_i\delta_z^2v_j-\delta_zv_i\delta_zv_j-3zv_{-1}v_iv_j.%\label{eq:Fn}
\end{gather*}
The general solution of ODE \eqref{eq:v-1} reads
\begin{gather}\label{eq:v-1-solution}
v_{-1}=\frac{2C_2^2C_1z^{C_2-1}}{(1-C_1z^{C_2})^2},
\end{gather}
where $C_1$ and $C_2$ are constants of integration. Having in mind the first equation in conditions~\eqref{eq:v1-initialdata}
we find that $C_2=2$ and $C_1=1/8$, so that finally
\begin{gather}\label{eq:v-1-formula}
v_{-1}=\frac{z}{(1-z^2/8)^2}.
\end{gather}
In view of equation~\eqref{eq:v-1-formula} equation~\eqref{eq:v-n} is a special case of
the inhomogeneous Gauss hyper\-geo\-metric equation. The general solution of its homogeneous part can be written as follows
\begin{gather*}
{\tilde C}_1\frac{z\big(z^2+8\big)}{(1-z^2/8)^3}+{\tilde C}_2\frac{z\big(z^2\ln{z}+8\ln{z}+16\big)}{(1-z^2/8)^3},
\end{gather*}
where ${\tilde C}_1$ and ${\tilde C}_2$ are constants of integration. Since we are looking for the singlevalued solution
of equation~\eqref{eq:v-n} we can put ${\tilde C}_2=0$ and apply the method of variation of constants to ${\tilde C}_1$
and initial conditions~\eqref{eq:v1-initialdata} to find the unique generating function $v_n$. It is easy to prove,
inductively, that all functions $F_n$ are rational functions of $z$, so that all functions $v_n$ are also rational functions
of $z$. The first few functions $v_n$ are as follows (we return to the original notation):
\begin{gather}
v_{1,0}=\frac{z^2\big(24^2-16z^2+z^4\big)}{24^2\big(1-z^2/8\big)^3},\nonumber\\
v_{1,1}=\frac{z^3\big(4091904+285696z^2-1920z^4+80z^6-z^8\big)}{2592\big({-}8+z^2\big)^4},\nonumber\\
v_{1,2}=\frac{z^2\big(\alpha_0-\alpha_2z^2-\alpha_4z^4-\alpha_6z^6+198784z^8-3128z^{10}-25z^{12}\big)}{1166400\big({-}8+z^2\big)^5},\label{eq:v2-formula}\\
\alpha_0=12740198400, \qquad \alpha_2=6834585600, \qquad \alpha_4=946999296, \qquad \alpha_6=10810368,\nonumber\\
v_{1,3}=\frac{z^3\big({-}\beta_0+\beta_2z^2+\beta_4z^4+\beta_6z^6+\beta_8z^8-\beta_{10}z^{10}+\beta_{12}z^{12}-7365z^{14}+25z^{16}\big)}{2687385600\big({-}8+z^2\big)^6},\label{eq:v3-formula}\\
\beta_0=131784612249600, \qquad \beta_2=1890263236608, \qquad \beta_4=14900362739712,\nonumber\\
\beta_6=215484420096, \qquad \beta_8=5050358784, \qquad \beta_{10}=153847552, \qquad \beta_{12}=188436.\nonumber
\end{gather}
With the help of Maple code one can easily continue the list of the generating functions $v_{1,n}$ for $n>3$.

Now, consider application of the generating functions $v_{1,k}(z)$ to calculation of the resi\-dues~$\gamma_{1,i}(n)$. Consider
the Laurent expansion of the coefficients $u_{2n}(a)$:
\begin{gather*}%\label{eq:def-junior-gamma}
u_{2n}(a)=\sum\limits_{m=0}^{+\infty}\frac{\gamma_{1,n_1-m}(n)}{\xi_1^{2(n_1-m)}},
\end{gather*}
where $\xi_1$ is defined in equation~\eqref{eq:xi1} and the numbers $\gamma_{1,n_1-k}(n)$ for $k=0,\ldots, n_1-1$ coincide
with the residues in equation~\eqref{eq:u-n-parfrac}.
\begin{Proposition}\label{prop:small-gamma1k}
Put $n_1=k$, then for $k=1,2,\ldots$,
\begin{gather}\label{eq:gamma-n1}
\gamma_{1,k}(2k-1)=\frac{k}{8^{k-1}},\qquad
\gamma_{1,k}(2k)=\frac{(2k+1)^2}{9\cdot8^{k-1}},
\\
\label{eq:gamma-n1-1-small}
\gamma_{1,0}(1)=0,\qquad
\gamma_{1,0}(2)=-\frac13,\qquad
\gamma_{1,1}(3)=\frac{37}{96},\qquad
\gamma_{1,1}(4)=-\frac{17}{576},
\end{gather}
for $k=3,4\ldots$,
\begin{gather}
\gamma_{1,k-1}(2k-1)=\frac{(128k-3)(k+1)^2}{162\cdot8^{k}},\nonumber\\
\gamma_{1,k-1}(2k)=\frac{\big(3200k^2-6625k-582\big)(2k+1)^2}{109350\cdot8^{k-1}}.\label{eq:gamma-n1-1-general}
\end{gather}
\end{Proposition}
\begin{proof}
By the arguments analogous to those in Sections~\ref{sec:generating-A} and \ref{sec:generating-B} we prove that
$v_{1,k}$ are the generating functions for the residues $\gamma_{1,l}$, more precisely that means
\begin{gather*}
v_{1,-\epsilon+2m}=\sum\limits_{k=1}^\infty\gamma_{1,k-m}(2k-\epsilon)z^{2k-\epsilon},\qquad \epsilon=0,1,\qquad m=0,1,\ldots.
\end{gather*}
Developing explicit formulae~\eqref{eq:v-1-formula}--\eqref{eq:v2-formula} into the Taylor series we finish the proof.
\end{proof}
\begin{Remark}There is one more explicit formula~\eqref{eq:v3-formula} for the generating function $v_{1,3}$. Using it one proves,
\begin{gather*}
\gamma_{1,-1}(1)=0,\qquad
\gamma_{1,0}(3)=-\frac{431}{2304},\qquad
\gamma_{1,1}(5)=-\frac{62743}{552960},\qquad
\gamma_{1,2}(7)=-\frac{222359}{11059200},
\end{gather*}
and for $k\geq5$
\begin{gather*}
\gamma_{1,k-2}(2k-1)=\frac{\big(13107200k^3-41164800k^2-22621088k+3402171\big))(k+1)^2}{1968300\cdot8^{k+2}}.
\end{gather*}
The above formula together with the results obtained in Proposition~\ref{prop:small-gamma1k} allows one to make the following inductive
\begin{Conjecture}
\begin{gather*}
\gamma_{1,k-m}(2k-\epsilon)=\sum\limits_{l=0}^{2m-\epsilon}\alpha_l(\epsilon,m)k^l\frac{\big((2-\epsilon)k+1\big)^2}{8^k},
\qquad k\geq 2m+1,\qquad 2m\geq\epsilon,
\end{gather*}
where $\alpha_l(\epsilon,m)$ some rational numbers.
\end{Conjecture}
\end{Remark}

Explicit results for the residues allows us to get some consequences for the coefficients $p_k(n)$ (see equation~\eqref{eq:polynom-P}).
\begin{Corollary}\label{cor:gamma1n-pkn}
\begin{gather*}
\sum\limits_{k=0}^{m(n)}(-1)^kp_k(n)=\gamma_{1,n_1}(n)\prod\limits_{k=2}^{n}\big(k^2-1\big),\\
\sum\limits_{k=0}^{m(n)}(-1)^{k+1}kp_k(n)=\left(\gamma_{1,n_1-1}(n)+\gamma_{1,n_1}(n)\sum\limits_{k=2}^n\frac{n_k}{k^2-1}\right)
\prod\limits_{k=2}^{n}\big(k^2-1\big),
\end{gather*}
where $\gamma_{1,n_1-1}(n)$ and $\gamma_{1,n_1}(n)$ are given by equations~\eqref{eq:gamma-n1}, \eqref{eq:gamma-n1-1-small} and \eqref{eq:gamma-n1-1-general}, respectively.
\end{Corollary}

Now we describe a construction of the generating functions for $\gamma_{k,i}$ for the fixed $k>1$. For this purpose we introduce
an auxiliary parameter
\begin{gather*}
\xi_k=\big(a^2+k^2\big)^{\frac{1}{k+1}}.
\end{gather*}
In this case the super generating function $V_k(\xi_k,z)$ is defined as the rescaling of the solution $U(a,x)$:
$V_k(\xi_k,z)=U\big(\sqrt{\xi_k^{k+1}-k^2},\xi_k z\big)$.
Therefore, the function $V_k$ solves the following ODE
\begin{gather}\label{eq:Vk}
\big(\delta_z^2+\xi_k^{k+1}-k^2\big)V_k=\xi_k z(1+V_k)^3+(\delta_zV_k)^2-V_k\delta_z^2V_k.
\end{gather}
The super generating function $V_k$ has the following asymptotics
\begin{gather}\label{eq:Vkxi}
V_k\underset{\xi_k\to 0}=\sum\limits_{l=-1}^\infty v_{k,l}(z)\xi_k^l,
\end{gather}
which has the same form as for the function $V_1$ (cf.\ equation~\eqref{eq:V1xi}). However, the generating functions
$v_{k,l}(z)$ differ with $v_{1,l}(z)$. Comparing equations~\eqref{eq:V1} and \eqref{eq:Vk} it is obvious that the
functions $v_{k,n}(z)$ satisfy the same system of equations~(\eqref{eq:v-1} and~\eqref{eq:v-n})
as the functions $v_n$ but with some minor change of the inhomogeneous contribution, the function $F_n$:
\begin{gather}
F_0 =\delta_z^2v_{-1}-k^2v_{-1}-3zv_{-1}^2,\label{eq:F0k}\\
F_1 =(1+v_0)\delta_z^2v_0-k^2v_0-(\delta_zv_0)^2-3zv_{-1}(1+v_0)^2,\label{eq:F1k}
\end{gather}
and for $n\geq2$
\begin{gather}
F_n =(1+v_0)\delta_z^2v_{n-1}\! -2\delta_zv_0\delta_zv_{n-1}\! + \big(\delta_z^2v_0-6zv_{-1}(v_0+1)-k^2\big)v_{n-1}\!
+v_{n-k-2}\! -3zv_{n-2}\nonumber\\
\hphantom{F_n =}{}+\sum_{\substack{i+j+m=n-2\\i,j,m\geq0}} v_iv_jv_m-3z\sum_{\substack{i+j=n-2\\i,j,k\geq0}}v_iv_j+
\sum_{\substack{i+j=n-1\\i,j\geq1}}v_i\delta_z^2v_j-\delta_zv_i\delta_zv_j-3zv_{-1}v_iv_j,\!\!\!\label{eq:Fnk}
\end{gather}
where for $n=2,\ldots, k$ we put formally $v_{n-k-2}\equiv0$.

Now we present explicit formulae which shows that functions $v_{k,l}(z)$ generate the resi\-dues $\gamma_{k,i}(n)$.
As in the case $k=1$, it is convenient to generalize residues $\gamma_{k,i}$ and define them for all integer $i\leq n_k$
as the coefficients of the Laurent expansion at $\xi_k=0$:
\begin{gather}\label{eq:def-junior-gamma-k}
u_{2n}(a)=\sum\limits_{m=0}^{+\infty}\frac{\gamma_{k,n_k-m}(n)}{\xi_k^{(k+1)(n_k-m)}}.
\end{gather}

Define nonnegative integers, $p$ and $q\leq k$ such that $n=p(k+1)+q$, then $n_k=p$ for $q<k$ and $n_k=p+1$ for $q=k$.
For each nonnegative integer $i$ and $q=0,1,\ldots,k-1$ define $l=i(k+1)+q$, then
\begin{gather}
v_{k,l}(z) =\sum\limits_{p=0}^\infty\gamma_{k,p-i}\big(p(k+1)+q\big)z^{p(k+1)+q},\label{eq:v-kl-q<k}\\
\mathrm{for}\quad
 q=0,1,\ldots,k-1,\quad
i=0,1,\ldots,\quad
l=i(k+1)+q;\label{eq:l-for-q<k}\\
v_{k,l}(z) =\sum\limits_{p=0}^\infty\gamma_{k,p+1-i}\big(p(k+1)+q\big)z^{p(k+1)+q},\label{eq:v-kl-q=k}\\
\mathrm{for}\quad
 q=k,\quad
i=0,1,\ldots,\quad
l=i(k+1)-1.\label{eq:l-for-q=k}
\end{gather}

As an example consider the case $k=2$. As is explained above the functions: $v_{2,-1}$, $v_{2,0}$, and $v_{2,1}$, are
defined by equations~\eqref{eq:v-1}--\eqref{eq:v-n}, which formally coincide with the equations for the functions
$v_{-1}$, $v_0$, and $v_1$, respectively. However, the first set of functions are different from the second one.
The functions $v_{2,-1}$ and $v_{-1}$ are different because of the initial conditions. The functions $v_{2,0}$ and $v_{2,1}$
differ from the corresponding functions $v_0$ and $v_1$ since the inhomogeneous terms $F_0$ and $F_1$ after substitution
the functions $v_{2,-1}$ instead of $v_{-1}$ and then $v_{2,0}$ instead of $v_0$ differ with the case $k=1$.

The function $v_{2,-1}$ is given by equation~\eqref{eq:v-1-solution} where we have to choose properly the constants
of integration: $C_1$ and $C_2$. In this case, $C_2=3$ and $C_1=-1/18$. These constants are defined from the fact that
for the first time the factor $a^2+4$ appears in $u_4$, see Remark~\ref{rem:coeff3}, therefore $C_2-1=4/2$, and the Laurent
expansion of $u_4$ reads
\begin{gather*}
u_4=-\frac{1}{\xi_2}-\frac13-\frac19\xi_2-\cdots,
\end{gather*}
the coefficient of the leading term is $-1$,
which means (see equation~\eqref{eq:v-1-solution}) that we have to put $2C_2^2C_1=-1$. Thus, we get
\begin{gather*}
v_{2,-1}(z)=-\frac{z^2}{(1+z^3/18)^2}=\sum_{p=0}^\infty(-1)^{p+1}\frac{p+1}{18^p}z^{3p+2},
\end{gather*}
therefore (equations~\eqref{eq:v-kl-q=k} and \eqref{eq:l-for-q=k})
\begin{gather}\label{eq:gamma2p+1-general}
\gamma_{2,p+1}(3p+2)=(-1)^{p+1}\frac{p+1}{18^p}, \qquad p=0,1,\ldots.
\end{gather}
Note that for $k=2$ and $n=3p+2$ we have $n_k=p+1$. To calculate $\gamma_{2,n_2}$ for $n=3p$ and $n=3p+1$ we have to find
the generating functions $v_{2,0}$ and $v_{2,1}$, respectively. To find them we have to solve successively two linear
inhomogeneous ODEs of the second order
\begin{gather}
v_{2,-1}\delta_{z}^2v_{2,0}-2\delta_{z}v_{2,-1}\delta_{z}v_{2,0}+\big(\delta_{z}^2v_{2,-1}-3zv_{2,-1}^2\big)v_{2,0}\nonumber\\
\qquad{}+\delta_z^2v_{2,-1} -4v_{2,-1}-3zv_{2,-1}^2=0,\label{eq:F02}\\
v_{2,-1}\delta_{z}^2v_{2,1}-2\delta_{z}v_{2,-1}\delta_{z}v_{2,1}+\big(\delta_{z}^2v_{2,-1}-3zv_{2,-1}^2\big)v_{2,1}
+(1+v_{2,0})\delta_z^2v_{2,0}\nonumber\\
\qquad{}-4v_{2,0}-(\delta_zv_{2,0})^2-3zv_{2,-1}(1+v_{2,0})^2=0.\label{eq:F12}
\end{gather}
The solution of equation~\eqref{eq:F02}, $v_{2,0}$, should be a rational function of $z^3$. This condition uniquely determines our function
\begin{gather*}
v_{2,0}=-\frac{3z^3\big(5184-18z^3+z^6\big)}{4(18+z^3)^3}.
\end{gather*}
Taylor expansion at $z=0$ reads
\begin{gather*}
v_{2,0}=\sum_{p=1}^\infty(-1)^p\frac{3(3p+1)^2}{4\cdot18^p}z^{3p}.
\end{gather*}
Comparing the last equation with equations~\eqref{eq:v-kl-q<k} and \eqref{eq:l-for-q<k} for $i=0$ and $q=0$ we obtain
\begin{gather}\label{eq:gamma2p-general}
\gamma_{2,p}(3p)=(-1)^p \frac{3(3p+1)^2}{4\cdot18^p},\qquad p=1,2,\ldots.
\end{gather}
The suitable solution of equation~\eqref{eq:F12} is uniquely defined by the condition that it is a rational
function of $z^3$ multiplied by $z$. Explicitly, it reads
\begin{gather*}
v_{2,1}=\frac{z\big(25z^{15}+4464z^{12}-468180z^9+284788224z^6-2687385600z^3-6046617600\big)}{172800\big(18+z^3\big)^4}.
\end{gather*}
Series expansion at $z=0$ reads
\begin{gather}\label{eq:v21-series}
v_{2,1}=-\frac13z-\frac2{27}z^4+\sum_{p=2}^\infty(-1)^p\frac{9(50p-31)(3p+2)^2}{3200\cdot18^p}z^{3p+1}.
\end{gather}
Comparing equation~\eqref{eq:v21-series} with equations~\eqref{eq:v-kl-q<k} and \eqref{eq:l-for-q<k}
for $i=0$ and $q=1$ we find
\begin{gather}
\gamma_{2,0}(1)=-\frac13,\label{eq:gamma201}\\
\gamma_{2,1}(4)=-\frac2{27},\qquad
\gamma_{2,p}(3p+1)=(-1)^p\frac{9(50p-31)(3p+2)^2}{3200\cdot18^p},\qquad
p=2,3,\ldots.\label{eq:gamma2p3p+1-general}
\end{gather}
The coefficient $\gamma_{2,0}(1)$ (equation~\eqref{eq:gamma201}) is nothing but the first coefficient of the
Taylor expansion of $u_2=1/\big(\xi_2^3-3\big)$ at $\xi_2=0$, see equation~\eqref{eq:def-junior-gamma-k} and explicit
formula for $u_2$ in Remark~\ref{rem:coeff3}. Coefficients $\gamma_{2,p}(3p+1)$ for $p=1,2,\ldots$, are
senior residues $\gamma_{2,n_2}(n)$ for $n=3p+1$. Thus,
equations~\eqref{eq:gamma2p+1-general}, \eqref{eq:gamma2p-general}, and~\eqref{eq:gamma2p3p+1-general}
deliver explicit expressions for senior residues $\gamma_{2,n_2}(n)$ for all~$n$. One can surely continue to
calculate the functions $v_{2,l}$ for $l=2,3,\ldots,$ and obtain explicitly general formulae for the junior
residues at any fixed distance from the senior ones.

The scheme for computation of $\gamma_{k,i}(n)$ with $k=1$ and $2$ presented above is working for any fixed
integers $k>0$ and $i<n_k$. For a fixed value of $k$ to find the senior residue $\gamma_{k,n_k}(n)$ for all $n$
one have to find the tuple of $k+1$ functions, $\{v_{k,-1}(z),v_{k,0}(z),\ldots,v_{k,k-1}(z)\}$. For junior
coefficients $\gamma_{k,n_k-i}(n)$ one have to find the $(k+1)$-tuple of functions with the second subscript
shifted by $+i(k+1)$. For the particular values of $k$ and $i$ this is a bit tedious but a straightforward
procedure related with the successive solution of linear second order inhomogeneous ODEs.
The homogeneous part of the ODEs is the degenerate hypergeometric equation. It follows from the fact that
the coefficients of this equation are defined by the function $v_{k,-1}(z)$ which is a proper solution of
the universal ODE \eqref{eq:v-1} and therefore, as follows from equation~\eqref{eq:v-1-solution}, reads
\begin{gather}\label{eq:vk-1-solution}
v_{k,-1}(z)=\frac{2(k+1)^2C_{1,k}z^{k}}{\big(1-C_{1,k}z^{k+1}\big)^2},
\end{gather}
where we put $C_2=k+1$ and $C_1\to C_{1,k}$ since the constant of integration depends on $k$. As follows from
equations~\eqref{eq:F0k}--\eqref{eq:Fnk} the inhomogeneous part of the linear ODEs is a rational function of $z$.
Solutions of such equations can be found by the standard procedure of variation of constants of integration,
however for large values of $k$ the problem becomes tedious. Within this approach to find a formula which
would be valid for all $k$ and/or $i$ seems to be a complicated problem. It is not trivial even to find a
general formula for numbers $C_{1,k}$ in equation~\eqref{eq:vk-1-solution}, which is an important step towards
the general formula for generating functions $v_{k,l}$. Below we present first terms of the sequence $C_{1,k}$:
\begin{gather*}
C_{1,1}=\frac18,\qquad
C_{1,2}=-\frac1{18},\qquad
C_{1,3}=\frac9{1024},\qquad
C_{1,4}=-\frac1{1350},\qquad
C_{1,5}=\frac{625}{15925248},\\
C_{1,6}=-\frac9{6272000},\qquad
C_{1,7}=\frac{117649}{3057647616000},\qquad
C_{1,8}=-\frac{2}{2531725875},\qquad \ldots,
\end{gather*}
which allows us to make the following
\begin{Conjecture}\label{con:C1k}
\begin{gather*}
C_{1,k}=\frac{(-k)^{k-1}}{2^{k}(k+1)^2\big((k-1)!\big)^3},\qquad
k=1,2,\ldots.
\end{gather*}
\end{Conjecture}
Assuming Conjecture~\ref{con:C1k} is true we can launch the iterative process
(see equations~\eqref{eq:v-n}, \eqref{eq:F0k}--\eqref{eq:Fnk}) of calculation the
functions $v_{k,l}$, for $q=0,1,\ldots$. For example, the direct consequence of
Conjecture~\ref{con:C1k} is the following
\begin{Conjecture}\label{con:vk0}
For $k=1,2,\ldots$,
\begin{gather}\label{eq:vk0-solution}
v_{k,0}(z)=\frac{z^{k+1}\big(a_kz^{2k+2}+b_kz^{k+1}+c_k\big)}{\big(1-C_{1,k}z^{k+1}\big)^3},
\end{gather}
where $C_{1,k}$ is given in Conjecture~{\rm \ref{con:C1k}} and
\begin{gather*}
a_k =\frac{2^{2-3k}(-k)^{3k-3}}{(k+2)^2(k+1)^5\big((k-1)!\big)^9},\\
b_k =\frac{2^{2-2k}(-k)^{2k-2}(k^2-3)}{(k+2)^2(k+1)^3\big((k-1)!\big)^6},\\
c_k =\frac{2^{2-k}(-k)^{k-1}}{(k+1)\big((k-1)!\big)^3}.
\end{gather*}
\end{Conjecture}
Based on the example for $k=2$ it is not complicated to calculate a few more functions $v_{k,q}$, however,
to get, say, formulae for senior residues $\gamma_{n_k}(n)$ for all $n$ we have to find functions $v_{k,q}$
for $q=1,\ldots,k-1$.

Now we present the formulae for the residues that can be obtained with the help of
generating functions \eqref{eq:vk-1-solution} and \eqref{eq:vk0-solution}.
\begin{Conjecture}For $k=1,2,\ldots$ and $p=1,2,\ldots$
\begin{gather*}
\gamma_{k,p+1}\big(p(k+1)+k\big)=\frac{(p+1)(-k)^{(p+1)(k-1)}}{2^{(p+1)k-1}(k+1)^{2p}\big((k-1)!\big)^{3(p+1)}},\\
\gamma_{k,p}\big(p(k+1)\big)=\frac{(pk+p+1)^2(-k)^{p(k-1)}}{2^{pk-2}(k+2)^2(k+1)^{2p-1}\big((k-1)!\big)^{3p}}.
\end{gather*}
\end{Conjecture}
We finish this section by the following
\begin{Remark}
Explicit construction for generating functions $v_{1,k}(z)$ and $v_{2,k}(z)$ in fact provide us a proof
that $n_1=\left[(n+1)/2\right]$ and $n_2=\left[(n+1)/3\right]$, respectively. This justify
Conjecture~\ref{con:expansion-u-main} in its part concerning numbers $n_k$ for $k=1$ and $2$
(see the first equation~\eqref{eq:n-k-m(n)}). To make analogous proof for general $k$ it is enough to
prove existence of rational generating func\-tions~$v_{k,q}(z)$ for $q=-1,0,\ldots,k-1$ satisfying suitable
initial condition at $z=0$.
\end{Remark}

\section[Polynomials $P_{m(n)}(x)$]{Polynomials $\boldsymbol{P_{m(n)}(x)}$}\label{sec:polynom-conjectures}
In this section I will not write any proofs therefore all statements are formulated as conjectures.

Another interesting property of the coefficients $u_{2n}$ is the greatest common divisor (g.c.d.) of coefficients of
polynomials $P_{m(n)}(x)$. We recall that these coefficients are positive integers.
\begin{Conjecture}\label{con:gen}
\begin{gather}\label{eq:gcd}
{\rm g.c.d.}\left\{p_{m(n)}(n), p_{m(n)-1}(n),\ldots,p_1(n), p_0(n)\right\}=
\begin{cases}
(n+1) 3^{z_n},& {\rm iff} \ n+1 \ \text{is odd},\\
\dfrac{n+1}2 3^{z_n},& {\rm iff} \ n+1 \ \text{is even},
\end{cases}
\end{gather}
where $z_n$ is a nonnegative integer sequence.
\end{Conjecture}
Our goal is to define the sequence $z_n$. First we define the subsequence of zeroes, i.e., those $n=a_k$,
$k=1,2,\ldots$ for which $z_{a_k}=0$.
Consider the triangular decomposition of $n$, namely, the pair of positive integer numbers $(q,l)$,
where $q$ is the triangular floor and $l$ is the triangular reminder
\begin{gather*}
n=\frac{q(q+1)}2+l,\qquad0\leq l\leq n,
\end{gather*}
where
\begin{gather*}
q=\max\left\{\hat{q}\in\mathrm Z_{+}\colon \frac{\hat{q}(\hat{q}+1)}2\leq n\right\},
\end{gather*}
i.e., $q$ define the largest triangular number which is not larger than $n$. Clearly, for any given $n\geq1$ the triangular decomposition is uniquely defined. The first triangular decompositions are
\begin{gather*}
1=(1,0),\quad
2=(1,1),\quad
3=(2,0),\quad
4=(2,1),\quad
5=(2,2),\quad
6=(3,0),\quad \ldots.
\end{gather*}
Now for any given $k=1,2,\ldots$ with the triangular decomposition $(q,l)$ we define
\begin{gather}\label{eq:a-k}
a_k=\frac{3^q+3^l}2-1.
\end{gather}
The first members of the sequence $a_k$ are as follows
\begin{gather*}
a_1=1,\quad
a_2=2,\quad
a_3=4,\quad
a_4=5,\quad
a_5=8,\quad
a_6=13,\quad
a_7=14,\\
a_8=17,\quad
a_9=26,\quad\ldots.
\end{gather*}
The numbers $a_k$ have a simple presentation in the ternary (base-3) numerical system, the numbers whose digits in this system
are in nondecreasing order
\begin{gather*}
\underset{q}{\underbrace{1\ldots1\overset{l}{\overbrace{2\ldots2}}}}.
\end{gather*}
This sequence can be found in \cite{OEIS} as A023745 and called ``plaindromies'', do not mix with palindromies!
\begin{Conjecture}\label{con:zeroes}
All zeroes of the sequence $z_n$ are enumerated by the monotonically increasing sequence~$a_k$, i.e.,
\begin{gather*}
z_{a_k}=0,\qquad{\rm and} \qquad z_n\neq0 \quad {\rm if} \quad n\neq a_k.
\end{gather*}
\end{Conjecture}
So, for all other values of $n$ our g.c.d.~\eqref{eq:gcd} is divisible by~3. In particular,
\begin{gather*}
z_3=z_6=z_7=z_9=z_{10}=z_{11}=1.
\end{gather*}
Moreover, from time to time appear the higher powers of $3$. The second natural question is at what $n$ happens
the first occurrence of the factor $3^k$, for $k=1,2,\ldots$ in g.c.d.~\eqref{eq:gcd}? \\
To give the answer on this question we define the sequence
\begin{gather}\label{eq:b-k}
b_k=\frac32\big(3^k-1\big).
\end{gather}
The numbers $b_k$ have simple presentation in base-3 numerical system
\begin{gather*}
\underset{k}{\underbrace{1\ldots1}0}.
\end{gather*}
The first members of the sequence
\begin{gather*}
b_1=3,\quad
b_2=12,\quad
b_3=39,\quad
b_4=120,\quad
b_5=363,\quad\ldots.
\end{gather*}
This sequence can be found in \cite{OEIS} as A029858 and A031988.
\begin{Conjecture}\label{con:b-k}
\begin{gather*}
z_n<k\quad{\rm for}\quad n<b_k,\qquad z_{b_k}=k.
\end{gather*}
\end{Conjecture}
Clearly, to completely define the sequence $z_n$ it is enough to describe for every $k=1,2,\ldots$ all solutions of equation
$z_n=k$ for $n>b_k$. This, however, appear to be a complicated problem. We begin with the description of solutions
$z_n=k$ for $b_k<n<b_{k+1}$.
\begin{Conjecture}\label{con:b-k-k+1}
For every positive integer $k=1,2,3,\ldots$, there are exactly $\frac{(k+2)(k+3)}2$ solutions of equation $z_n=k$ for $n<b_{k+1}$.
These solutions are given by numbers
\begin{gather}\label{eq:p-nforn<b-k+1}
b_k^{(k-m)(k+m+5)/2+l}=\frac32\big(3^{k+1}-3^{m}-3^{m-l}-1\big),
\end{gather}
where $m=k,k-1,\ldots,0,-1$ and $l=0,\ldots,m+1$.
\end{Conjecture}
The first successive members of the sequence~\eqref{eq:p-nforn<b-k+1} for $m=k\leq1$, $l=0$, $l=1$, $l=2$:
\begin{gather*}
b_k^0=\frac32\big(3^k-1\big)=b_k,\qquad
b_k^1=b_k+3^k,\qquad
b_k^2=b_k^1+3^{k-1},
\end{gather*}
The last successive members of the sequence~\eqref{eq:p-nforn<b-k+1} for $m=0$, $l=0$, $l=1$; $m=-1$, $l=0$:
\begin{gather*}
b_k^{k(k+5)/2}=\frac32\big(3^{k+1}-1\big)-3=b_{k+1}-3,\qquad
b_k^{k(k+5)/2+1}=b_{k+1}-2,\\
b_k^{(k+1)(k+4)/2}=b_{k+1}-1.
\end{gather*}
The total number of terms in sequence~\eqref{eq:p-nforn<b-k+1} is
\begin{gather}\label{eq:finite-b-k}
\frac{(k+1)(k+4)}{2}+1=\frac{(k+2)(k+3)}{2},
\end{gather}
we add $1$ because we start enumeration with $0$. So, it is the sequence of the triangular numbers (A000217 in~\cite{OEIS})
without the first two terms.
The formula~\eqref{eq:p-nforn<b-k+1} reflects a simple recurrence construction of the sequence $b_k^{(k-m)(k+m+5)/2+l}$.
This construction can be described as follows: For all $m=k,\ldots, 1$ (not $-1$ as above!) consider $(m+1)\times(m+1)$
unit matrix $I_{m+1}$ we will treat the rows of this matrix as the numbers written in the base-3 numerical system and also
one more number
\begin{gather*}
1\underset{m-1}{\underbrace{0\ldots0}}1/2=\underset{m-1}{\underbrace{1\ldots1}}2.
\end{gather*}
Now, take $b_k^0$ and successively add $k+2$ numbers defined above for $m=k$. Then, put $m=k-1$ and add $k+1$ corresponding
numbers, until we arrive to $m=1$ where we have to add three numbers: $10$, $1$, $2$ (base-3!). Finally, we formally consider $m=0$ to which
we associate two numbers, both equal to $1$. All in all, counting together with $b_k^0$, we get the finite sequence of
size~\eqref{eq:finite-b-k}, which coincides with sequence~\eqref{eq:p-nforn<b-k+1}.

To find solutions of equation $z_n=k$ for $n>b_{k+1}$ looks a complicated problem. Instead we present the answer in
a geometric form. We consider the plot of the function $n\rightarrow z_n$ on the $(x,y)$-plane and connect by
the segments the neighbouring points $(n, z_n)$ and $(n+1, z_{n+1})$. Together with the $x$-axis we get a figure
that we call the {\it fence}. We are going to describe
how one can build this fence. We, actually, present two equivalent constructions. For the first one we need to define two shapes,
$\mathcal{A}$ and $\mathcal{B}$:
\begin{center}
\begin{picture}(420,80)
\put(10,30){\circle*{2}}
\put(10,30){\line(1,1){10}}
\put(20,40){\circle*{2}}
\put(20,40){\line(1,0){10}}
\put(30,40){\line(1,0){10}}
\put(30,40){\circle*{2}}
\put(40,40){\line(1,1){10}}
\put(40,40){\circle*{2}}
\put(50,50){\line(1,-1){10}}
\put(50,50){\circle*{2}}
\put(60,40){\line(1,0){10}}
\put(60,40){\circle*{2}}
\put(70,40){\line(1,1){10}}
\put(70,40){\circle*{2}}
\put(80,50){\line(1,0){10}}
\put(80,50){\circle*{2}}
\put(90,50){\line(1,-1){10}}
\put(90,50){\circle*{2}}
\put(100,40){\line(1,1){10}}
\put(110,50){\circle*{2}}
\put(110,50){\line(1,0){10}}
\put(120,50){\circle*{2}}
\put(100,40){\circle*{2}}
\put(120,50){\line(1,0){10}}
\put(130,50){\circle*{2}}
\put(130,50){\line(1,1){10}}
\put(140,60){\circle*{2}}
\put(160,30){Shape $\mathcal{A}$ is the tuple of $14$ points with the respective}
\put(5,10){heights: $z,z+ 1,z+ 1,z+ 1,z+ 2,z+ 1,z+ 1,z+ 2,z+ 2,z+ 1,z+ 2,z+ 2,z+ 2,z+ 3.$}
\end{picture}

\begin{picture}(420,80)
\put(20,40){\line(1,0){10}}
\put(20,40){\circle*{2}}
\put(40,50){\line(1,0){10}}
\put(40,50){\circle*{2}}
\put(30,40){\line(1,1){10}}
\put(30,40){\circle*{2}}
\put(50,50){\line(1,-1){10}}
\put(50,50){\circle*{2}}
\put(60,40){\line(1,1){10}}
\put(60,40){\circle*{2}}
\put(70,50){\line(1,0){10}}
\put(70,50){\circle*{2}}
\put(80,50){\line(1,0){10}}
\put(80,50){\circle*{2}}
\put(90,50){\line(1,1){10}}
\put(90,50){\circle*{2}}
\put(100,60){\line(1,-1){10}}
\put(100,60){\circle*{2}}
\put(110,50){\line(1,0){10}}
\put(110,50){\circle*{2}}
\put(120,50){\line(1,1){10}}
\put(120,50){\circle*{2}}
\put(130,60){\line(1,0){10}}
\put(130,60){\circle*{2}}
\put(140,60){\circle*{2}}
\put(150,40){Shape $\mathcal{B}$ is the tuple of $13$ points with the respective}
\put(20,20){heights: $z,z,z+1,z+1,z,z+1,z+1,z+1,z+2,z+1,z+1,z+2,z+2.$}
\end{picture}
\end{center}

These shapes define the upper edge of our fence. As long as we know coordinates, $(n,z_n)$, of any point of the shapes
we immediately now the coordinates of all their other points, by using the scheme presented on the corresponding figures.
Practically, it means that we put this point of the shape into the right position and orient the shape such that the points
with equal $y$-coordinates (heights) would be parallel to the $x$-axis.

Since our fence is semi-infinite to the right direction we present inductive construction starting from $n=1$ to the right side.
The very first step is irregular, we take shape $\mathcal{A}$, cut the first two segments and attach its third point (which
after the cutting becomes the first one) to the point with the coordinates $(1,0)$. Then the end point of the shape will have
the coordinates $(12,2)$. After that fall down by 2 units to the point $(13,0)$ and we attach the left point of shape $\mathcal{B}$
to the last point. Now the last point of the shape $\mathcal{B}$ is $(25,2)$ and again we get a fall down by two units at the
point $(26,0)$. We attach the left point of shape $\mathcal{A}$ to the point $(26,0)$. The last point of shape $\mathcal{A}$
has he coordinates $(39,3)$ and we have fall down by three units at point $(40,0)$ and attach to this point the left end of
shape $\mathcal{B}$, and so on. We can present the construction of this fence as the following symbolic sequence:
\begin{gather}\label{eq:ABsequence}
\mathcal{A}' 2 \mathcal{B} 2 \mathcal{A} 3^+\mathcal{B} 2 \mathcal{A} 2 \mathcal{B} 3^+
\mathcal{A} 2 \mathcal{B} 2 \mathcal{A}3^+\mathcal{B} 2 \mathcal{A} 2 \mathcal{B} 3^+\ldots.
\end{gather}
The prime in the first symbol $\mathcal{A}'$ clearly denotes the cutting procedure explained above. The symbol $3^+$
may denote any integer $\geq3$. In the above sequence numbers $2$ and $3^+$ denotes the points with the coordinates
$z_{n+1}=z_n-2$ and $z_{n+1}=z_n-3^+$, respectively. To finally define sequence~\eqref{eq:ABsequence} we have to know at what
points $3^+>3$ may happen? Suppose $n$ is the $x$-coordinate of the right endpoint of the shapes $\mathcal{A}$ or
$\mathcal{B}$, from these points we suffer the falls.
We call the point $n+1$ resonant if it coincides with one of the members of the sequence $a_k$ (see equation~\eqref{eq:a-k}),
i.e., $n+1=a_k$ for some $k$. If the point is nonresonant, then $3^+=3$. In the resonant case $3^+=z_n$, which means
that we fall down on the $x$-axis. Note that at the places, where $3^+$ are located in sequence~\eqref{eq:ABsequence}
always $z_n\geq3$. The resonances sometimes happen at those points where we have the fall down by two units.
It means that at this points $z_n=2$ and we fall down on the $x$-axis. In the nonresonant cases we have a drop down by $2$ or $3$ units according
sequence~\eqref{eq:ABsequence} but we still remain higher the $x$-axis.

We call the fence constructed, as explained above, accordingly symbolic sequence~\eqref{eq:ABsequence} the
{\it quasiperiodic fence} $\mathcal{P}$.
\begin{Conjecture}\label{con:main}
The heights of quasiperiodic fence $\mathcal{P}$ at positive integers $n=1,2,\ldots$ coincide with the sequence $ z_n$
defined in Conjecture~{\rm \ref{con:gen}}.
\end{Conjecture}

We can define fence $\mathcal{P}$ in a different way. Consider a new shape $\mathcal{C}$, which is obtained from shape~$\mathcal{A}$ by cutting off two segments from each end:
\begin{center}
\begin{picture}(420,80)
\put(40,40){\line(1,0){10}}
\put(40,40){\circle*{2}}
\put(50,40){\line(1,1){10}}
\put(50,40){\circle*{2}}
\put(60,50){\vector(1,-1){6}}
\put(60,50){\line(1,-1){10}}
\put(60,50){\circle*{2}}
\put(70,40){\line(1,0){10}}
\put(70,40){\circle*{2}}
\put(80,40){\line(1,1){10}}
\put(80,40){\circle*{2}}
\put(90,50){\line(1,0){10}}
\put(90,50){\circle*{2}}
\put(100,50){\vector(1,-1){6}}
\put(100,50){\line(1,-1){10}}
\put(100,50){\circle*{2}}
\put(110,40){\line(1,1){10}}
\put(120,50){\circle*{2}}
\put(120,50){\line(1,0){10}}
\put(130,50){\circle*{2}}
\put(110,40){\circle*{2}}
\put(180,40){Shape $\mathcal{C}$ is the tuple of $10$ points}
\put(10,20){with the respective heights: $z,z,z+1,z,z,z+1,z+1,z,z+1,z+1$.}
\end{picture}
\end{center}

On two descending sides of the triangles in shape $\mathcal{C}$ we put arrows which means nothing but the direction
of their deformations:

\begin{center}
\begin{picture}(420,90)
\put(20,50){\line(1,0){10}}
\put(20,50){\circle*{2}}
\put(30,50){\line(1,1){10}}
\put(30,50){\circle*{2}}
\put(40,60){\vector(1,-2){6}}
\put(40,60){\line(1,-2){10}}
\put(40,60){\circle*{2}}
\put(50,40){\line(1,0){10}}
\put(50,40){\circle*{2}}
\put(60,40){\line(1,1){10}}
\put(60,40){\circle*{2}}
\put(70,50){\line(1,0){10}}
\put(70,50){\circle*{2}}
\put(80,50){\vector(1,-1){6}}
\put(80,50){\line(1,-1){10}}
\put(80,50){\circle*{2}}
\put(90,40){\line(1,1){10}}
\put(90,40){\circle*{2}}
\put(100,50){\line(1,0){10}}
\put(110,50){\circle*{2}}
\put(100,50){\circle*{2}}
\put(120,80){Shape $\mathcal{C}_1$ is the deformation of shape $\mathcal{C}$ such that projection}
\put(120,60){of the first arrow on the $y$-axis becomes $-2$ instead of $-1$.}
\put(120,40){We denote such deformation, more precisely, as $\mathcal{C}_1^2$. Its points}
\put(10,20){have the following respective heights: $p,p,p+1,p-1,p-1,p,p,p-1,p,p$.}
\end{picture}
\end{center}

\noindent
In analogous way one defines deformations $\mathcal{C}_1^n$ for all positive integers $n$.
In this notation $\mathcal{C}=\mathcal{C}_1^1$. in analogous way we define deformation shape $\mathcal{C}_2$.
\begin{center}
\begin{picture}(420,90)
\put(20,60){\line(1,0){10}}
\put(20,60){\circle*{2}}
\put(30,60){\line(1,1){10}}
\put(30,60){\circle*{2}}
\put(40,70){\vector(1,-1){6}}
\put(40,70){\line(1,-1){10}}
\put(40,70){\circle*{2}}
\put(50,60){\line(1,0){10}}
\put(50,60){\circle*{2}}
\put(60,60){\line(1,1){10}}
\put(60,60){\circle*{2}}
\put(70,70){\line(1,0){10}}
\put(70,70){\circle*{2}}
\put(80,70){\vector(1,-3){6}}
\put(80,70){\line(1,-3){10}}
\put(80,70){\circle*{2}}
\put(90,40){\line(1,1){10}}
\put(90,40){\circle*{2}}
\put(100,50){\line(1,0){10}}
\put(110,50){\circle*{2}}
\put(100,50){\circle*{2}}
\put(120,80){Shape $\mathcal{C}_2$ is the deformation of shape $\mathcal{C}$ such that projection}
\put(120,60){of the second arrow on the $y$-axis becomes $-3$ instead of $-1$.}
\put(120,40){We denote such deformation, more precisely, as $\mathcal{C}_2^3$. Its points}
\put(10,20){have the following respective heights: $p,p,p+1,p,p,p+1,p+1,p-2,p-1,p-1.$}
\end{picture}
\end{center}

\noindent
Again we can define $\mathcal{C}_2^n$ for all positive integers $n$, in particular, $\mathcal{C}_2^1=\mathcal{C}$.

Now we consider the following symbolic sequence
\begin{gather}\label{eq:C-symbolic}
\mathcal{C} \mathcal{C}_1^2 \mathcal{C}_2^2 \mathcal{C} \mathcal{C}_1^{3^+}\mathcal{C}_2^2
\mathcal{C} \mathcal{C}_1^2 \mathcal{C}_2^{3^+}
\mathcal{C} \mathcal{C}_1^2 \mathcal{C}_2^2 \mathcal{C} \mathcal{C}_1^{3^+}\mathcal{C}_2^2
\mathcal{C} \mathcal{C}_1^2 \mathcal{C}_2^{3^+}
\ldots.
\end{gather}
Since the order of the shapes is preserved we can simplify notation because the sequence of the upper subscripts\footnote{For $\mathcal{C}$ we put $1$ because by definition $\mathcal{C}=\mathcal{C}^1=\mathcal{C}_1^1=\mathcal{C}_2^1$.}
immediately restore the whole symbolic sequence~\eqref{eq:C-symbolic}:
\begin{gather}\label{eq:C-digital}
\underbrace{12213^+2123^+}\underbrace{12213^+2123^+}12213^+2123^+\ldots.
\end{gather}
We see that both sequences~\eqref{eq:C-symbolic} and \eqref{eq:C-digital} are quasiperiodic with the quasiperiod underbraced
in sequence~\eqref{eq:C-digital}. Again the definition of $3^+$ is exactly the same as in sequence~\eqref{eq:ABsequence}.
Which means that $3^+>3$ at the resonances. Because the resonances occurs not in every period and the ``depths'' of these
resonances are different we call the sequence quasiperiodic.
To get the fence $\mathcal{P}$ from sequences~\eqref{eq:C-symbolic} and~\eqref{eq:C-digital} is simple: we put the left
end of the first shape $\mathcal{C}$ at the point with the coordinates $(1,0)$ and successively glue together the right
end of the previous shape with the left end of the following one.

Let us calculate $3^+$ in resonant points. Consider, first, this calculation with the help of sequence~\eqref{eq:C-symbolic}.
We begin with $3^+$ resonances in shapes $\mathcal{C}_1$. Assume that the resonance happens when $\mathcal{C}_1$ appears $N$th
time in sequence~\eqref{eq:C-symbolic}. Then the resonance happens at the point with $x$-coordinate $n=10+9\cdot3(N-1)+3$.
This point should coincide with one of the members of the sequence $a_k$. We see that at this resonance $a_{k+1}=a_{k}+1$.
This may happen only in case the triangular decomposition of $k$ reads as $(q,0)$. Thus, using equation~\eqref{eq:a-k} we arrive
at the following condition for $N$:
\begin{gather}\label{eq:C1-resonances}
13+3^3(N-1)=\frac{3^q+1}{2}-1 \quad \Longrightarrow \quad
N=\frac{3^{q-3}-1}{2}+1 \quad \Longrightarrow \quad
n=b_{q-1}+1,
\end{gather}
where $b_{q-1}$ is defined in equation~\eqref{eq:b-k}. Conjecture~\ref{con:b-k} implies that $q-1>3$ if we want to get $3^+>3$.
So, we get $\mathcal{C}_1$-resonances with $3^+>3$ iff $q=5,6,\ldots$, namely,
\begin{gather*}
n =121, 364, 1093, 3280, 9841, 29524,\ldots,\\ %\label{eq:1resonances}\\
z_{n-1} = 4,\quad 5,\quad 6,\quad 7,\quad 8,\quad 9,\quad n \quad \ldots.\nonumber
\end{gather*}
We remind the reader that according to our definition $z_n=0$ and $z_{n+1}=0$ for $n$ given by equation~\eqref{eq:C1-resonances}.

Now consider $3^+$ resonances in shapes $\mathcal{C}_2$. Let us use the symbolic sequence~\eqref{eq:ABsequence}. Obviously,
the case $3^+$-deformations $\mathcal{C}_2$, occurs after each third appearance of shape~$\mathcal{B}$ in \eqref{eq:ABsequence}.
Therefore, we are interested in $3N$th appearance of $\mathcal{B}$ in sequence~\eqref{eq:ABsequence}: it happens when $n=-1+81N$.
To study resonances we have to consider equation
\begin{gather}\label{eq:Nql}
-1+81N=\frac{3^q+3^l}{2}-1 \quad \Longrightarrow\quad
N=\frac{3^{q-4}+3^{l-4}}{2}\quad \Longrightarrow\quad
4\leq l\leq q,
\end{gather}
the last condition comes from the fact that $N$ is a positive integer. Finally, substituting $N$ into equation, $n=-1+81N$ we
arrive at the conclusion that $3^+$ resonances of type $\mathcal{C}_2$ occurs at
\begin{gather*}
n=a_k,\qquad{\rm where}\qquad
k=(q,l),\qquad
4\leq l\leq q.
\end{gather*}
It is more complicated to distinguish cases when $3^+>3$. We have only the points $b_k^{k'}$ (see equation~\eqref{eq:p-nforn<b-k+1})
for which Conjecture~\ref{con:b-k-k+1} says that height of the fence equals $k$. Therefore it is natural to look
whether $\mathcal{C}_2$-resonance may happen after some of these points?
We have the following condition for numbers~\eqref{eq:p-nforn<b-k+1}
\begin{gather}\label{eq:Nklm}
b_k^{(k-m)(k+m+5)/2+l}+1=-1+81N \quad \Longrightarrow\quad l=m+1, \qquad
N=\frac{3^{k-2}+3^{m-3}}2,
\end{gather}
where $l$ is a dummy variable, which has nothing to do with $l$ in equation~\eqref{eq:Nql}. Now comparing formulae for $N$
obtained in equations~\eqref{eq:Nql} and \eqref{eq:Nklm} one proves that $k=m$, $q=l$, and $q=k+1$. Note that condition
$3^+>3$ in this notation reads as $k\geq4$. Thus actually for every $k$ we found one $\mathcal{C}_2$-resonance
\begin{gather}\label{eq:C2-resonances}
n=b_{q-1}^q+1=3^q-1, \qquad z_{n-2}=z_{n-1}=q-1,\qquad q=5,6,\ldots.
\end{gather}
The first $\mathcal{C}_2$-resonances with $3^+>3$ defined by equation~\eqref{eq:C2-resonances} are as follows
\begin{gather*}
n =242, 728, 2186, 6560, 19682, 59048,\ldots,\\ %\label{eq:2resonances}\\
z_{n-1} = 5,\quad 6,\quad 7,\quad 8,\quad 9,\quad 10,\quad \ldots.\nonumber
\end{gather*}
Contrary to $\mathcal{C}_1$-resonances which are completely defined by the last equation in~\eqref{eq:C1-resonances}, it is not
clear whether equation~\eqref{eq:C2-resonances} describes all $\mathcal{C}_2$-resonances with $3^+>3$.

There is one more interesting property of quasiperiodic fence~$\mathcal{P}$. By definition after each
$\mathcal{C}_1$-resonance fence~$\mathcal{P}$ suffer a~gap\footnote{The height of the fence equals $0$ at the resonance and next point.} of the length $1$. Thus the fence
consists of infinite number of the connected parts. The first point of the $k$th-part is $b_{k-1}+2$ and the last
$b_{k}+1$, where we put formally $b_0=0$. By definition each connected part begins at a point with the zero height
and finishes at some other point with the zero height. Inside of the connected parts there are other points with
zero heights but they do not destroy the connectedness of these parts, they are just ``some faults'' in the
construction. Denote the area of the $k$th connected part of $\mathcal{P}$ as $S_k$. It is easy to observe that
$S_k$, $k=1,2,\ldots$, is the integer sequence. One finds its first terms:
\begin{gather*}
1, \ 7, \ 34, \ 142, \ 547, \ \ldots.
\end{gather*}
There is only one sequence A014915 in OEIS~\cite{OEIS} with these first terms. Therefore, it is natural to assume
\begin{Conjecture}
\begin{gather*}
S_k=\frac{(2k-1)3^k + 1}4.
\end{gather*}
\end{Conjecture}
In OEIS there is a recurrence relation, $S_{k+1}=(k+1)3^k+S(k)$, $S(1)=1$. This relation allow one to make
a conjecture that the first $b_k-(b_{k-1}+2)$ elements of the $(k+1)$-th connected part of the fence $\mathcal{P}$
exactly coincides with its $k$-th connected part without the very last unit segment. That means that the heights of
the fence $\mathcal{P}$ on the segments $[b_{k-1}+2, b_k]$ in the $k$-th part coincide with the corresponding
heights on the segment $[b_k+2,2b_k-b_{k-1}]$ of the $(k+1)$-th part. The ``corresponding'' heights mean the heights
measured at the points of the segments equidistant from the left ends of the segments. Since the lengths of the
segments coincide these points will be equidistant also from the right ends of the segments. The heights on the
left ends of the segments vanishing by construction, the heights on their right ends equal by these conjecture:
\begin{gather*}
z_{(b_{k-1}+2)}=z_{(b_k+2)}=0,\qquad z_{b_k}=z_{(2b_k-b_{k-1})}=k.
\end{gather*}
We can formulate our last conjecture in a bit different form: the $(k+1)$-th connected part consists of the ``old'' fragment, i.e.,
the fence built on the segment $[b_k+2,2b_k-b_{k-1}]$, and the ``new'' fragment, it is the fence built on the segment
$[2b_k-b_{k-1}+1,b_{k+1}+1]$. The old part coincides with the previous $k$-th connected part of the fence without the very
last segment. The total length of the $(k+1)$-th connected part is
\begin{gather*}
b_{k+1}+1-(b_k+2)+1=3(b_{k}+1)-b_{k}=2b_{k}+3.
\end{gather*}
The length of the old fragment is
\begin{gather*}
2b_k-b_{k-1}-(b_k+2)+1=b_{k}-b_{k-1}-1.
\end{gather*}
The length of the new fragment is
\begin{gather*}
b_{k+1}+1-(2b_k-b_{k-1}+1)+1=3(b_{k}+1)+1-(2b_k-b_{k-1})=b_{k}+b_{k-1}+4.
\end{gather*}
So we see that the new fragment of every connected part of the fence is longer than the old fragment, moreover,
its area is asymptotically two times larger than the area of the old fragment.

Now we can turn back to $\mathcal{C}_2$-resonances satisfying the condition $3^+>3$ and make a reasonable conjecture about
their location. It is clear that $\mathcal{C}_2$-resonances with $3^+>3$ which are defined in equation~\eqref{eq:C2-resonances}
belong to the new fragments of the connected parts of $\mathcal{P}$. On the other hand it is clear that
if we have $\mathcal{C}_2$-resonances with $3^+>3$ in $k$th connected part of~$\mathcal{P}$, then they reappear in the
old fragment of $(k+1)$-th connected part. So that the number of such resonances linearly grow with $k$. More precisely,
on $k$-th connected part located exactly $k+1$ points of the sequence~$a_{k'}$, $k'=(k,l)$, $l=0,\ldots,k$, including the end points.
It means that the following $(k+1)$-th part contains the images of the resonances from the $k$-th part and one more resonance
in the new part given by equation~\eqref{eq:C2-resonances}. At this stage it would be convenient to count all $\mathcal{C}_2$-resonances
not necessary those with $3^+>3$. Then we have the following
\begin{Conjecture}\label{con:C2-resonances}
All $\mathcal{C}_2$-resonances of the $k+1$-th connected part $(k=1,2,\ldots)$ of fence $\mathcal{P}$ with depth $l-1$ are
given by the sequence $a_{k'}$ where $k'=(k+1,l)$, $l=2,\ldots, k+1$.
\end{Conjecture}
In analogous way we can formulate our study of $\mathcal{C}_1$-resonances.
\begin{Conjecture}\label{con:C1-resonances}
All $\mathcal{C}_1$-resonances with the depth $k=1,2,\ldots$ are given by the sequence
$a_{k'}=b_k+1$ where $k'=(k+1,0)$. They coincide with the right end points of the $k$-th connected part of $\mathcal{P}$.
\end{Conjecture}
\begin{Conjecture}\label{con:nonresonance3+}
All nonresonant $3^+$-numbers of symbolic sequences~\eqref{eq:ABsequence}, \eqref{eq:C-symbolic}, and \eqref{eq:C-digital}
equal $3$.
\end{Conjecture}
Conjectures~\ref{con:C2-resonances}--\ref{con:nonresonance3+} completely define symbolic
sequences~\eqref{eq:ABsequence}--\eqref{eq:C-digital} and thus our quasipe\-rio\-dic fence $\mathcal{P}$.

\section{Monodromy data}\label{sec:monodromy}
This section is based on paper~\cite{KV2004}.
Here we explain how to use the results of \cite{KV2004} to get information about
asymptotics as $\tau\to\infty$ of a solution of equation~\eqref{eq:dp3} defined by its expansion as $\tau\to0$.
Sure we consider only the solution which is the main hero of this paper. In~\cite{KV2004} we studied the
general solution, while our case is a very degenerate one, therefore some additional efforts are required to
specify our solution.

The facts and notation we need from the paper~\cite{KV2004} would take a few pages, therefore
here we recall only some basic definitions, which allow the reader to follow the schemes of proofs and
understand the main statements. For the complete understanding of this section the reader should address
the corresponding places in paper~\cite{KV2004} we reference below.

We recall that according to \cite{KV2004} the pair of functions $\{u(\tau),\varphi(\tau)\}$, where $u(\tau)$ is
any solution of equation~\eqref{eq:dp3} and $\varphi(\tau)$ is defined as the general solution of the following ODE
\begin{gather}\label{eq:varphi}
\varphi'(\tau)=\frac{2a}{\tau}+\frac{b}{u(\tau)},
\end{gather}
can be uniquely parameterized with the points of {\it the manifold of monodromy data}, or just
{\it the monodromy manifold}, which, for a given parameter $a$, is an algebraic
variety of the complex dimension $3$.

Obviously, for a given $u(\tau)$, the function $\varphi(\tau)$ is defined
by equation~\eqref{eq:varphi} up to an additive parameter, $\varphi_0\in\mathbb C$,
$\varphi(\tau)\to\varphi(\tau)+\varphi_0$. In principle, it is not complicated to exclude the function $\varphi(\tau)$,
and contract (consider bilinear combinations of some coordinates) manifold of the monodromy data to the complex
dimension $2$, so that it would parameterize solely solutions of equation~\eqref{eq:dp3}, however it is not done in
\cite{KV2004} and we follow that definitions not to confuse the reader.
The parameter $a$, which is the coefficient of equation~\eqref{eq:dp3} is called, sometimes, the formal monodromy
and in the form ${\rm e}^{\pi a}$ it enters the algebraic equations defining the monodromy manifold.

Consider ${\mathbb C}^8$ with the coordinates (monodromy data) denoted as $a$, $s_0^0$, $s_0^\infty$, $s_1^\infty$,
$g_{11}$, $g_{12}$, $g_{21}$, and $g_{22}$. The monodromy manifold is defined by the following system of algebraic
(if we turn from $a$ to ${\rm e}^{\pi a}$) equations (system~(33) in \cite{KV2004})
\begin{gather}
 s_0^{\infty}s_1^{\infty}=-1-{\rm e}^{-2\pi a}-is_0^0{\rm e}^{\pi a}, \qquad
 g_{22}g_{21}-g_{11}g_{12}+s_0^0g_{11}g_{22}=i{\rm e}^{-\pi a},\nonumber\\
 g_{11}^2-g_{21}^2-s_0^0g_{11}g_{21}=i{\rm e}^{-\pi a}s_0^\infty,\qquad
 g_{22}^2-g_{12}^2+s_0^0g_{22}g_{12}=i{\rm e}^{\pi a}s_1^\infty,\nonumber\\
 g_{11}g_{22}-g_{12}g_{21}=1.\label{monodromy-manifold}
\end{gather}
The main goal of this section is to find for our solution, $u(\tau)$, the monodromy data. In the next section we
use them to get asymptotics of $u(\tau)$ as $\tau\to\infty$.

Since we know the behavior of the solution at $\tau=0$, we have to check whether Theorems~3.4 and 3.5 of \cite{KV2004} describing asymptotics as $\tau\to0$
of solutions of equation~\eqref{eq:dp3} are applicable to it. Below, until system~\eqref{sys:monodromi-pm-i/4}
we discuss how one can get the monodromy data for $u(\tau)$ with the help of these theorems.

Asymptotics as $\tau\to0$ of the general solution of equation~\eqref{eq:dp3} is given by equation~(45)
of~\cite{KV2004} (Theorem~3.4 of~\cite{KV2004}). In this equation is assumed that
\begin{gather}\label{eq:|a|<1}
|\operatorname{Im}(a)|<1.
\end{gather}
The equation contains a parameter $\rho$, which defines branching ($\tau^{\pm4\rho}$) of
the general solutions as $\tau\to0$.

Our solution is holomorphic at $\tau=0$, therefore, at first glance,
the last equation implies that the branching parameter $\rho$ should vanish.
Since equation~(45) of \cite{KV2004} is valid when $\rho\neq0$, we have to use Theorem~3.5 of~\cite{KV2004}.
In this case we have to ``kill'' the logarithmic terms in equation~(51) of~\cite{KV2004}. It is equivalent,
see system~(48) of \cite{KV2004} where $z_1=z_2=0$, to the condition $\det\{g_{ij}\}=0$, while according to the last
equation in system~\eqref{monodromy-manifold} (see above) this determinant equals $1$. So, the simplest natural
assumption is wrong.

The second natural assumption, which also leads to a singlevalued solution at $\tau=0$ is $\rho=\pm1/4$.
In this case equation~(44) of \cite{KV2004} implies
\begin{gather*}%\label{eq:s-for-a=i/2}
s_0^0=0.
\end{gather*}
In fact, all equations in Theorem~3.4 of~\cite{KV2004} are symmetric with respect to the reflection $\rho\to-\rho$,
we put $\rho=1/4$. Since we are interesting in the solution vanishing at $\tau=0$, we have to impose an
additional condition on the parameters of equation~(45) of~\cite{KV2004}
\begin{gather*}
\varpi_1^\natural(\varepsilon_1,\varepsilon_2;-1/4)\varpi_2^\natural(\varepsilon_1,\varepsilon_2;-1/4)=0.
\end{gather*}
Note that, in our case $\varepsilon_1=\varepsilon_2=0$.
Assume $\varpi_2^\natural(\varepsilon_1,\varepsilon_2;-1/4)=0$, then the first equation of system~(48)
of \cite{KV2004}, implies, $g_{12}=-g_{22}$. Now, if we put $s_0^0=0$, and $g_{12}\leftrightarrow-g_{22}$
into the second equation of the first row of system~\eqref{monodromy-manifold}, then we find,
$\det\{g_{ij}\}=i{\rm e}^{-\pi a}$. Comparing it with the last equation of system~\eqref{monodromy-manifold}
and condition~\eqref{eq:|a|<1} we get $a=i/2$. Lemma~\ref{lem:odd} says that in this case there might exists
one parameter family of the solutions holomorphic at $\tau=0$. The Suleimanov solution discussed in Introduction
belongs to this family. Our goal is to find the corresponding monodromy data uniquely characterizing this solution.
Using the relations on the monodromy data which are already obtained above, we find from the second equation
in the second row of system~\eqref{monodromy-manifold} that $s_1^{\infty}=0$.

We continue to analyze equation~(45) of \cite{KV2004}. Solution $u(\tau)$ has the leading term of
asympto\-tics~$ib\tau$, see equation~\eqref{eq:zero-expansion} for $a=i/2$. This result can be reproduced
via equation~(45) of \cite{KV2004}: if after the straightforward calculations with the help of
equations~(46)--(48) of \cite{KV2004} we demand $(g_{11}+g_{21})(g_{12}-g_{22})=-2$. This equation holds
in our case because: $g_{12}=-g_{22}$ and the last equation of system~\eqref{monodromy-manifold}.
Thus we have two complex parameters, say $g_{11}$ and $s_0^\infty$, for characterization of two functions
$u(\tau)$ and $\varphi(\tau)$.

The next term to the leading one in equation~(45) of \cite{KV2004} is defined by the following sum
\begin{gather}\label{eq:correction-term}
\frac{b\tau}{16\pi}{\rm e}^{\frac{i\pi}4}
\big(\varpi_1^\natural(\varepsilon_1,\varepsilon_2;1/4)\varpi_2^\natural(\varepsilon_1,\varepsilon_2;1/4)
\tau^{4\rho}+O\big(\tau^\delta\big)\big),\qquad\delta>0.
\end{gather}
Since $4\rho=1$, then, in case we would know that $\delta>1$, we can equate the coefficient of the leading term
$\tau^{4\rho}$ in equation~\eqref{eq:correction-term} to the parameter $c_0$ of the Taylor
expansion~\eqref{eq:tau-Taylor}.
With the help of this equation we would be able to determine all the monodromy data uniquely characteri\-zing~$u(\tau)$ for $a=i/2$. In the odd case (the Suleimanov solution) $c_0=0$ and we arrive at the following equation
\begin{gather}\label{eq:varepsilon-1-1/4}
\varpi_1^\natural(\varepsilon_1,\varepsilon_2;1/4)=0.
\end{gather}
The second possibility, $\varpi_2^\natural(\varepsilon_1,\varepsilon_2;1/4)=0$ contradicts
system~\eqref{monodromy-manifold}.

In fact, one more monodromy parameter can be (correctly!) fixed with the help of equation~\eqref{eq:varepsilon-1-1/4},
although Theorem 3.4 of \cite{KV2004} declares only inequality $\delta>0$, so that, strictly speaking, we are not
allowed to use equation~\eqref{eq:varepsilon-1-1/4}.
In view of the technique used in \cite{KV2004}, it is, most probably, possible either to get a more accurate error
estimate for the general solutions, or, at least, for our special one; however, it would require much more efforts
in the general case, or separate consideration of our solution.
It is the manifestation of the degeneracy of the solution $u(\tau)$ mentioned in the beginning of this section:
the leading term of asymptotics at $\tau=0$ does not allow to determine the complete set of the monodromy parameters
uniquely characterizing the solution.

Let us assume that equation~\eqref{eq:varepsilon-1-1/4} is valid and obtain the corresponding set of the monodromy
data. Equation~\eqref{eq:varepsilon-1-1/4} implies, $g_{21}=g_{11}$. Then, the first equation in the second row of
system~\eqref{monodromy-manifold} gives $s_0^{\infty}=0$. We omit analogous considerations for $a=-i/2$ and
formulate the final result for two solutions:
\begin{gather}
a=\frac{i}2, \qquad g_{12} =-g_{22},\qquad g_{21}=g_{11},\qquad g_{11}g_{22}=-g_{12}g_{21}=\frac12,\qquad s_0^0=s_0^\infty=s_1^\infty=0,\nonumber\\
 a=-\frac{i}2, \qquad g_{12} =g_{22},\qquad g_{21}=-g_{11},\qquad g_{11}g_{22}=-g_{12}g_{21}=\frac12,\nonumber\\
 s_0^0=s_0^\infty=s_1^\infty=0. \label{sys:monodromi-pm-i/4}
\end{gather}
At this stage equations~\eqref{sys:monodromi-pm-i/4} are not rigorously confirmed. Below we give another rigorous
derivation valid for all values of $a$ satisfying condition~\eqref{eq:|a|<1}.

Now we turn to the case of general $a$ (restriction~\eqref{eq:|a|<1} is revoked). Since the value of the parameter
$\rho$ in this case is not obvious, we begin with the fact that our solution $u(\tau)$ is holomorphic in a
neighbourhood of $\tau=0$. With the help of equation~\eqref{eq:varphi} one confirms that the same is true for the
function $\varphi(\tau)$. This means that after analytic continuation around $\tau=0$, system~(12) of~\cite{KV2004}
does not change.

Unfortunately, definition~(15) of \cite{KV2004} of the canonical solutions of this system contains inaccuracy,
namely, the correct definition should read
\begin{gather*}%\label{eq:canonical-solution-infinity}
Y_k^\infty(\mu)\underset{\substack{\mu\to\infty\\\mu\in\Omega_k^\infty}}=\left( I+\frac{\Psi^{(1)}}{\mu}+
\frac{\Psi^{(2)}}{\mu^2}+\cdots \right)
\exp\left( -i\left(\tau\mu^2+\left(a-\frac{i}{2}\right)\ln\mu-\frac{a}{2}\ln\tau\right)\sigma_3 \right),
\end{gather*}
where $\sigma_3=\left(\begin{smallmatrix}1&0\\0&-1\end{smallmatrix}\right)$ and the sectors $\Omega_k^\infty$ are defined on p.1170
of \cite{KV2004}.
The last term, $-\frac{a}{2}\ln\tau$ in the exponent above is absent in~\cite{KV2004} and~\cite{KV2010}. This
incorrectness, does not have any effect on the definitions and results presented in Sections~2 and 3 of these
papers.\footnote{The incorrectness appeared because of the change of notation, one can simplify definition of
the canonical solution given above and use the original one given in~\cite{KV2004}, but with a simultaneous
gauge transformation of system~(12) of~\cite{KV2004}. Finally, this difference in the definition of the canonical
solution resulted only in a possible appearance of the term $a\ln\tau$ in asymptotics of the function
$\varphi(\tau)$, which is not included in the list of the main results. In the next publication on special solutions
of equation~\eqref{eq:dp3} we are going to check what definition of the canonical solution $Y_k^\infty(\mu)$ was
used for actual calculation of asymptotics.}
If we take the canonical solution $Y_k^\infty(\mu)=Y_k^\infty(\mu,\tau)$ of this
system, then after the analytic continuation with respect to $\tau$, $\tau\to\tau\cdot {\rm e}^{-2\pi i}$,
the variable $\mu$ belongs to the same sector where the canonical solution
$Y_{k+2}^\infty(\mu)=Y_{k+2}^\infty(\mu,\tau)$ is defined;
see definition of the sectors $\Omega_k^\infty$ at p.~1170 of~\cite{KV2004}, where the argument of $\tau$ enters
the definition of the sectors. Since at every value of $\tau$ during this continuation $Y_k^\infty(\mu)$ keeps the
same (canonical) asymptotics, therefore after arriving at the sector, where the canonical solution
$Y_{k+2}^\infty(\mu)$ is defined, it has exactly the same asymptotics as~$Y_{k+2}^\infty(\mu)$,
but with $\tau\to\tau\cdot {\rm e}^{-2\pi i}$, and solves system~(12) of~\cite{KV2004} with exactly the same coefficients.
Therefore, $Y_{k+2}^\infty(\mu,\tau)=Y_k^\infty\big(\mu,\tau\cdot {\rm e}^{-2\pi i}\big){\rm e}^{-\pi a\sigma_3}$ for all $k\in\mathbb Z$. This relation between
the canonical solutions immediately implies the following relation for the Stokes matrices
\begin{gather*}
S_{k+2}^\infty={\rm e}^{\pi a\sigma_3}S_k^\infty {\rm e}^{-\pi a\sigma_3}.
\end{gather*}
Comparing the above equation with equation~(23) of \cite{KV2004}
\begin{gather}\label{eq:23ofKV2004}
S_{k+2}^\infty=\sigma_3{\rm e}^{-\pi(a-i/2)\sigma_3}S_k^\infty {\rm e}^{\pi(a-i/2)\sigma_3}\sigma_3=
{\rm e}^{-\pi a\sigma_3}S_k^\infty {\rm e}^{\pi a\sigma_3},
\end{gather}
we get that for every integer $k$
\begin{gather*}
S_k^\infty={\rm e}^{2\pi a\sigma_3}S_k^\infty {\rm e}^{-2\pi a\sigma_3}.
\end{gather*}
Each Stokes matrix is known to have the triangular structure with units on the diagonal and one (generally
nontrivial) off-diagonal element, $s_k^\infty$, called the Stokes multiplier. The last equation implies for the
Stokes multipliers the following equation
\begin{gather*}
s_k^\infty={\rm e}^{\pm4\pi a}s_k^\infty.
\end{gather*}
Therefore, we arrive at the conclusion that for $a\neq in$ and $a\neq i/2+in$, $n\in\mathbb Z$, the
meromorphic solution of equation~\eqref{eq:dp3} vanishing at the origin has vanishing Stokes multipliers
\begin{gather}\label{eq:s0=s1=0}
s_0^\infty=s_1^\infty=0.
\end{gather}

We have to cope with the two remaining cases of $a$: $a=in$ and $a=i/2+in$. Theorem~\ref{th:existence & uniqueness}
says, that for $a=i/2+in$ there exists the unique odd meromorphic solution vanishing at $\tau=0$. Sure, by the
continuity argument (monodromy data depends analytically on $a$) we can prove that equation~\eqref{eq:s0=s1=0}
holds also for this case.
As for the case $a=in$, we know that for $a=0$, $\pm i$,\ldots,$\pm5i$ (see formulae for $u_{2n}$ underneath
Remark~\ref{rem:coeff3}) the meromorphic solution vanishing at $\tau=0$ does not exists. Sure explicit calculations
can be continued further and the nonexistence can be confirmed for the larger values of $n$.
Since Section~\ref{sec:residues} says that the nontrivial function generating residues at $a=in$ can be constructed
for any $n$, such solutions do not exists for all $n\in\mathbb Z$. The proof of Lemma~\ref{lem:odd} shows that if
holomorphic solution at $\tau=0$ exists for some $a=in_0$, $n_0\in\mathbb Z$, it is odd and not unique.
The continuity argument, analogous to the one given above for the case $a=i/2+in$, shows that
equation~\eqref{eq:s0=s1=0} should hold for at least one limiting case as $a\to in_0$. However, it contradicts
to the last equation of system~\eqref{sys:monodromy-data-u-varphi-unique} defining the monodromy manifold.

Before going further, let us consider the symmetry for system~(12) of \cite{KV2004}
related with the odd solutions $u(\tau)$ and reproduce condition~\eqref{eq:s0=s1=0} for all such solutions.
This symmetry is considered in \cite[Section~6.2, item~6.2.1, p.~1199]{KV2004} and requires correction.

So, we assume that for some solution $u(\tau {\rm e}^{i\pi})=-u(\tau)$ in a neighborhood of $\tau=0$. Then obviously,
$u(\tau)=-u(\tau {\rm e}^{-i\pi})$. The coefficients of system~(12) of~\cite{KV2004} are denoted as $A(\tau)$, $B(\tau)$,
$C(\tau)$, and $D(\tau)$. Equations in Proposition~1.2 of~\cite{KV2004} shows that
\begin{gather*}
A(\tau)=A(-\tau),\qquad
B(\tau)=B(-\tau),\qquad
C(\tau)=-C(-\tau),\qquad
D(\tau)=-D(-\tau).
\end{gather*}
These equations implies the following relations for the canonical solutions:
\begin{gather}
Y_k^\infty\big({\rm e}^{i\pi}\mu,{\rm e}^{-i\pi}\tau\big)={\rm e}^{-\frac{i\pi}4\sigma_3}Y_k^\infty(\mu,\tau){\rm e}^{\pi a\sigma_3},\label{eq:Yk=k}\\
Y_{k+2}^\infty\big({\rm e}^{i\pi}\mu,{\rm e}^{i\pi}\tau\big)={\rm e}^{-\frac{i\pi}4\sigma_3}Y_k^\infty(\mu,\tau).\label{eq:Yk=k+2}
\end{gather}
Equations~\eqref{eq:Yk=k} and \eqref{eq:Yk=k+2} imply for the following relations Stokes matrices
\begin{align*}
S_k^\infty={\rm e}^{-\pi a\sigma_3}S_k^\infty {\rm e}^{\pi a\sigma_3},\qquad
S_{k+2}^\infty=S_k^\infty,\qquad
k=0,\pm1,\ldots,
\end{align*}
respectively. With the help of equation~\eqref{eq:23ofKV2004} we see that both above equations are equivalent
to the following equation for the Stokes multipliers
\begin{gather}\label{eq:stokes-odd}
s_k^\infty=s_k^\infty {\rm e}^{2\pi a}.
\end{gather}
Equation~\eqref{eq:stokes-odd} implies condition~\eqref{eq:s0=s1=0} for all odd solutions of equation~\eqref{eq:dp3}
in the neighborhood of $\tau=0$ and $a\neq in$, $n\in\mathbb Z$.

Analogous reasoning does not work for the canonical solutions $X_k(\mu)$, in the neighbourhood of $\mu=0$,
because their asymptotics contains terms $\tau^{1/2}$ and $\tau^{1/4}$. Actually, the first equation of
system~\eqref{monodromy-manifold} implies
\begin{gather*}
s_0^0=2i\cosh(\pi a)=i{\rm e}^{\pi a}+i{\rm e}^{-\pi a}.
\end{gather*}
Moreover, equation~(44) of \cite{KV2004}, together with the restrictions (43) of \cite{KV2004}, implies
$\rho=\pm\frac{ai}2$. Again we can choose here any sign of $\rho$ because asymptotic formula~(45) of \cite{KV2004}
is symmetric with respect~$\rho$. Calculations slightly simpler with $\rho=-\frac{ai}2$. After that a simple
analysis of system~\eqref{monodromy-manifold} allows us to prove the following
\begin{Proposition}\label{prop:monodromy-final}
For $a\notin i\mathbb Z/2$ there exists the only one solution of equation~\eqref{eq:dp3} such that both functions
$u(\tau)$ and ${\rm e}^{i\varphi(\tau)}$ are meromorphic. It is an odd function of $\tau$. For $a=i/2+in$ with
$n\in\mathbb Z$ there exists a unique odd solution of equation~\eqref{eq:dp3} such that both functions~$u(\tau)$ and~${\rm e}^{i\varphi(\tau)}$ are meromorphic.
Their monodromy parameters are as follows
\begin{gather}
s_0^\infty=s_1^\infty=0,\qquad s_0^0=2i\cosh(\pi a),\nonumber\\
g_{21}=-i{\rm e}^{\pi a}g_{11},\qquad
g_{12}=i{\rm e}^{\pi a}g_{22},\qquad
g_{11}g_{22}\big(1-{\rm e}^{2\pi a}\big)=1.\label{sys:monodromy-data-u-varphi-unique}
\end{gather}
\end{Proposition}
\begin{proof}The main part of the proof is given before Proposition. The last equation of system~\eqref{sys:monodromy-data-u-varphi-unique} implies $a\notin i\mathbb Z$, which is consistent with
equation~\eqref{eq:u2ka}.

To finish the proof we have to notice that the monodromy data given by system~\eqref{sys:monodromy-data-u-varphi-unique} contain one parameter, say, $g_{12}$ or $g_{21}$.
As soon as this parameter is fixed the others are uniquely defined. This parameter defines the constant of
integration in equation~\eqref{eq:varphi} and does not effect on the function~$u(\tau)$. The function $u(\tau)$
is uniquely determined by the bilinear combinations of the monodromy data: $g_{11}g_{22}$ and $g_{11}g_{12}$.
In our case these combinations are uniquely defined as long as the parameter~$a$ is fixed.
\end{proof}

\begin{Remark}Note, that system~\eqref{sys:monodromy-data-u-varphi-unique} for $a=\pm i/2$ coincides with
the corresponding system~\eqref{sys:monodromi-pm-i/4}, therefore the latter systems are proved.
\end{Remark}
\begin{Remark}The functions $u(\tau)$ and $\varphi'(\tau)$ are odd and meromorphic. Note that function $\varphi(\tau)$ is not a~meromorphic function.
\end{Remark}
\begin{Remark}The solution $u(\tau)$ is not the only meromorphic solution of equation~\eqref{eq:dp3}. For the other
meromorphic solutions the corresponding functions ${\rm e}^{i\varphi(\tau)}$ are not single-valued.
\end{Remark}

\section[Asymptotics as $\tau\to+\infty$]{Asymptotics as $\boldsymbol{\tau\to+\infty}$}\label{sec:asymptotics}

Here we apply the results obtained in the previous section and \cite{KV2004} to get asymptotics of $u(\tau)$ for
the large values of~$\tau$. In this section we assume the following restrictions on the coefficients of equation~\eqref{eq:dp3}
\begin{gather*}
|\operatorname{Im}a|<1,\qquad b>0.
\end{gather*}
First of all we have to check the conditions on the monodromy data for applicability of Theorem~3.1 of~\cite{KV2004}.
There are two such conditions
\begin{gather}
g_{11}g_{12}g_{21}g_{22}\neq0,\label{ineq:cond1theorem31}\\
\left|\operatorname{Re}\left(\frac{i}{2\pi}\ln(g_{11}g_{22})\right)\right|<1/6.\label{ineq:cond2theorem31}
\end{gather}
Inequality~\eqref{ineq:cond1theorem31} is an obvious consequence of the second line of equations
of system~\eqref{sys:monodromy-data-u-varphi-unique}. The second condition should be examined more carefully.
Substituting the last equation of system~\eqref{sys:monodromy-data-u-varphi-unique} into
condition~\eqref{ineq:cond2theorem31} we find
\begin{gather}\label{ineq:cond2a}
\left|\operatorname{Re}\left(\frac{i}{2\pi}\ln\big(1-{\rm e}^{2\pi a}\big)\right)\right|<1/6.
\end{gather}
\begin{Remark}\label{rem:weaker-condition}
The leading term of asymptotics $u(\tau)$ obtained in \cite{KV2004} contains the $\cosh$-function which by
definition can be written as the half-sum of two exponents. Restriction~\eqref{ineq:cond2theorem31} obtained
in \cite{KV2004} guarantee that both exponents are greater than the correction term. In Appendix $B$ of our
subsequent paper~\cite{KV2010} we have corrected the phase-shift in the $\cosh$-function obtained in~\cite{KV2004}
and also pointed out that in case we require that only the largest exponent of the $\cosh$-function is greater than
the correction term, then restriction~\eqref{ineq:cond2theorem31} is weaker, namely, $1/6$ in the r.h.s.\ should be
changed by $1/2$. For a better numerical correspondence of the leading term of asymptotics with the exact solution
outside restriction~\eqref{ineq:cond2a} one should use the result obtained in \cite{KV2010}.
From the point of view of the complete asymptotic expansions, which are not yet considered, it means
rearrangement of the corresponding series.
\end{Remark}

\begin{Remark}
The function $\ln$ is multivalued: the sense of equation~\eqref{ineq:cond2a} is that there should exists the
branch of $\ln$-function such that the restriction holds. If such branch of $\ln(\cdot)$ exists, then it is
fixed uniquely and it is the branch which should be used in the corresponding asymptotic formula. This remark
applies also to the weaker condition discussed in Remark~\ref{rem:weaker-condition}.
\end{Remark}

So modulo these remarks the asymptotics of $u(\tau)$ according to Theorem~3.1 of~\cite{KV2004} and the
correction made in Appendix~B of~\cite{KV2010} reads
\begin{gather}
u(\tau)\underset{\tau\to+\infty}{=}u_{\rm as}(\tau)+o\big(\tau^{-\delta}\big),
\qquad \delta>0,\label{eq:asymp-main-1}\\
u_{\rm as}(\tau)=\frac{b^{1/2}}{3^{1/4}}
\left(\sqrt{\frac{\vartheta(\tau)}{12}}
+\sqrt{\nu+1} {\rm e}^{\frac{3\pi i}{4}}\cosh(i\vartheta(\tau)+(\nu+1)\ln\vartheta(\tau)+z)\right),\label{eq:asymp-main-2}\\
\vartheta(\tau)=3^{3/2}b^{1/3}\tau^{2/3},\qquad
\nu+1=\frac{i}{2\pi}\ln(g_{11}g_{22}),\label{eq:asymp-main-3}\\
z=\frac{\ln(2\pi)}{2}-\frac{\pi i}{2}-\frac{3\pi i}{2}(\nu+1)+ia\ln(2+\sqrt{3})+(\nu+1)\ln12\nonumber\\
\hphantom{z=}{}-\ln\left(\omega\sqrt{\nu+1}\Gamma(\nu+1)\right),\qquad \omega=g_{11}g_{12}.\label{eq:asymp-main-4}
\end{gather}
The error estimate is written in equation~\eqref{eq:asymp-main-1} in a different form comparing with
Theorem~3.1 of~\cite{KV2004} because here it is convenient to introduce notation $u_{\rm as}(\tau)$ for
the leading term of asymptotics. Since $\delta>0$ does not fixed the present formulation is equivalent
to the original one in~\cite{KV2004}. In equations~\eqref{eq:asymp-main-1}--\eqref{eq:asymp-main-4},
we put $\varepsilon_1=\varepsilon_2=0$ and $\varepsilon=1$. The very last equation in~\eqref{eq:asymp-main-4}
is corrected accordingly Appendix~B of~\cite{KV2010} (see equation~(B.6) at p.~53 of~\cite{KV2010}).

Consider condition \eqref{ineq:cond2a} and the weaker one mentioned in Remark~\ref{rem:weaker-condition} more
explicitly, note that we have to take into account also equation~\eqref{eq:|a|<1}. Here, we address two specific
cases of the parameter $a$, namely, $a$ is purely imaginary and $a$ is real.

\subsection[$a=i\alpha$, $\alpha\in\mathbb R$]{$\boldsymbol{a=i\alpha}$, $\boldsymbol{\alpha\in\mathbb R}$}
In this case condition~\eqref{ineq:cond2a} reads
$1/6<|\alpha|<5/6$, the corresponding condition from Remark~\ref{rem:weaker-condition} is $0<|\alpha|<1$.
As we mentioned above in both cases the leading term of asymptotics is given by
equations~\eqref{eq:asymp-main-1}--\eqref{eq:asymp-main-4}. However, the leading term numerically describes the solution better when~$\alpha$ is closer to $1/2$.
Below I present the real and imaginary part of the leading term of asymptotics of the function $u(\tau)$.
To make it more explicit we consider the case $\alpha\in (0,1/2)\cup(1/2,1)$:
\begin{gather*}
\operatorname{Re}u_{\rm as} =\frac{b^{2/3}}2\tau^{1/3}+\frac{b^{1/2}}{3^{1/4}}\sqrt{\nu_1}
\left(\sin(\pi/4+\psi)\cosh(\chi)\cos(\phi)+\cos(\pi/4+\psi)\sinh(\chi)\sin(\phi)\right) ,\\
\operatorname{Im}u_{\rm as} =\frac{b^{1/2}}{3^{1/4}}\sqrt{\nu_1}
\left(\sin(\pi/4+\psi)\sinh(\chi)\sin(\phi)-\cos(\pi/4+\psi)\cosh(\chi)\cos(\phi)\right),
\end{gather*}
where
\begin{gather*}
\nu_1 =\frac12\sqrt{(\alpha-1/2)^2+(\ln(2\sin(\pi\alpha)))^2/\pi^2},\\
\psi =-\frac12\arctan\left(\frac{\ln(2\sin(\pi\alpha))}{\pi(\alpha-1/2)}\right) +\frac{\pi}4(\operatorname{sign}(\alpha-1/2)-1),\\
\chi =\frac12(\alpha-1/2)\ln\left(3\sqrt{3} b^{1/3}\tau^{2/3}\right)+\chi_0,\\
\chi_0 =\frac{\ln(2\pi)}{2}+\frac{\ln(2\sin(\pi\alpha))}{4}-\alpha\ln(2+\sqrt{3})+\frac{\ln(12)}{2}(\alpha-1/2)\\
\hphantom{\chi_0 =}{} -\ln\left(\sqrt{\nu_1}\left|\Gamma\left(\frac12(\alpha-1/2)
-\frac{i}{2\pi}\ln(2\sin(\pi\alpha))\right)\right|\right),\\
\phi =3\sqrt{3} b^{1/3}\tau^{2/3}-\frac{\ln(2\sin(\pi\alpha))}{2\pi}\ln\big(3\sqrt{3} b^{1/3}\tau^{2/3}\big)+\phi_0,\\
\phi_0 =\frac{\pi}{2}-\psi-\frac{3\pi}4(\alpha-1/2)-\frac{\ln(12)}{2\pi}\ln(2\sin(\pi\alpha)))\\
\hphantom{\phi_0 =}{} -\arg\left(\Gamma\left(\frac12(\alpha-1/2)-\frac{i}{2\pi}\ln(2\sin(\pi\alpha)\right)\right).
\end{gather*}
In the formulae above all roots and fractional powers of positive numbers are positive, the function
$\ln$ with positive argument is positive, $\arctan$ is the principle branch of the corresponding inverse
function with the values in $(-\pi/2,\pi/2)$, and $\operatorname{sign}(z)=1$ for $z>0$ and $-1$ for $z<0$, respectively.

\begin{figure}[h!]\centering
\includegraphics[width=100mm]{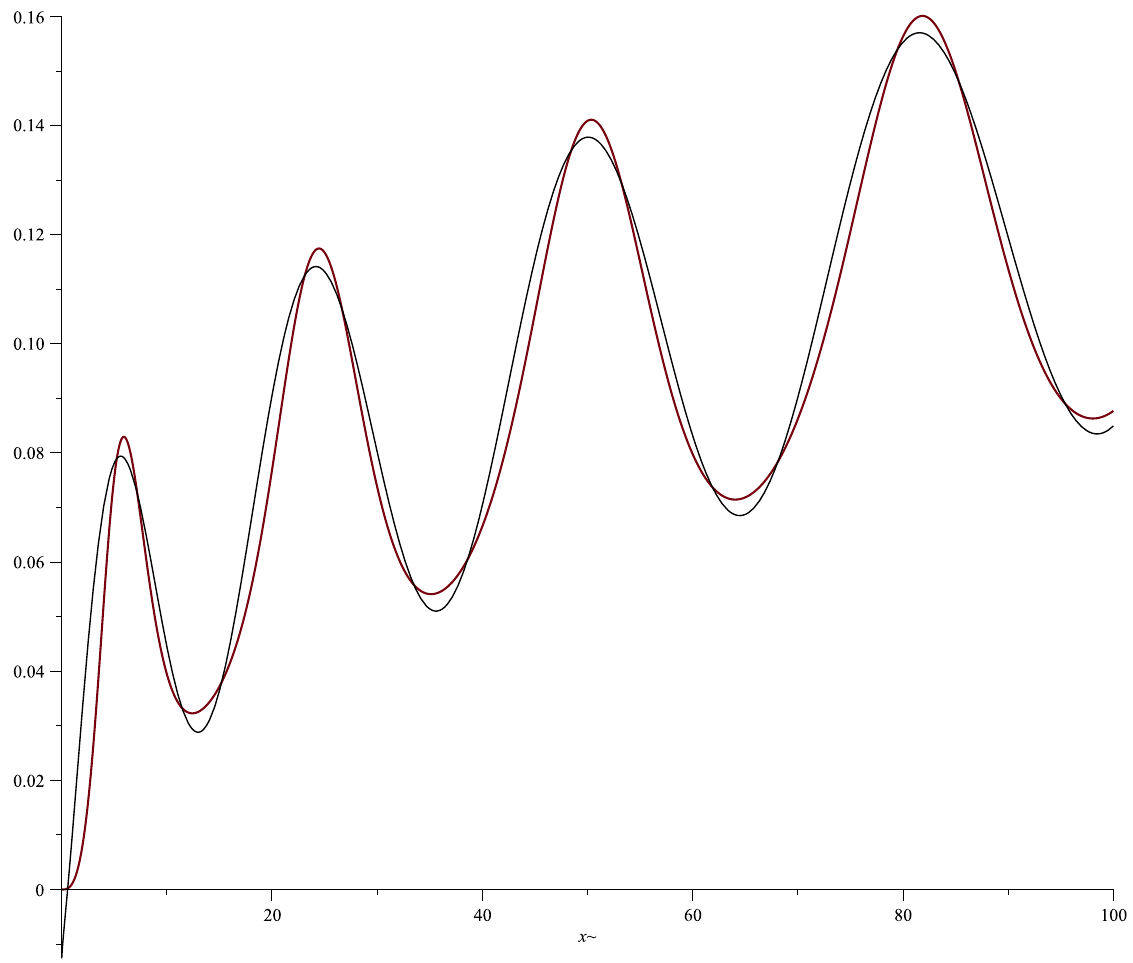}
\caption{The plots of $\operatorname{Re}u(\tau)$ and $\operatorname{Re}u_{\rm as}(\tau)$ for $\alpha=2/7$ and $b=1/80$.
The first one has higher corresponding extrema.}\label{fig:Re(u)alpha2d7b1d80}
\end{figure}

\begin{figure}[h!]\centering
\includegraphics[width=100mm]{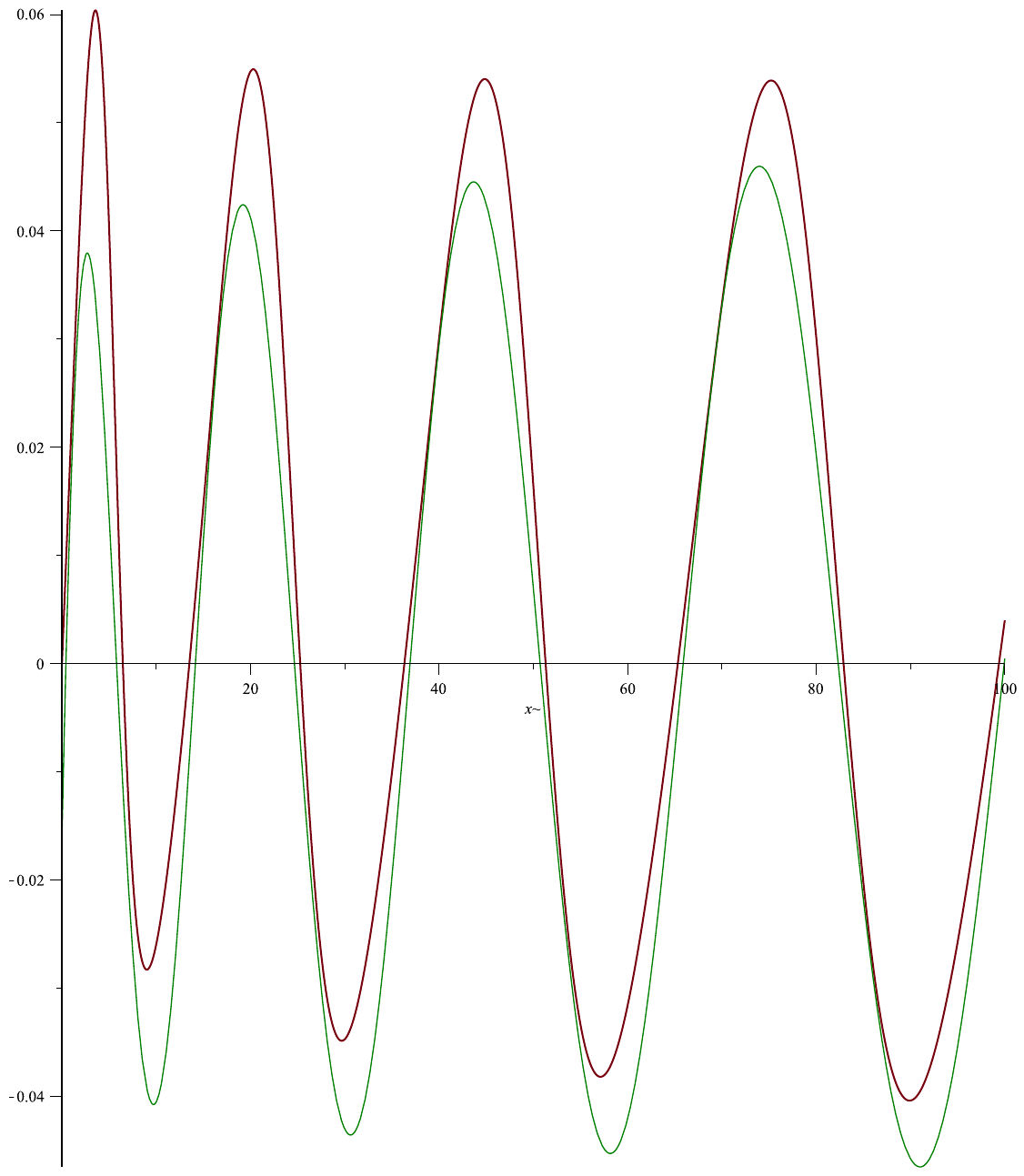}
\caption{The plots of $\operatorname{Im}u(\tau)$ and $\operatorname{Im}u_{\rm as}(\tau)$ for $\alpha=2/7$ and $b=1/80$.
The first one has higher corresponding extrema.}\label{fig:Im(u)alpha2d7b1d80}
\end{figure}

\begin{figure}[h!]\centering
\includegraphics[width=100mm]{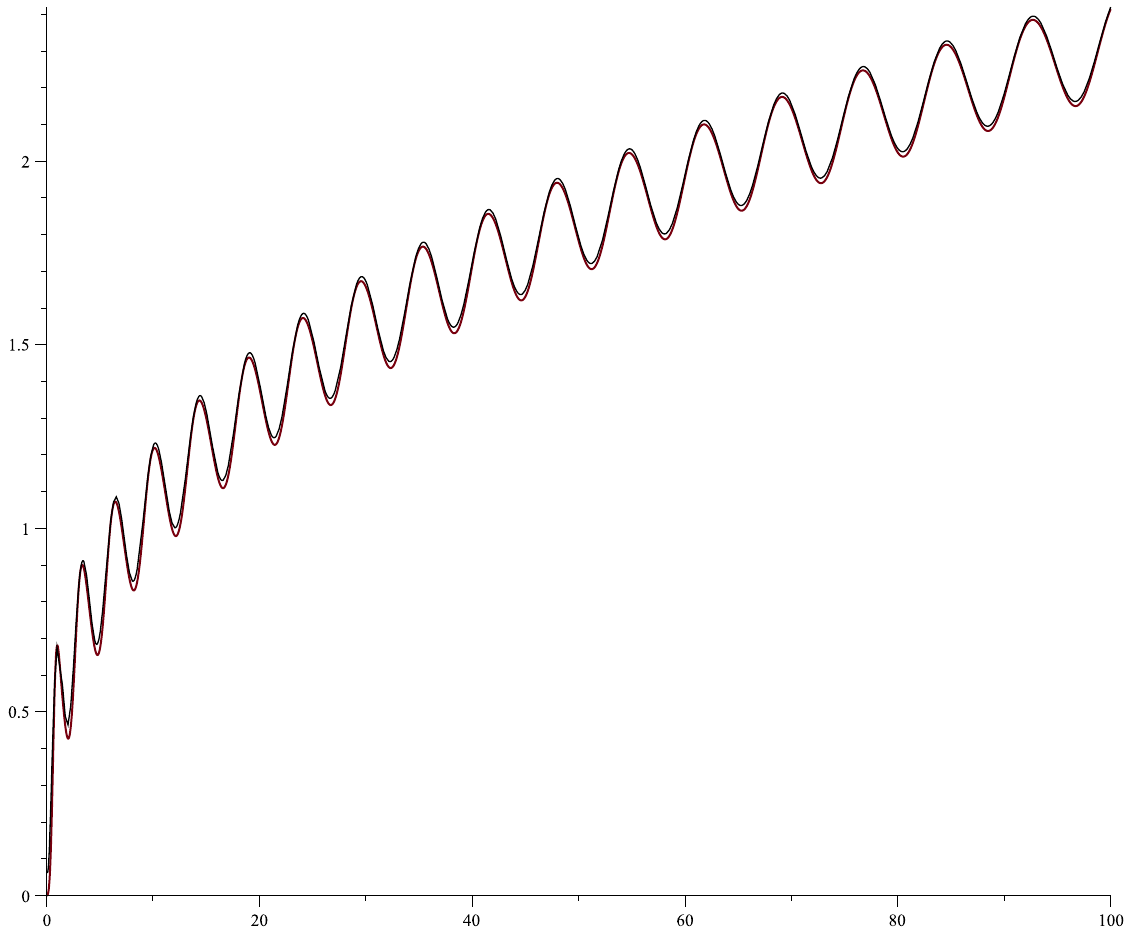}
\caption{The plots of $\operatorname{Re}u(\tau)$ and $\operatorname{Re}u_{\rm as}(\tau)$ for $\alpha=5/7$ and $b=1$.
The corresponding extrema are higher for $\operatorname{Re}u(\tau)$.}\label{fig:Re(u)alpha5d7b1}
\end{figure}

\begin{figure}[h!]\centering
\includegraphics[width=100mm]{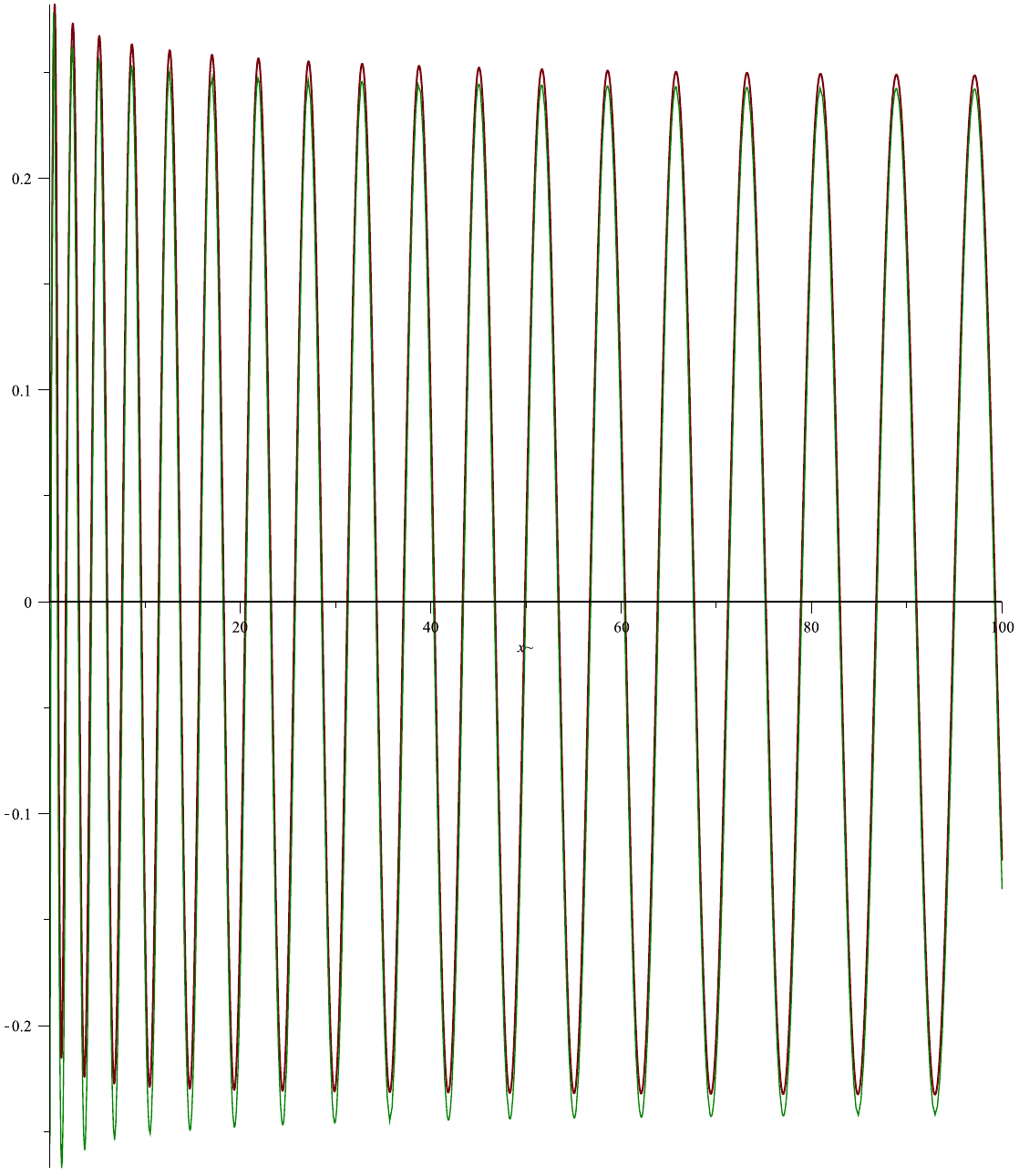}
\caption{The plots of $\operatorname{Im}u(\tau)$ and $\operatorname{Im}u_{\rm as}(\tau)$ for $\alpha=5/7$ and $b=1$.
The corresponding extrema are higher for $\operatorname{Im}u(\tau)$.}\label{fig:Im(u)alpha5d7b1}
\end{figure}

Below we present a comparison of the solution with its asymptotic approximation, with the help of the
above asymptotic formulae and Maple 16. As it often happens for the Painlev\'e equations the leading term of
asymptotics for large values of independent variable gives quite good numerical approximation of the solution
when the independent variable takes relatively small values.
For illustrative purposes I especially show the plots
(see Figs.~\ref{fig:Re(u)alpha2d7b1d80} and~\ref{fig:Im(u)alpha2d7b1d80}) for a small value of
the parameter $b$, at such values the difference between the solution and its asymptotics can be easily observed.
The increase of $b$ means a rescaling of the argument and we, in some sense, observe how the same solution would
behave for the larger values of the argument although we keep the same interval for our plots
(see Figs.~\ref{fig:Re(u)alpha5d7b1} and~\ref{fig:Im(u)alpha5d7b1}). Of course, when $\alpha$ approaches the
boundaries of the validity of asymptotics~\eqref{eq:asymp-main-1} the numerical correspondence with the solution
for the small and finite values of $\tau$ becomes worse, however by increasing~$b$, one can observe that
asymptotics is working on the whole interval~$(0,1)$.

Now consider the Suleimanov solution, i.e., the odd meromorphic solution for $\alpha=1/2$. In this case, the formulae presented above can be considerably simplified
\begin{gather*}
\operatorname{Im}u_{\rm as}(\tau) =-\sqrt{\frac{b\sqrt{3}\ln2}{4\pi}} \cos{\phi},\\
\operatorname{Re}u_{\rm as}(\tau) =\frac{b^{2/3}\tau^{1/3}}2-\sqrt{\frac{b\ln2}{4\pi\sqrt{3}}} \sin{\phi},
\end{gather*}
where
\begin{gather*}
\phi =3\sqrt{3} b^{1/3}\tau^{2/3}-\frac{\ln2}{2\pi}\ln\big(3\sqrt{3} b^{1/3}\tau^{2/3}\big)+\phi_0,\\
\phi_0 =\frac34\pi-\frac{\ln2\ln{12}}{2\pi}-\arg\Gamma\left(-i\frac{\ln2}{2\pi}\right).
\end{gather*}
Qualitatively the plots of the Suleimanov solution and its asymptotics look similar to the ones
presented on Figs.~\ref{fig:Re(u)alpha2d7b1d80}--\ref{fig:Im(u)alpha5d7b1}, with the same comment
concerning dependence on the parameter~$b$.
\begin{Remark}
The asymptotics of Suleimanov's solution for $\alpha=-1/2$ has the same real part as for the case $\alpha=1/2$
presented above, whilst its imaginary part differs by sign, namely,
$\operatorname{Im}u_{\rm as}(\tau)=\sqrt{\frac{b\sqrt{3}\ln2}{4\pi}} \cos{\phi}$.
\end{Remark}

\subsection[$a<0$]{$\boldsymbol{a<0}$}
\begin{figure}[h!]\centering
\includegraphics[width=100mm]{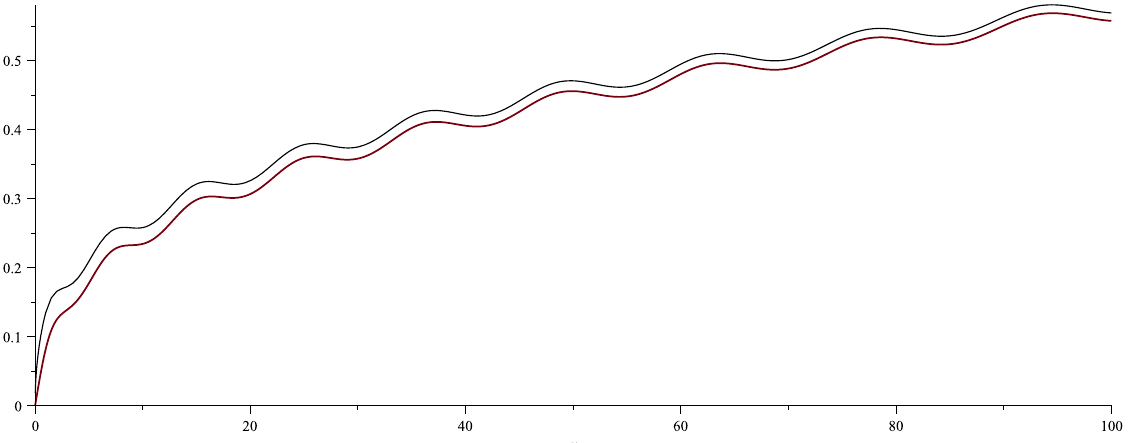}
\caption{The plots of $u(\tau)$ and $u_{\rm as}(\tau)$ for $a=-2/3$ and $b=1/8$. The second plot is higher.}\label{fig:u-a2d3b1d8}
\end{figure}

\begin{figure}[h!]\centering
\includegraphics[width=100mm]{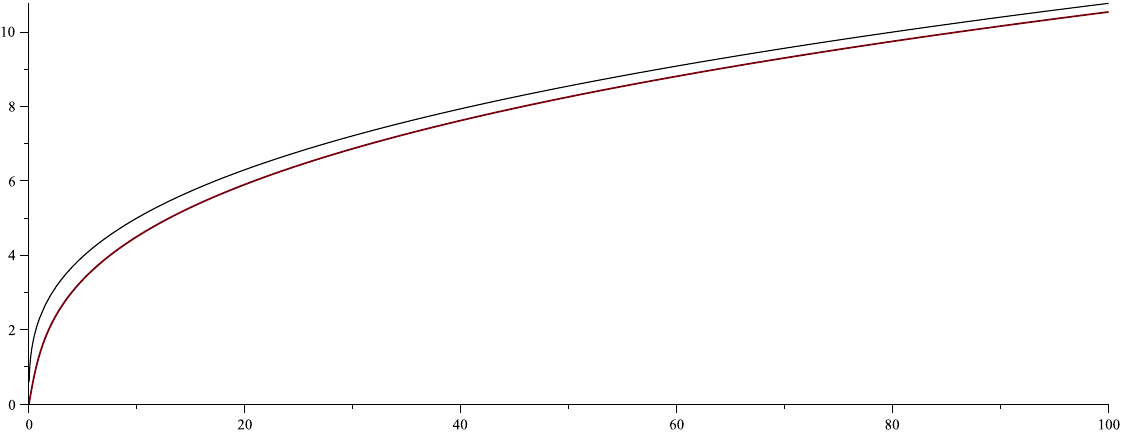}
\caption{The plots of $u(\tau)$ and $u_{\rm as}(\tau)$ for $a=-3$ and $b=10$. The second plot is higher.}\label{fig:u-a3b10}
\end{figure}

Assume $a\in\mathbb R$. Both conditions, \eqref{ineq:cond2a} and the weaker one discussed in
Remark~\ref{rem:weaker-condition}, give the same result $a<0$, which is assumed below. In this
case the solution $u(\tau)$ is real. After some calculation one proves that
asymptotics~\eqref{eq:asymp-main-2}--\eqref{eq:asymp-main-4} is also, as it should be, real and
can be rewritten as follows
\begin{gather}
u_{\rm as}=\frac{b^{2/3}}{2}\tau^{1/3}\nonumber\\
\hphantom{u_{\rm as}=}{} -\frac{b^{1/2}}{3^{1/4}}\sqrt{\frac{-\ln(1-{\rm e}^{2\pi a})}{2\pi}}
\cos\left(3^{3/2}b^{1/3}\tau^{2/3}-\frac{\ln\big(1-{\rm e}^{2\pi a}\big)}{2\pi}\ln\big(3^{3/2}b^{1/3}\tau^{2/3}\big)
+\phi_0\right),\nonumber\\
\phi_0=a\ln\big(2+\sqrt{3}\big)-\ln(12)\frac{\ln\big(1-{\rm e}^{2\pi a}\big)}{2\pi}-\frac{\pi}{4}
-\arg\left(\Gamma\left(-i\frac{\ln(1-{\rm e}^{2\pi a})}{2\pi}\right)\right).\label{eq:uas-a<0-1}
\end{gather}
On Figs.~\ref{fig:u-a2d3b1d8} and \ref{fig:u-a3b10} we present examples of asymptotics~\eqref{eq:uas-a<0-1}.
The smaller values of $a$ the better approximation of $u(\tau)$ by $u_{\rm as}(\tau)$, it can be seen comparing
the scale of the $y$-axes and also taking into account that, as mentioned above, by enlarging $b$ we in a
sense enlarging $\tau$ therefore approximation of $u(\tau)$ by $u_{\rm as}(\tau)$ becomes better. In making
numerics this, however, means that if we consider plots on the same $\tau$-segment, then the larger $b$
the higher accuracy of calculations are required, because, it is equivalent to consideration of the solution
on a longer interval. We also see that for the large negative values of $a$, say, $a=-3$ is already ``large'', the amplitude of oscillation
becomes so small that the plot visually looks like the cubic parabola.
We have chosen the parameters $a$ and $b$ such that on one hand one can see the difference between
$u(\tau)$ and $u_{\rm as}(\tau)$ and on the other hand the fact the convergence of the asymptotics to solution.

\section[Positiveness of $\operatorname{Re}u(\tau)$]{Positiveness of $\boldsymbol{\operatorname{Re}u(\tau)}$}\label{sec:positive}
Looking on Figs.~\ref{fig:u-a2d3b1d8} and \ref{fig:u-a3b10} one arrives at the following
\begin{Proposition}\label{prop:real-positive}
The function $u(\tau)$ with $a<0$ is bounded for $\tau\in\mathbb{R}$. Moreover, $u(\tau)$ is an odd
function and $u(\tau)>0$ for $\tau>0$.
\end{Proposition}
As follows from expansion~\eqref{eq:zero-expansion} the solution possesses the following symmetries:
\begin{gather}\label{eq:u-symmetries}
u(\tau)=-u(-\tau),\qquad u(\tau,a)=-iu(i\tau,-a).
\end{gather}
\begin{proof}
It is enough to prove boundedness of our solution for $\tau>0$. Expansion~\eqref{eq:zero-expansion} implies (recall $a<0$) $u(\tau)>0$ for $0<\tau<\delta$ at least for some
rather small $\delta>0$. Assume that $\tau_z>0$ is a zero of~$u(\tau)$, then substituting the Taylor expansion
of $u(\tau)$ at $\tau_z$ with the leading term $A(\tau-\tau_z)^n$ with \mbox{$n>0$} into equation~\eqref{eq:dp3} one
proves that $n=1$ and $b^2/A^2=-1$. Therefore, real solutions of equation~\eqref{eq:dp3} have no zeroes
on the real line. It means that our solution is positive on the positive semiaxis. Analogously, substituting
the Laurent expansion of $u(\tau)$ at pole $\tau_p$ with the leading term $A_1/(\tau-\tau_p)^m$ with $m>0$
one finds $m=2$ and $2\tau_p=-8A_1$, If we assume $\tau_p>0$ we find that
\mbox{$u(\tau)\underset{\tau\to\infty}{\sim}-\tfrac{\tau_p}{4(\tau-\tau_p)^2}<0$}. Thus, our positive solution
does not have poles on positive semiaxis.
\end{proof}

\begin{Lemma}If a solution of equation~\eqref{eq:dp3} has positive real part for $\tau>0$, then it does not have poles for $\tau>0$.
\end{Lemma}
\begin{proof}The Laurent expansion of solutions of equation~\eqref{eq:dp3} at the pole $\tau_p$ reads
\begin{gather*}
-\frac{\tau_p}{4(\tau-\tau_p)^2}+O(1),\qquad\tau\to\tau_p.\tag*{\qed}
\end{gather*}\renewcommand{\qed}{}
\end{proof}

Since $\tau$ and $\tau_p$ are assumed real, then the imaginary part of any solution does not have
poles on the real axis except possibly $\tau_p=0$. Moreover, if we assume that $\tau_p>0$, then real
part of any solution with a pole is negative in some its neighborhood, however we assume that the real part is positive.
\begin{Proposition}\label{prop:positiveness-epsilon}For $\operatorname{Re}a<0$ there exists $\varepsilon>0$ such that for $|\operatorname{Im} a|<\varepsilon$ and $\tau>0$
\mbox{$\operatorname{Re}u(\tau)>0$}.
\end{Proposition}
Assume there is a sequence $a_n$ with $\operatorname{Re}a_n\equiv\beta<0$ and $\operatorname{Im} a_n\to0$, such that
the functions $u(\tau)=u(\tau,a_n)$ has a sequence of zeroes $\tau_n$. If the sequence $\tau_n$ has a bounded
subsequence, then it leads to immediate contradiction with Proposition~\ref{prop:real-positive}, since in this case
$u(\tau_0)=u(\tau_0,0)=0$ where $\tau_0$ is a finite limiting point of the sequence $\tau_n$.

The case when the sequence $\tau_n$ is unbounded also cannot happen. It follows from the general property of local
uniformness of asymptotics of the Painlev\'e equations. This property means the following: assume that an asymptotic
formula is valid for some open subset in the space of the monodromy parameters, as it happens for
asymptotics~\eqref{eq:asymp-main-1}--\eqref{eq:asymp-main-4}, which is valid in the open strip, $a=\beta+i\alpha$
with $\beta<0$ and $|\alpha|<1$. We consider any compact subset of this strip, for our purposes
it is enough to fix $\beta=\beta_0<0$ and a segment $|\alpha|\leq1-\epsilon$ for some positive $\epsilon<1$.
Then there exist $\Upsilon=\Upsilon(\beta,\epsilon)$ such that the error estimate in equation~\eqref{eq:asymp-main-1}
holds for all $\tau>\Upsilon$ and all $a$ from the above compact subset. Thus, if we choose the error estimate in~\eqref{eq:asymp-main-1} small enough by fixing~$\Upsilon$, then the zeroes $\tau_n>\Upsilon$ cannot exist because~$u(\tau)$ is close to the first term in asymptotics~\eqref{eq:asymp-main-1} and, therefore, cannot vanish for any
values of $\tau>\Upsilon$. The following conjecture is the extension of Proposition~\ref{prop:positiveness-epsilon}.
\begin{Conjecture}\label{con:positiveness}
If $\operatorname{Re}a\leq0$ and $0<|\operatorname{Im}a|<1$, then $\operatorname{Re}u(\tau)>0$ for $\tau>0$.
\end{Conjecture}
\begin{Remark}
Solutions for the parameter $a$ from Conjecture~\ref{con:positiveness} have neither zeroes nor poles on
the positive semiaxis. To prove the absence of poles one can surely apply Zhou's vanishing lemma~\cite{Z}.
In our case the Riemann--Hilbert problem corresponding to the monodromy data defined in Section~\ref{sec:monodromy}
can be formulated on the positive real semiaxis and the circle centered at~$0$. The jump matrices are
\begin{gather*}
\begin{pmatrix}1&s_0\\s_0&1+s_0^2\end{pmatrix},\qquad
{\rm e}^{2\pi(a-i/2)\sigma_3},\qquad {\rm and}\qquad G
\end{gather*}
inside, outside, and on the circle, respectively. There are singular points at
zero and infinity with the proper behavior, which we do not discuss here.
The vanishing lemma implies that the jump matrices on the real axis are positively defined;
we get the following condition for the trace of the diagonal matrix
\begin{gather*}
\cosh(2\pi(a-i/2))+\cosh(2\pi(\bar{a}+i/2))>0 \quad\Rightarrow\quad
\cos(2\pi\operatorname{Im}a)<0 \quad \Rightarrow \quad \frac14<\operatorname{Im}a<\frac34,
\end{gather*}
and arbitrary $\operatorname{Re}a$. The trace of the nondiagonal matrix is $2+s_0^2$ its positiveness
exactly coincide with the above condition on $a$. There is one more condition that demand
the vanishing lemma, $G^\dagger=G^{-1}$. I recall that matrix $G$ (see the second line of equations in
system~\eqref{sys:monodromy-data-u-varphi-unique}), contains one free complex parameter, which is
related with the constant of integration in equation~\eqref{eq:varphi}, so that for a given $u(\tau)$
it can be chosen arbitrarily. If we choose it to satisfy the condition on matrix $G$ we get
$\operatorname{Im}a=\pm 1/2$ with arbitrary $\operatorname{Re} a$.

Thus, in particular, the Suleimanov solution is regular on the positive semiaxis,
in fact, on the coordinate cross, because of the symmetries~\eqref{eq:u-symmetries}. Unfortunately, this
remark does not shed any light on the proof of Conjecture~\ref{con:positiveness}.
\end{Remark}

\subsection*{Acknowledgements}
The author is grateful to P.D.~Miller and B.I.~Suleimanov for discussions of the papers~\cite{BLM, S}.
The author is indebted to the referees for their significant contribution to improving the quality of
the original version of this paper.

\pdfbookmark[1]{References}{ref}
\LastPageEnding


\begin{thebibliography}{99}
\addcontentsline{toc}{section}{References}
\footnotesize\itemsep=0pt

\bibitem{BLM}
Bilman D., Ling L., Miller P.D., Extreme superposition: roague waves of
 infinite order and the {P}ainlev\'e-III hierarchy, \href{https://arxiv.org/abs/1806.00545}{arXiv:1806.00545}.

\bibitem{BE}
Erd\'{e}lyi A., Magnus W., Oberhettinger F., Tricomi F.G., Higher
 transcendental functions, {V}ol.~{I}, McGraw-Hill Book Company, Inc., New
 York~-- Toronto~-- London, 1953.

\bibitem{GIL2013}
Gamayun O., Iorgov N., Lisovyy O., How instanton combinatorics solves
 {P}ainlev\'{e} {VI}, {V} and {III}s, \href{https://doi.org/10.1088/1751-8113/46/33/335203}{\textit{J.~Phys.~A: Math. Theor.}}
 \textbf{46} (2013), 335203, 29~pages, \href{https://arxiv.org/abs/1302.1832}{arXiv:1302.1832}.

\bibitem{G}
Garnier R., Sur des \'{e}quations diff\'{e}rentielles du troisi\`eme ordre dont
 l'int\'{e}grale g\'{e}n\'{e}rale est uniforme et sur une classe
 d'\'{e}quations nouvelles d'ordre sup\'{e}rieur dont l'int\'{e}grale
 g\'{e}n\'{e}rale a ses points critiques fixes, \href{https://doi.org/10.24033/asens.644}{\textit{Ann. Sci. \'{E}cole
 Norm. Sup.~(3)}} \textbf{29} (1912), 1--126.

\bibitem{HL}
Hardy G.H., Littlewood J.E., Tauberian theorems concerning power series and
 {D}irichlet's series whose coefficients are positive, \href{https://doi.org/10.1112/plms/s2-13.1.174}{\textit{Proc. London
 Math. Soc.}} \textbf{13} (1914), 174--191.

\bibitem{HW}
Hardy G.H., Wright E.M., An introduction to the theory of numbers, 5th ed., The
 Clarendon Press, Oxford University Press, New York, 1979.

\bibitem{KV2004}
Kitaev A.V., Vartanian A., Connection formulae for asymptotics of solutions of
 the degenerate third {P}ainlev\'{e} equation.~{I}, \href{https://doi.org/10.1088/0266-5611/20/4/010}{\textit{Inverse Problems}}
 \textbf{20} (2004), 1165--1206, \href{https://arxiv.org/abs/math.CA/0312075}{arXiv:math.CA/0312075}.

\bibitem{KV2010}
Kitaev A.V., Vartanian A., Connection formulae for asymptotics of solutions of
 the degenerate third {P}ainlev\'{e} equation:~{II}, \href{https://doi.org/10.1088/0266-5611/26/10/105010}{\textit{Inverse Problems}}
 \textbf{26} (2010), 105010, 58~pages, \href{https://arxiv.org/abs/1005.2677}{arXiv:1005.2677}.

\bibitem{KVdP3int}
Kitaev A.V., Vartanian A., Asymptotics of integrals of some functions related
 to the degenerate third {P}ainlev\'e equation, \href{https://arxiv.org/abs/1811.05276}{arXiv:1811.05276}.

\bibitem{OEIS}
Sloane N.J.A., {S}equences {A}001764, {A}023745, {A}029858, {A}031988, and
 {A}014915, {T}he on-line encyclopedia of integer sequences,
 \url{http://oeis.org}.

\bibitem{S}
Suleimanov B.I., Effect of a small dispersion on self-focusing in a spatially
 one-dimensional case, \href{https://doi.org/10.1134/S0021364017180126}{\textit{JETP Lett.}} \textbf{106} (2017), 400--405,
 \href{https://arxiv.org/abs/1706.06849}{arXiv:1706.06849}.

\bibitem{W}
Weisstein E.W., Dirichlet divisor problem, {W}olfram MathWorld,
 \url{http://mathworld.wolfram.com/DirichletDivisorProblem.html}.

\bibitem{Z}
Zhou X., The {R}iemann--{H}ilbert problem and inverse scattering,
 \href{https://doi.org/10.1137/0520065}{\textit{SIAM~J. Math. Anal.}} \textbf{20} (1989), 966--986.

\end{thebibliography}
\end{document}